\documentclass[11pt]{article}
\usepackage{enumerate}
\usepackage[T1]{fontenc}
\usepackage{amsmath,amssymb}
\usepackage{natbib}
\usepackage[usenames, dvipsnames]{color}
\usepackage[dvipsnames]{xcolor}
\usepackage{dsfont}
\usepackage{smile}      
\usepackage{multirow}

\usepackage{lmodern}
  
\usepackage[linesnumbered,ruled]{algorithm2e}
\usepackage{bm}
\usepackage{amsmath}

\RequirePackage[colorlinks,
            linkcolor=green,
            anchorcolor=blue,
            citecolor=blue
            ]{hyperref}
            
\usepackage{thmtools,thm-restate} 
\usepackage{amssymb}

\def\biggiven{\,\big{|}\,}

\def\tr{\mathop{\text{tr}}\kern.2ex}

\def\E{{\mathbb E}}

\def\supp{\mathop{\text{supp}}}

\long\def\comment#1{}

\def\tr{\mathop{\text{Tr}}}

\newcommand{\bel}{\begin{eqnarray}\label}
\newcommand{\eel}{\end{eqnarray}}
\newcommand{\bes}{\begin{eqnarray*}}
\newcommand{\ees}{\end{eqnarray*}}

\let\emptyset\varnothing
\def\mid{\,|\,}

\def\E{{\mathbb E}}

\def\supp{\mathop{\text{supp}\kern.2ex}}
\def\argmin{\mathop{\text{\rm arg\,min}}}

\def\tr{{\rm{Tr}}}

\def\supp{\mathop{\text{supp}}}

\def\tr{\mathrm{Tr}}

\usepackage{mathrsfs}

\usepackage{bbm}

\usepackage{fullpage}

\usepackage{algpseudocode}

\usepackage{marginnote}

\def \tanh {\text{tanh}}

\def\##1\#{\begin{align}#1\end{align}}
\def\$#1\${\begin{align*}#1\end{align*}}

\usepackage[shortlabels]{enumitem}

\usepackage{undertilde}
\theoremstyle{plain}

\theoremstyle{mytheoremstyle}

\usepackage{bookmark}
\bookmarksetup{startatroot}

\definecolor{darkblue}{rgb}{0,0,0.6}

\newcommand{\RN}[1]{%
  \textup{\uppercase\expandafter{\romannumeral#1}}%
}

\begin{document}

\title{\huge The Wreaths of KHAN: Uniform Graph Feature Selection with False Discovery Rate Control}

\author{
    Jiajun Liang$^{1, *}$ \ 
    Yue Liu$^{2, *}$ \ 
    Doudou Zhou$^{3}$ \ 
    Sinian Zhang$^{4}$  \ 
    Junwei Lu$^{3, \dagger}$
}
 
\footnotetext[1]{Department of Statistics, Purdue University, West Lafayette, IN 47906}
\footnotetext[2]{Department of Statistics, Harvard University, Cambridge, MA 02138}
\footnotetext[3]{Department of Biostatistics, Harvard Chan School of Public Health, Boston, MA 02130}
\footnotetext[4]{School of Statistics, Renmin University of China, Beijing, China 100872}
\renewcommand{\thefootnote}{\fnsymbol{footnote}}
\footnotetext{* Equal contribution}
\footnotetext{$\dagger$ Corresponding author (junweilu@hsph.harvard.edu) }

\date{ }
\maketitle

\begin{abstract}

Graphical models find numerous applications in biology, chemistry, sociology, neuroscience, etc.
While substantial progress has been made in graph estimation, it remains largely unexplored how to select significant graph signals with uncertainty assessment, especially those graph features related to topological structures including cycles (i.e., wreaths), cliques, hubs, etc. These features play a vital role in protein substructure analysis, drug molecular design,  and brain network connectivity analysis. 
To fill the gap, we  propose 
a novel inferential framework for general high dimensional graphical models to select graph features with false discovery rate controlled.
Our method is based on the maximum of $p$-values from single edges that comprise the topological feature of interest, thus is able to detect weak signals. 
Moreover, we introduce the $K$-dimensional persistent Homology Adaptive selectioN (KHAN) algorithm to  select all the  homological features  within $K$ dimensions with the uniform control of the false discovery rate over continuous filtration levels. The KHAN method applies a novel discrete Gram-Schmidt algorithm to select statistically significant generators from the homology group. 
 We apply the structural screening method to identify the important residues of the SARS-CoV-2 spike protein during the binding process to the ACE2 receptors. We score the residues for all domains in the spike protein by the $p$-value weighted filtration level in the network persistent homology for the closed, partially open, and open states and identify the residues crucial for protein conformational changes and thus being potential targets for inhibition.
\end{abstract}


\section{Introduction}
Graphical models are a prevalent tool for modeling relationships between variables. They are widely used in various fields, especially certain sub-structures of the  graphs between  variables is of great interest. In biology, proteins are vital to life and their practical function is highly influenced by their unique sub-structure, which is determined by the combination of amino acids. 
For example, the Y-shape of an antibody enables it to bind to antigens like bacteria and viruses with one end and to other immune-system proteins with the other end \citep{janeway2001structure}. Similarly, the ``wreath'' shape of the DNA polymerase III holoenzyme facilitates the efficient formation of a ring around DNA, enabling fast DNA synthesis \citep{kelman1995dna, 2013v}. 
Alphabet/Google is developing AlphaFold to predict the structures of proteins with a reinforcement learning approach \citep{jumper2021highly, evans2021protein}.   
In chemistry, a critical task is modeling molecular structures, which determine their pharmacological and ADME/T properties (absorption, distribution, metabolism, excretion, and toxicity). For instance, the presence of the cyanide group makes hydrogen cyanide (HCN) highly toxic, leading to potential death within minutes \citep{national2002acute}.  
High-Density Polyethylene (HDPE) has a linear structure, offering density and strength for products requiring rigidity, such as milk jugs and water pipes. Conversely, Low-Density Polyethylene (LDPE), with its branched structure, is suited for flexible applications, including plastic bags and squeeze bottles \citep{fried2014polymer}. In neuroscience, cliques within the brain network are vital as they are neurologically strongly associated with efficient support of behavior \citep{bassett2017network, wang2021learning}.
In sociology, social networks play an important role in understanding behavior among individuals or organizations and making social recommendations or marketing strategies. For instance, the presence of a ``hub'' structure in the social network indicates that the central user is active and influential, as they are connected to many other users. Therefore, the advertising industry often targets these individuals with free samples to increase their reputation, rather than selecting people at random \citep{ ilyas2011distributed, li2018influence,lee2019discovering}.
 In this paper, we mainly focus on the application of the SARS-CoV-2 spike protein sub-structure analysis. To identify the  potential targets for inhibition in the vaccine design, we aim to select the protein sub-structures persistently vary in the binding process to the ACE2 receptors \citep{ou2020characterization}.

Despite the vast applications of the sub-structure detection of graphical models mentioned above, the majority of existing research focused on the edge-wise inference \citep{cai2011constrained, fan2016overview, cai2016estimating, ding2020estimation}. The majority literature of the statistical methods on the false discovery control  \citep{benjamini2010discovering, van2014asymptotically, li2021ggm}  worked on selecting the parametric signals which cannot be directly applied to the discrete graph sub-structures selection.
This paper aims to bridge the gap by proposing the selecting the graph topological features with the false discovery rate control.
We also propose to screen  topological features persistent under continuous filtration levels utilizing the framework of the persistent homology \citep{horak2009persistent, aktas2019persistence}.  
In specific, let $G^* = (V, E^*)$ represents the true undirected graph, where $V=\{1,2 \cdots, d\}$ is the collection of nodes and $E^*$ is the edge set of the graph.
Suppose $\cF=\{F_1,F_2,\ldots,F_J\}$ is the set consisting of all  graph features of interest in $G^*$ that we aim to select, where $J$ is the cardinality of candidate  graph features. Examples of such  graph features include cliques, which are fully connected subsets of nodes, and loops, which are closed paths within the graph.  We consider the multiple hypotheses testing whether each  graph feature  $F_j\in {\cal F}$ can be embedded into the true graph $G^*$, i.e.,
\begin{equation}\label{eqn:hypo}
H_{0j}: F_j \not\subseteq  G^* \quad \mbox{v.s.}\quad H_{1j}: F_j \subseteq G^*,\quad 1\leq j\leq J,
\end{equation}   
where $F \subseteq G^*$ means $F$ is a subgraph of $G^*$. 
In order to test these hypotheses 
with false discovery rate controlled,  we consider the maximum $p$-value among all edges within $F_j$. This approach differs significantly from previous methods in combinatorial inference \citep{shen2023combinatorial, liu2023lagrangian, zhang2021startrek}, which adjust the multiplicity by the maximum statistic and maximal Gaussian multiplier bootstrap. In comparison, 
our proposed method only requires single edge $p$-values and therefore can detect weaker signals and is also computationally more efficient. 

\begin{figure}[h!]
	\begin{center}
\includegraphics[width=0.5\textwidth,angle=0]{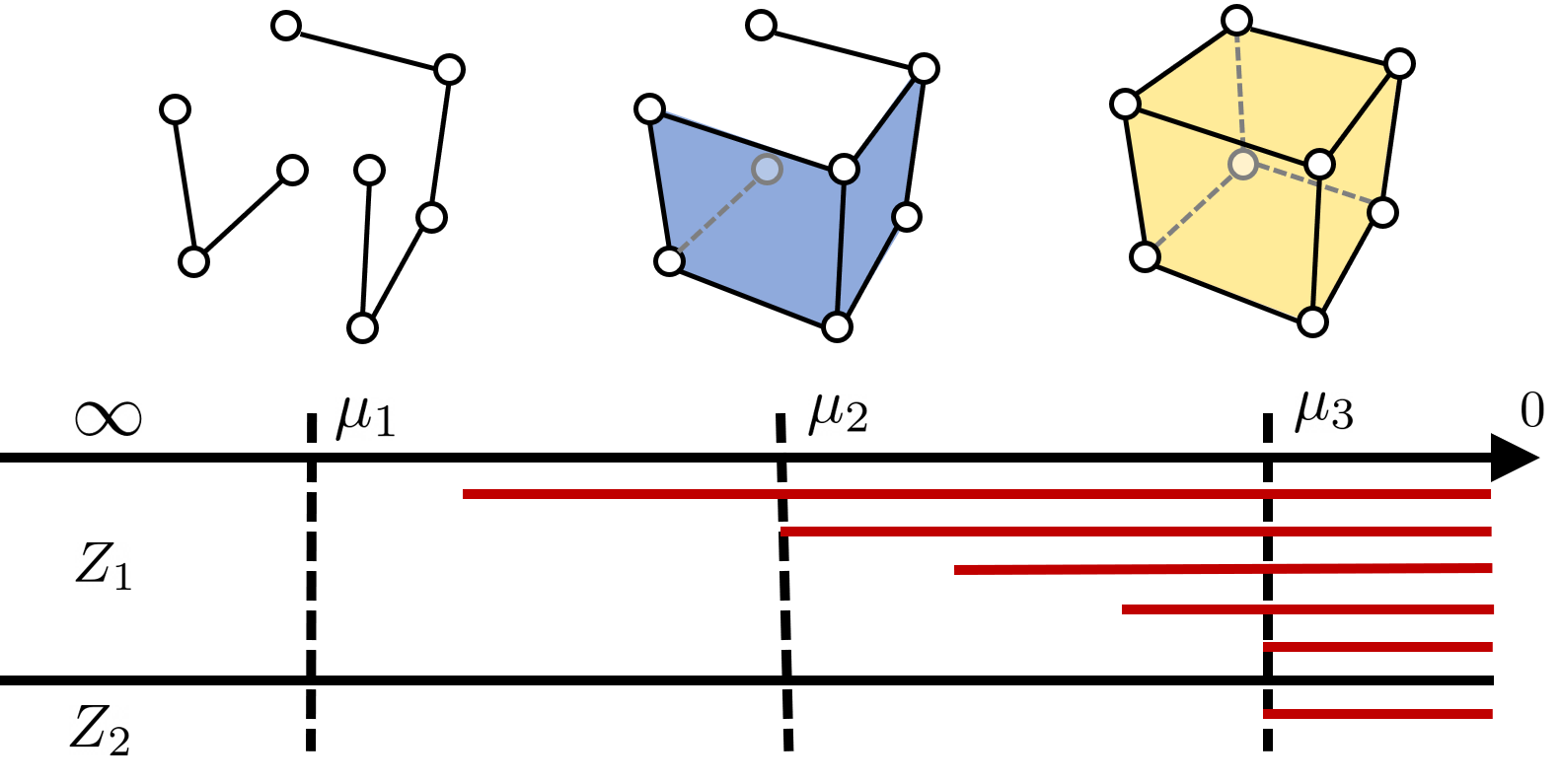} 
	\end{center}
\vspace{-10pt}
	\caption{
 Illustration of the filtered graph, homological features, and the persistent barcode. Two $1$-dimensional homological features (blue, i.e., $\mathrm{rank}(Z_1(E(\mu_2))) = 2$) show up at $\mu_2$ and a $2$-dimensional homological feature (yellow, i.e., $\mathrm{rank}(Z_2(E(\mu_3))) = 1$) appears at $\mu_3$.}
	\label{fig:filtration}
\end{figure}

Besides selecting the given graph features, we also consider to evaluate the strength of these signals at multiple scales.  The network persistent homology \citep{horak2009persistent, lee2012persistent, aktas2019persistence} and the persistent barcode \citep{ghrist2008barcodes, kovacev2016using} are powerful tools for analyzing and measuring the persistence of multiple homological features, such as triangles, tetrahedrons,  and higher-dimensional polytopes. However, few methods have been developed to quantify the uncertainty when inferring the network persistent homology from noisy datasets.
We solve this problem  by considering the multiple hypotheses at the filtration level $\mu$ as 
\begin{equation}\label{eqn:3}
\begin{aligned}
& H_{0j}(\mu):\text{the $j$-th homological feature} \nsubseteq G^*(\mu),\\
     \text{v.s.}  \ 
     & 
     H_{1j}(\mu):\text{the $j$-th homological feature} \subseteq G^*(\mu),
\end{aligned}
\end{equation}
where  $G^*(\mu)$ is the filtered graph at level $\mu$ and  the  homological feature  is a group  generator of the graph homology group which we will define in details in Section~\ref{sec:method}. Intuitively, a $K$-dimensional homological feature is a $K$-dimensional ``wreaths'' on the graph as a discrete analogue of a $K$-dimensional sphere in Euclidean space. Figure~\ref{fig:filtration} illustrates the $1$-dimensional and $2$-dimensional homological features and how they appear with the change of the filtration level.

Our goal is to screen all the homological features under the given dimension $K$, with detailed definition in Section~\ref{sec:setup}. We could not directly apply the approach for testing graph features in \eqref{eqn:hypo}, as the candidate group generators could be linearly dependent.  To tackle the double dependency: the probabilistic dependency among testing statistics and linear dependency among the generators in the homology group, we propose a novel discrete Gram-Schmidt $p$-value screening method  to select the linearly independent homological features at any given filtration level.  The algorithm is a homological analogue of the Gram-Schmidt algorithm which iteratively selects the generator with the smallest $p$-value among the subgroup which is linearly independent of the currently selected generators. The other challenge is to guarantee the false discovery rate of the selected homological features is uniformly controlled over continuous filtration levels. Instead of conducting the screening for each filtration level $\mu$, our method adaptively identifies the change points  in the persistent barcode in Figure~\ref{fig:filtration} which is both computationally and statistically efficient. This results in the   $K$-dimensional persistent Homology Adaptive selectioN (KHAN) algorithm. See Section~\ref{sec:barcode} for the details. We introduce the new concept of uniform False Discovery Rate (uFDR) to quantify the performance of persistent homology over the continuous filtration levels and prove that our algorithm can handle the double dependency and control uFDR at the given level.

\subsection{Related Work}

Learning  graph structures  has been widely studied in the existing literature. There have been many significant progress on recovering the underlying graph by estimation.
For the Gaussian graphical model, many works estimate the
graph through estimating the inverse covariance matrix $\Theta^*$ where $\Theta^*_{uv}\neq 0$  if and only if the edge $(u, v) \in E^*$
\citep{meinshausen2006high, friedman2008sparse, lam2009sparsistency, peng2009partial, ravikumar2011high, cai2011constrained, shen2012likelihood}.
For the Ising model, which consists of discrete variables ($+1$ or $-1$), the estimation of parameters is achieved by using penalized logistic regression \citep{ravikumar2010high}, and  the study of parameter estimation and inference is furthered through the maximum pseudo-likelihood estimate \citep{bhattacharya2018inference}.


Persistent homology is a powerful tool to analyze topological features over these graphs, which involves  tracking the birth and death of topological features along filtrations of graphs \citep{horak2009persistent, aktas2019persistence}. However, the primary focus of these studies has been on extracting the persistent topological features of a deterministic networks without the uncertainty quantification of the graph estimation error.  An inferential method is proposed \citep{fasy2014confidence} to construct
confidence set for the persistence diagrams for the landscape of probability density functions. However, this method is restricted to the persistent homology of Euclidean spaces,  and cannot be applied to discrete networks.

For the inferential methods on graphical models, most literature has focused on continuous quantities and local properties of graph (e.g., presence of edges) \citep{cai2013optimal,jankova2015confidence,ren2015asymptotic, cai2016inference,  chen2016asymptotically, neykov2018unified}. For instance, hypothesis tests for the covariance matrix have been constructed \citep{cai2013optimal}, tests for the presence of edges in the Gaussian graphical model have been developed \citep{ren2015asymptotic}, and considerations for the estimation and testing of high-dimensional differential correlation matrices have been proposed \citep{cai2016inference}.
Regarding multiple tests in graphical models, efforts have been made to select edges with false discovery rate (FDR) control using the Benjamini-Hochberg method \citep{benjamini1995controlling}, and methodologies have been developed to learn the structure of Gaussian graphical models with FDR control based on the knockoff framework \citep{li2021ggm}. However, a critical gap in this corpus of work is its focus on methodologies primarily for continuous parameters, and it is challenging to apply those techniques to select discrete graph features addressed in this paper.

Research on the combinatorial structure and global properties of graphs has led to proposals for testing monotone properties within the Gaussian graphical model, such as connectivity, the maximum degree, and cycle presence \citep{lu2017adaptive}. Similar efforts have been directed towards testing graph properties in the Ising model  \citep{neykov2019property, jin2020computational}. However, these methodologies are limited as they are designed for conducting single tests without controlling the FDR. Addressing the challenge of simultaneous multiple testing, the StarTrek algorithm for selecting hub nodes based on a pre-specified degree threshold in graphs has been proposed \citep{zhang2021startrek}. However, the StarTrek procedure \citep{zhang2021startrek} is primarily designed for selecting hub nodes, whereas our goal is to select general graph features of interest including the hub nodes. In comparison,  is much more complicated and challenging, with broader applications.
In methodology, the StarTrek procedure involves the use of the maximal statistic for uniform control over edges and the Gaussian multiplier bootstrap for quantile estimation, which are statistically conservative and computationally intensive. In contrast, our proposed method does not require the use of either the maximal statistic or the Gaussian multiplier bootstrap, making it simpler and more efficient.
Moreover, our method can be applied to study the persistence of graph features within the framework of persistent homology and offers stronger false discovery control across multiple  scales.

\subsection{Our Contributions}

In methodology,  we propose a novel inferential framework to select general graph structural features with false discovery rate control. One of the challenge problem is to assign each structural feature with a $p$-value.  Compared to the estimation of the $p$-values of edges, the estimation of the $p$-value of a graph structure requires combinatorial and uniform control on involved edges for multiple cases. As a result, the estimation of the $p$-values of edges within and between structures could be highly dependent.  To address this issue, maximal statistics have been used in the literature \citep{zhang2021startrek}. However, computing $p$-values for maximal statistics can be computationally intensive and statistically conservative.
In comparison, our method is based on the maximum of single-edge $p$-values which becomes faster in computation and more efficient in detecting weaker signals. Based on that, this paper also introduces a novel KHAN algorithm which can adaptively  select homological features from the network persistent barcode over continuous filtration level. The KHAN algorithm is based on a discrete Gram-Schmidt $p$-value screening  method which untangles the probabilistic and algebraic   dependency among the candidate generators in the  homology group. The method also provides adaptive selection procedure for the change points on the filtration of the persistent barcode, and thus is computationally efficient. To the best of our knowledge, it is the first time for our paper to propose inferential method for selecting homological features from the network persistent homology with FDR controlled.

In theory, the proposed KHAN framework is guaranteed to achieve FDR control and be  powerful under reasonable sparsity and dependence conditions for general structural features and general graphical models. We introduce the uniform false discovery rate (uFDR) on continuous scales and show the proposed method can have the uFDR controlled at the given level only using the single-edge $p$-values, and developing  uniform false discovery control.  Compared to the Cram\'{e}r-type inequality for the maximal statistics in Gaussian graphical model in \cite{zhang2021startrek}, we propose novel theoretic analysis for FDR for t-statistics and covers the general family of graphical models. Our proof technique can handle the dependency among $p$-values of structural features and control the FDR as long as their correlations are not too strong. This analysis approach can be applied to the proof of uFDR control even the dependency among infinite filtration levels are involved. 

In application, we implement the proposed method to select residues of interest in the SARS-CoV-2 Spike protein during the binding process. Compared the existing works \citep{ou2020characterization, ray2021distant, liu2020potent} 
using marginal correlations to identify the vital residues involving the receptor-binding domain, our analysis applies the graphical model incorporating the conditional dependency and the network persistent homology considering the  structural changes. Therefore, we are able to take into account all four domains: NTB, RBD, linker, and S2 in the protein and their interactions during the conformational changes.

\subsection{Notations}
 
We denote set $A \times B=\{(a,b): a \in A, b \in B\}$. 
Let $G^*$ be a true graph with vertex set $V$ and edge set $E^*$, and let $G$ be an arbitrary graph that may differ in different instances.
We define $\bar E=\{(x,y):x,y\in V,x\neq y\}$ to be the complete edge set. 
We define $E(\cdot)$ to be a  function where   $E(G)$ gives the  edge set of graph $G$, and we define $V(G)$  be the vertex set  of graph $G$. We say $G \subseteq G'$ if and only if $V(G)\subseteq V(G')$ and $E(G)\subseteq E(G')$.   
$\Phi(\cdot)$ is the cumulative distribution  function of the standard Gaussian. For a sequence of random variables \(\{X_n\}_{n=1}^{\infty}\) and a scalar \(a\), we say \(X_n \leq a + o_P(1)\) if and only if  \(\lim_{n \to \infty} \PP(X_n - a > \epsilon) = 0\) for all \(\epsilon > 0\). We say $f=O(g)$ if $f\leq C g$ for some constant $C$, and   $f=o(g)$ and $g=\omega (f)$  if $C\rightarrow 0$. Define  the set  $[n] = \{1, 2, . . . , n\}$.  We use $|A|$ to denote the cardinality of a set $A$. For a vector $v=\left(v_1, \ldots, v_d\right)^T \in \mathbb{R}^d$, and $1 \leq q \leq \infty$, we define norm of $v$ as $\|v\|_q=$ $\left(\sum_{i=1}^d\left|v_i\right|^q\right)^{1 / q}$. In particular, $\|v\|_{0}=\sum_{1 \leq i \leq d} \mathbb{I} \{v_i \neq 0\}$.
%
For a matrix $\mathbf{M} \in \mathbb{R}^{r\times c}$, we use  $\mathbf{M}_{\cdot j}$ and $\mathbf{M}_{j \cdot}$ to denote the $j$-th column and row of $\mathbf{M}$ correspondingly. 
Additionally, let $\|\mathbf{M}\|_p=\max _{\|v\|_p=1}\|\mathbf{M} v\|_p$ for $p \geq 1$. In particular, 
$\|\mathbf{M}\|_{1}=\max _{1 \leq j \leq c}\sum _{i=1}^{r}|\mathbf{M}_{ij}|$.
If $\mathbf{M}$ is symmetric, let $\lambda_{\max }(\mathbf{M})$ and $\lambda_{\min }(\mathbf{M})$ denote the largest and smallest eigenvalues correspondingly. We denote  $\Phi(\cdot)$ as the standard normal cumulative density function.



\section{Graph Feature Selection}\label{sec:method}
In this section, we first provide some preliminaries on the concept of the graphical model and the graph feature selection problem. Then, we introduce the method for the sub-problem of selecting general graph features.

\subsection{Graphical Model and Persistent Homology}\label{sec:setup}

In our paper,    we consider the general graphical model.   
Let $\bX=\bigl(X_{1}, \cdots, X_{d}\bigr)^\top \in \mathbb{R}^{d}$ be a random vector indexed by the nodes of the true graph $G^* = (V, E^*)$, i.e., the node $j$ of the graph corresponds to the random variable $X_j$. 
The graph is assigned with the edge weights matrix $W^* =\{W^*_e\}_{e \in \bar{E}} \in \RR^{d\times d}$ and the graphical model implies the conditional dependency that $X_j$ is independent to $X_k$ conditioning on all other variables if and only if $(j,k) \not\in E^*$. 
The connection between the edges and the weights can typically be distilled into two distinct scenarios, both of which are covered in this paper:   
\begin{equation}\label{eq:sab}
\begin{aligned}
& \text{Scenario $(a)$: $E^*=\big\{e\in \bar E: |W^*_e| > 0\big\}$,}  \\
& \text{Scenario $(b)$: $E^*=\big\{e\in \bar E: W^*_e > 0\big\}$.}
\end{aligned}    
\end{equation}
Scenario (a) delineates a two-sided case, where both negative and positive values of $W^*$ are considered as edges. Conversely, Scenario (b) is indicative of a one-sided case, where only positive values of $W^*$ are considered edges.
For example, 
the Gaussian graphical model belongs to the Scenario (a).  In specific, $\bX$ follows a multivariate Gaussian distribution $N_d(\bm{0},\Sigma^*)$
with mean vector $\bm{0}$ and covariance matrix $\Sigma^*$. 
The weights are the entries of the precision matrix $W^*  = (\Sigma^*)^{-1}$ and $e \in E^*$ if and only if $W^*_e \neq 0$. 
On the other hand, the ferromagnetic Ising model belongs to the Scenario (b) with 
  the edge weights  $
     W^*_{uv}=\mathbb{E}[X_uX_v]-\tanh(\theta)$ and $e \in E^*$ whenever $W^*_e > 0$. 
We will discuss the details of these two models in  Examples  \ref{GS:ex:GGM} and \ref{ex:ising}.

Similar to \eqref{eq:sab}, we define the graph at filtration level $\mu$ as $G^*(\mu)=(V, E^*(\mu))$ under two scenarios:
\begin{equation*}
\begin{aligned}
& \text{Scenario $(a)$: $E^*(\mu)=\big\{e\in E^*: |W^*_e| >\mu\big\}$,}  \\
& \text{Scenario $(b)$: $E^*(\mu)=\big\{e\in E^*: W^*_e >\mu\big\}$.}
\end{aligned}    
\end{equation*}
Therefore, for a sequence of levels $\mu^{(1)}<\mu^{(2)}<\cdots<\mu^{(t)}$, we can obtain a filtration of edge sets $E^*(\mu^{(1)}) \supseteq E^*(\mu^{(2)}) \supseteq \ldots \supseteq E^*(\mu^{(t)})$ and a corresponding filtration of graphs
\begin{equation}\label{eq:ph:filtration}
  G^*(\mu^{(1)}) \supseteq G^*(\mu^{(2)}) \supseteq \ldots \supseteq G^*(\mu^{(t)}).  
\end{equation} 
We now define the network homology group. 
Given a graph $G =(V,E)$, we consider the $k$-th order chain group $\mathrm{C}_k(E)$ consisting of $k$-cliques and the cycle 
group  $ Z_k(E)$ is the kernel of the the boundary operator $\partial_k$. A persistent barcode is a graphical representation where horizontal line segments encode the birth-death filtration ranges of  generators across a given interval $[\mu_0, \mu_1]$, with the vertical axis represents the dimension $k$.
We illustrate the graph filtration and the barcode in Figure~\ref{fig:filtration} and refer the detailed definition of these concepts to Section \ref{sec:phapp} in the Supplementary Material. In this paper, we consider the cycle group as the homology group without taking the quotient in order to involve cliques as the features of interest as well.  In the following of the paper, we will refer the homology group to the cycle group and vice versa. The persistent homology groups are the sequence of groups $Z_k(E(\mu))$ for multiple filtration levels.
In this paper, we focus on selecting the generators of $Z_k(E(\mu))$, i.e.,  the homological features or the segments in the barcode, for all dimensions $k \le K$. We formulate the multiple hypotheses in \eqref{eqn:3} by the homological notations.  Denote the homology group  given edge set $E$ and the maximum dimension $K$ as 
$
 Z(E)=\bigoplus_{1 \le k\le K}Z_k (E),
 $
 where $K \le d$ is the maximum homology dimension of interest and define $Z(\mu) = Z(E^*(\mu))$. 
 In order to quantify the selection uncertainty,  we introduce the uniform False Discovery Proportion (uFDP) over a filtration interval $[\mu_0,\mu_1]$:  
 \begin{equation}\label{eq:ph:FDPFDR}
 \begin{aligned}
& \mathrm{uFDP} =\sup_{\mu \in [\mu_0,\mu_1]} \frac{\mathrm{rank}\big(\hat{Z}(\mu)\big)-\mathrm{rank}\big(Z(\mu)\cap \hat{Z}(\mu)\big)}{\max \bigl\{1, \mathrm{rank}\big(\hat{Z}(\mu)\big)\bigr\}},
 \end{aligned}
\end{equation}
where $\mathrm{rank}({Z}(\mu))$ is the group rank and the uniform False Discovery Rate becomes $\mathrm{uFDR} =\mathbb{E}[\mathrm{uFDP}]$.
Here at each filtration level $\mu$, $\mathrm{rank}\big(\hat{Z}(\mu)\big)$ counts the number of selected group generators, namely homological features, and $\mathrm{rank}\big(Z(\mu)\cap \hat{Z}(\mu)\big)$ counts that of correctly selected one. Our goal is to infer the homology group $\hat{Z}(\mu)$ for all $\mu \in [\mu_0,\mu_1]$ with uFDR controlled at the given level $q$.

\subsection{Graph Feature Selection}\label{sec:subalg}
We start the problem with the fixed filtration level $\mu$, which reduces the multiple hypotheses in \eqref{eqn:hypo} for the general graph features $F_j$. Let $\psi_{j}$ be the test for $H_{0j}: F_j \not\subseteq  G^*$ where $\psi_{j}=1$ if we reject $H_{0j}$ and $\psi_{j}=0$ otherwise. We aim to control the false discovery rate $\mathrm{FDR}=\mathbb{E}[\mathrm{FDP}]$ at given level $q$, where $
    \mathrm{FDP} = \big(\sum_{j \in \cH_{0}} \psi_{j} \big)/ \big(\max \bigl\{1, \sum_{j=1}^{J} \psi_{j}\bigr\}\big)$. Here  $\cH_0=\bigl\{1\leq j \leq J \biggiven  F_j \not\subseteq G^* \bigr\}$. 
In this paper, we propose a general multiple testing procedure to select graph features while  the  tail probability of  FDP and the FDR can be controlled below a given level $0<q<1$.

The main idea of our algorithm is to use the maximum $p$-values to select the significant graph features. 
Given any graph feature $F$, we assign its $p$-value as 
\begin{equation}\label{eq:pvalueF}
\alpha(F)=\max_{e\in E(F)}p_e.    
\end{equation}
Here  $p_e$ is any valid $p$-value for edge $e$. 
Recall that $W^*$ is the true edge weights matrix. In the following of the entire paper, we will present our method under the general graphical model. We assume that we have a generic edge weights estimator $\hat W$ which is asymptotically normal, i.e., $\sqrt{n} (\hat{W}_e - W^*_e) / \hat\sigma_e  \rightsquigarrow N(0,1)$, and $\hat\sigma_e^2$ is a generic estimated variance for $\sqrt{n}\widehat{W}_e$. We will present the general assumptions on these estimators in Assumption~\ref{assum:That:condition} and the concrete examples of such estimators will be discussed in Examples  \ref{GS:ex:GGM} and \ref{ex:ising}.
We then can use these generic estimators to obtain the $p$-values under the two scenarios: 
\begin{equation}\label{eq:That:condition0}
\begin{aligned}
    \text{Scenario $(a)$: }& p_e=2-2\Phi \big( \big|\sqrt{n}\hat{W}_e /\hat\sigma_e \big| \big);\\
    \text{Scenario $(b)$: }& p_e=1-\Phi(\sqrt{n}\hat{W}_e/\hat\sigma_e).
\end{aligned}
\end{equation}


Suppose we have the estimated $p$-values of all edges $\{p_e\}_{e \in \bar E}$.  We begin with an initial full edge  set and initialize graph feature $p$-values $\alpha(F_j)=1$ for all $j\in [J]$. 
We then filter out edges with $p$-values greater than a threshold $q$, denoting the resulting set of edges as $E_{0}(q)= \big\{e\in \bar E: p_e < q \big\}$. 
We only use \eqref{eq:pvalueF} to update those  $\alpha(F_j)$ values for graph features $F_j$ that can be embedded in the filtered graph $(V,E_0(q))$. Finally, we apply the Benjamini-Hochberg  procedure \citep{benjamini1995controlling}  to the set of updated $\alpha(F_j)$ values to select the significant graph features. 
 The Benjamini-Hochberg procedure ranks individual p-values from multiple tests and determines a threshold under which the null hypothesis can be rejected, thus controlling the FDR under $q$ among all significant results. Throughout this paper, we refer to this method as the BHq algorithm. The detailed pseudo-code is presented in  Algorithm~\ref{alg:1}. 

\begin{algorithm}[H]
    \caption{{General  Graph Feature Selection with FDR Control} \label{alg:1}}
    \SetKwInOut{Input}{Input}
    \SetKwInOut{Output}{Output}
    
    {\bf Input}: Edge $p$-values $p_e$ for $e\in \bar E$, FDR level $q$. 
    
    Initialize $\alpha_j = 1$ for $j \in [J]$\;
        
    Denote $E_{0}(q)= \big\{e\in \bar E: p_e < q \big\}$\;

    \For{$F \in \big\{F_j:E(F_j) \subseteq  E_0(q),j\in [J] \big\} $}{
    $\alpha(F)=\max_{e\in E(F)} p_e$\; 
    }
   
    Order $\alpha_{1},\ldots,\alpha_{J}$ as $\alpha_{(1)} \leq \alpha_{(2)} \leq \ldots \leq \alpha_{(J)}$ and set $\alpha_{(0)}=0$\;
    Let $j_{\max}=\max \bigl\{0 \leq j \leq J: \alpha_{(j)} < q j / J\bigr\}$ and $\hat  \alpha = q_{\max} j / J$\;
    
    {\bf Output}: {Reject $H_{0j}$, i.e., $\psi_j = 1$, for those $j$ such that $\alpha_{j} <\hat\alpha$.  Reject nothing if $j_{\max}= 0$.}
    
\end{algorithm}

Next we apply our general testing framework to two specific models: Gaussian graphical
model and Ising model. The core is how to estimate edge weights and $p$-values under these two models.

\subsection{Gaussian Graphical Model} \label{GS:sec:GGM}

\begin{example}[Gaussian Graphical Model (GGM)]\label{GS:ex:GGM}
The random vector $\boldsymbol{X}$ follows a multivariate Gaussian distribution $N_d\left(\mathbf{0}, \Sigma^*\right)$ with mean vector $\mathbf{0}$ and covariance matrix $\Sigma^*$. Let the precision matrix $\Theta^*=\Sigma^{*-1}$. We have $\Theta_{i j}^*=0$ if and only if $(i, j) \notin E^*$, therefore $E^*=\left\{(i, j) \in V \times V \mid i \neq j, \Theta_{i j}^* \neq 0\right\}$.
    
\end{example}

We consider the following parameter space for precision matrices
\begin{equation}\label{eqn:us}
    \mathcal{U}(s)=\Big\{\Theta \in \mathbb{R}^{d \times d} \mid 1 / \rho \leq \lambda_{\min }(\Theta) \le \lambda_{\max}(\Theta) \leq \rho, \max _{j \in[d]}\left\|\Theta_{\cdot j}\right\|_0 \leq s,\|\Theta\|_1 \leq M, \Theta=\Theta^{\top}\Big\}.
\end{equation}
Under such parameter space, we discuss how to estimate weights and  $p$-values for edges for the Gaussian graphical model. 
Given $n$ i.i.d. observations $\bX_1, \ldots, \bX_n$ following the distribution of $\bX$ and the GLasso estimator $\widehat{\Theta}$ for $\Theta^*$\citep{friedman2008sparse}.  We consider the unbiased estimator $\hat\Theta^{\mathrm{d}}$ and its variance estimate \citep{neykov2018unified} as
\begin{equation}\label{eq:thetad}
    \widehat{\Theta}_{uv}^{\mathrm{d}} =\widehat{\Theta}_{uv}- \frac{ {\widehat{\Theta}_{\cdot u}^{\top}\big(\widehat{\Sigma} \widehat{\Theta}_{\cdot v}-\mathbf{e}_{v}\big)} } {({\widehat{\Theta}_{\cdot u}^{\top} \widehat{\Sigma}_{\cdot u}})}, \qquad \hat \sigma_{uv}^2 = {\widehat{\Theta}_{uu}^{\mathrm{d}} \widehat{\Theta}_{vv}^{\mathrm{d}}  + (\widehat{\Theta}_{uv}^{\mathrm{d}} )^2 }
\end{equation}
where $\widehat{\Sigma}=1/n \sum_{i=1}^n \bX_i\bX_i^\top$ is the sample covariance matrix, $\widehat{\Theta}_{\cdot j}, \widehat{\Sigma}_{\cdot j}$ is the $j$-th column of matrix $\widehat{\Theta}, \widehat{\Sigma}$ respectively, and $\mathbf{e}_{k} \in \RR^d$ is the $k$-th canonical basis with only the $k$-th entry being 1.
$\hat\Theta^\mathrm{d}_{uv}$ can be decomposed into \citep{neykov2019combinatorial}
\begin{equation}\label{eqn:decomp}
 \sqrt{n} (\widehat{\Theta}^\mathrm{d}_{uv} - {\Theta}^*_{uv})=\frac{1}{\sqrt{n}} \sum_{i=1}^{n} \Theta_{\cdot u}^{* \top} \big(\bX_i\bX_i^\top \Theta^*_{\cdot v}- \mathbf{e}_{v}  \big)+o_{P}(1),  
\end{equation}
which can be shown to be asymptotically normal, i.e.,  $\sqrt{n}\big(\widehat{\Theta}_{uv}^{\mathrm{d}}- \Theta^*_{uv}\big) \rightsquigarrow N\big(0, \Theta^*_{uu}\Theta^*_{vv}+\Theta^{*^2}_{uv}\big)$.

So for  Gaussian graphical model, we take $ W_{uv}^*=  {\Theta}^*_{uv} $ and $\hat{W}_{uv}= \widehat{\Theta}_{uv}^{\mathrm{d}}$.
  And we estimate the $p$-value by $p_e=2\big(1-\Phi(|\sqrt{n}\hat{W}_e /\hat\sigma_e|)\big)$. {Throughout the paper, we interchangeably use $W_e$ and $W_{uv}$, $\Theta_e$ and $\Theta_{uv}$ for $e = (u,v)$.} The theoretical property of $\hat{W}$ is presented in   Proposition  \ref{lem:GGM}. 

\subsection{Ising Model}\label{GS:sec:Ising}

\begin{restatable}[Ferromagnetic Ising Model]{example}{exising}\label{ex:ising}
In the Ising model,  $\bX \in \{1,-1\}^{d}$ for each node has distribution
\begin{equation}\label{eqn:Ising0}
\mathbb{P}(\bX = \mathbf{x}) = \frac{1}{Z(\mathbf{w}^*)} \exp \Big( \sum_{(u, v) \in \bar E} w^*_{u v} x_{u} x_{v}\Big),
\end{equation} 
where $\mathbf{x} = (x_{1}, \cdots, x_{d})^\top$ is an observation of $\bX$ and $Z(\mathbf{w}^*) = \sum_{\mathbf{x} \in \{1,-1\}^{d}} \exp \Big( \sum_{(u, v) \in \bar E} w^*_{u v} x_{u} x_{v}\Big)$ is the partition function. The Ising model is ferromagnetic in the sense that the weights $\mathbf{w}_{uv}^*  \geq 0$ for all $(u, v) \in \bar E$.
\end{restatable}

Following  \cite{neykov2019property}, we assume there is a known lower bound $\theta$ for $w^*_e, e\in E^*$ and consider the following parameter space

\begin{equation}\label{model:Ising:para}
\mathcal{W}=\Big\{w^*: 
    \min_{e\in E^*}w^*_e\geq \theta,\|w^*\|_\infty \leq \Theta, s \tanh(\Theta)<\rho,s\geq 2\rho/(1-\rho)\Big\}, 
\end{equation}
where $\rho \in (0,1)$ is a constant. For the  ferromagnetic Ising model, we have $w^*_e > 0$ if and only if $e \in E^*$, and $w^*_e=0$ for $e \notin E^*$. 
However, estimating the weight $w^*_{uv}$ in the Ising model \eqref{eqn:Ising0} is not easy. 
 Lemma 3.1 in \cite{neykov2019property} shows that  the correlation between variables can be used to characterize the edge weight. Specifically,  $\E[X_uX_v]- \tanh \theta >0$ if and only if $w^*_{uv} >0$. Moreover,  correlation $\E[X_uX_v]$ is easier to estimate than $w^*_{uv}$.
 So we focus on studying  the weight  \begin{equation}\label{eq:Wstr-ising}
     W^*_{uv}=\mathbb{E}[X_uX_v]-\tanh(\theta). 
 \end{equation} 
Consider $n$ i.i.d. samples $\bX_1,\ldots, \bX_n \in \{1,-1\}^{d}$ from the Ising model, where $\bX_i=(X_{i1}, X_{i2}, \ldots X_{id})^{\top}$. We then estimate $W^*_{uv}$ in \eqref{eq:Wstr-ising} and the variance by the empirical average 
\begin{equation}\label{eq:That:Ising}
   \hat{W}_{uv}=\frac{1}{n}\sum_{i=1}^{n}  X_{iu} X_{iv}  - \tanh(\theta),\qquad \hat \sigma_{uv}^2 =  {1-\big(\frac{1}{n}\sum_{i=1}^{n}  X_{iu} X_{iv} \big)^2}.
\end{equation}
And we estimate the $p$-value by $p_e= 1-\Phi(\sqrt{n}\hat{W}_e /\hat\sigma_e)$. We have the theoretical guarantee for asymptotic normality of $\hat W_{uv}$ in  Proposition  \ref{lem:Ising}. 

\section{Inferential Analysis for  Persistent Homology}
\label{sec:ph}
In this section, we propose the method for selecting homological features with uFDR controlled. There are two major challenges compared to the graph feature selection. First, the selected graph features could be linearly dependent  in the homology group and we need to screen the independent  generators among them which involves complicated algebra operators. Second, in order to estimate the persistence of graph features, it is necessary to determine the generators cross continuous  filtration levels. 
To tackle the first issue, we apply the idea of Gram-Schmidt algorithm from linear algebra to homology group and combine it with the $p$-values.  For the second issue, we  introduce an upper-layer algorithm  that efficiently identifies finite-state homology groups across the continuous filtration levels.  This is to utilize the discrete structure of the graph and dynamically selecting the change point for persistent homology.

\subsection{Discrete Gram-Schmidt Algorithm for Graph Homology Group} 

We begin by introducing the method of selecting homological features at a given filtration level $\mu$. Similar to \eqref{eq:That:condition0},  we compute the $p$-values under  the two scenarios following
\begin{equation*}
\begin{aligned}
  \text{Scenario $(a)$: }& \hat{W}_e(\mu)= |\hat{W}_e|-\mu \text{ and }  p_e(\mu)=2-2\Phi({\sqrt{n}  \hat{W}_e(\mu)}/\hat{\sigma}_e); \\
  \text{Scenario $(b)$: }& \hat{W}_e(\mu)=\hat{W}_e-\mu \text{ and } p_e(\mu)=1-\Phi\big(\sqrt{n}  \hat{W}_e(\mu)/\hat{\sigma}_e\big).     
\end{aligned}    
\end{equation*}

The key idea of our algorithm is to keep track of a list of $p$-values corresponding to the selected generators in the homology group and then implement the BHq procedure.
First, we choose a prior edge set $E_0$, which by default considers these edges with $p$-values {smaller} than $q$. 
Let $\widetilde{E}$ be the edge set consisting the significant generators.  We initialize the algorithm with $\widetilde{E} = \varnothing$ and the list of $p$-values starts as empty as well. The algorithm aims to iteratively add edges from the candidate set $E_0$ to $\widetilde{E}$. In each iteration,  we choose the edge (denoted as $e^*$) with the smallest $p_e(\mu)$ over all $e\in E_0$. After adding $e^*$ to $\widetilde{E}$, we calculate the increase of the  rank of the 
homology group of $(V,\widetilde{E})$  and denote it by $\ell$. The algorithm to compute the rank of the homology group has been well established in the literature and it is computationally efficient with  linear complexity to the number of simplices \citep{fugacci2016persistent}.
 We then add $\ell$ repetitions of $p_{e^*}(\mu)$ into the $p$-value list. This can be interpreted as a discrete Gram-Schmidt procedure by adding the $p$-values for the $\ell$ generators independent to the selected generators in the previous iteration.  Compared to the Gram-Schmidt algorithm in the linear algebra, our method does not need to explicitly find the induced generators but simply add    $p_{e^*}(\mu)$ for $\ell$ times. We only need to identify the rank without selecting specific generators which  makes our approach computationally efficient. 
 The algorithm stops when  the candidate set $E_0$ becomes empty. The last step is to apply the BHq procedure on the $p$-value list and obtain a $p$-value threshold $\hat \alpha$ with the selected edge set {$\hat E(\mu) = \{e\in E_0:p_{e}(\mu)< \hat \alpha\}$} and the selected homology group $\hat Z(\mu) = Z(\hat E(\mu))$.  
 The detailed algorithm is presented in Algorithm~\ref{GS:alg:2}. 
 
 \vspace{10pt}
 
 \begin{algorithm}[H]\label{GS:alg:2}

    \caption{Discrete Gram-Schmidt (DGS) Algorithm}  \label{alg:ph:basis}
    \SetKwRepeat{Do}{do}{while}
    {\bf Function} DGS($\mu$)
    
    {\bf Input}: Filtration level $\mu$, FDR level $q$.  
    
    Initialize  $\widetilde{E}=\emptyset$, $E_0 = \{e\in \bar{E}:p_e(\mu) < q\}$ and $j=0$\; 

    \While{$E_0 \neq \emptyset$ }{
    Denote  $e^* = \argmin_{e\in E_0}p_e(\mu)$\;
    Calculate the increased rank $\ell=\mathrm{rank}(Z(\{e^*\}\cup \widetilde{E}))-\mathrm{rank}(Z(\widetilde{E}))$\;
  
    \If{ $\ell > 0$}{
     Assign $p$-values to the $\ell$ generators $\alpha_i=p_{e^*}(\mu)$ for $i = j+1, ..., j+\ell$\;
    $j=j+\ell$\;
    }

    $    \widetilde{E} = \widetilde{E} \cup \{e^*\}, E_0 = E_0 \setminus \{e^*\}$;
}

    Set $\alpha_{(0)}=0$ and let $j_{\max}=\max \Bigl\{0 \leq i \leq j: \alpha_{(i)} \leq q i /  \bar{J}\Bigr\}$, $\hat  \alpha = q j_{\max}  / \bar{J}$\;  
    
    {\bf Output}: Selected edge set {$\hat E(\mu) = \{e\in { \widetilde{E}}:p_{e}(\mu)< \hat \alpha\}$} and the selected homology group $\hat Z(\mu) = Z(\hat E(\mu))$.
\end{algorithm}

\subsection{KHAN Algorithm: FDR Selection for the Persistent Homology} \label{sec:barcode} To obtain the uniform selection for persistent homology over continuous filtration line $\mu \in[\mu_0,\mu_1]$, instead of applying Algorithm \ref{alg:ph:basis} point by point for infinite number of $\mu$'s, the key idea of the KHAN algorithm is to identify the finite change points of the barcode and apply the discrete Gram-Schmidt algorithm at each change point.

\begin{figure*}[htbp]	
\centering
\begin{tabular}{cc}
\includegraphics[width=1\textwidth,angle=0]{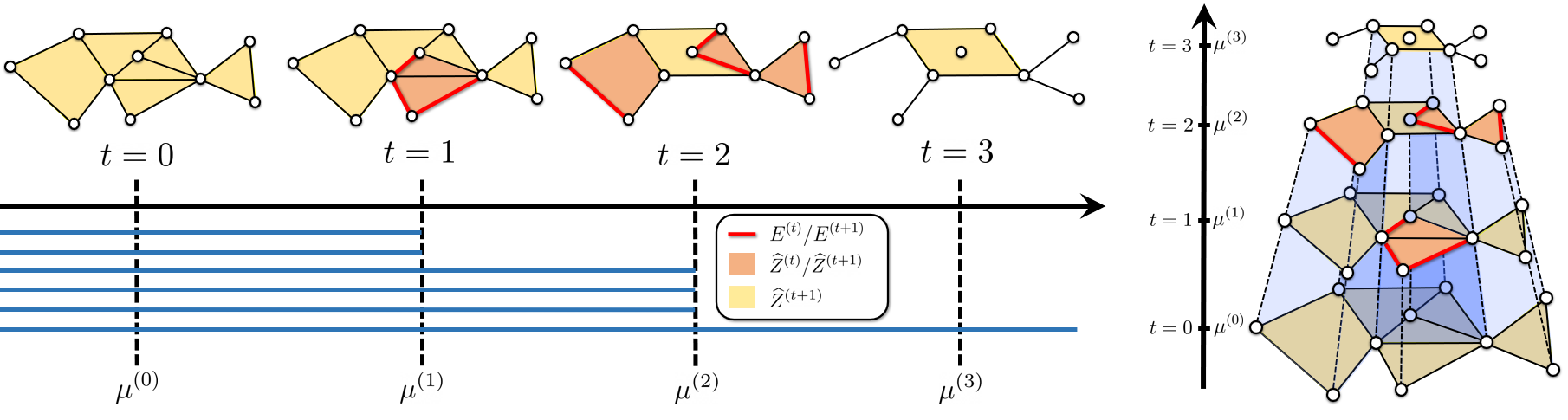}	
		\end{tabular}
\caption{ (Left) Illustration of the method of selecting homological features for persistent homology. At each iteration $t$, all the edges on the graph make the filtered edge set $E^{(t)}$, the weight of the red edge(s) represents the next filtration level $\mu^{(t+1)}$, the yellow cycles represent the remaining generator(s) in the next iteration, and the orange cycle(s) represents the disappearing generator(s) in the next iteration. The blue horizontal lines represent the life time of each cycle along the filtration of graph. (Right) Illustrations in 3D of the filtration of graphs when the filtration level $\mu$ increases.
} 
	\label{fig:dynamic}
\end{figure*}

We initialize the edge set and its homology group estimator at the minimum filtration level as $ E^{(0)}=\hat E(\mu_0)$ and $\hat Z^{(0)} = \hat Z(\mu_0)$ applying the discrete Gram-Schmidt algorithm.  We then iteratively identify the next change point by estimate the confidence lower bound for the edge weights. 
 In specific,  suppose at the $t$-th iteration, we have the estimated edge set ${E}^{(t)}$ and the  homology group $\hat{Z}^{(t)}$. 
We update the next change point under two scenarios as
\begin{equation}\label{GS:eqn:cp}
\begin{aligned}
\text{Scenario $(a)$: }& \mu^{(t+1)} =\min\Bigl\{|\hat{W}_e|- \Phi^{-1}(1-\alpha) {\hat{\sigma}_e}/{\sqrt{n}}:e\in E^{(t)}\Bigr\},  \quad\alpha = q \, \mathrm{rank}(\hat Z ^{(t)})/ (2\bar{J}), \\
\text{Scenario $(b)$: }& \mu^{(t+1)} =\min\Bigl\{\hat{W}_e- \Phi^{-1}(1-\alpha) {\hat{\sigma}_e}/{\sqrt{n}}:e\in E^{(t)}\Bigr\}, \quad\alpha = q \, \mathrm{rank}(\hat Z ^{(t)})/ \bar{J},      
\end{aligned}    
\end{equation}
where $\bar J = \mathrm{rank}(Z(\bar{E}))$ whose order is studied in Proposition \ref{prop:dim} of  the Supplementary Material.
Recall that $\hat{W}_e$ and $\hat{\sigma}_e^2$ are the generic weights and variance estimators and their examples under Gaussian graphical models and Ising models can be found in Sections \ref{GS:sec:GGM} and \ref{GS:sec:Ising}. 
We then update the edge set $E^{(t+1)} = \hat E(\mu^{(t+1)})$ and the homology group $\hat{Z}^{(t+1)} = \hat Z(\mu^{(t+1)})$ by the discrete Gram-Schmidt algorithm. Then we can interpolate the persistent homology group estimator $\hat Z(\mu) = \hat Z( \mu^{(t)})$ for all $\mu \in [\mu^{(t)},  \mu^{(t+1)})$. 
We repeat the process until $ \mu^{(t+1)} > \mu_1$. If $\mu_1 = \infty$, we stop the algorithm when the  estimated homology group becomes empty. 
 The detailed KHAN algorithm is presented in  Algorithm~\ref{GS:alg:3}. 
 See Figure \ref{fig:dynamic} for the illustration of the procedure. The algorithm is computationally efficient as the maximum  number of  iterations is $|E^{(0)}|$ and the implementation of discrete Gram-Schmidt only involves rank computation without basis selection.

 \vspace{7pt}

\begin{algorithm}[H] \label{GS:alg:3}
\caption{KHAN:  {\bf K}-Dimensional Persistent {\bf H}omology 
{\bf A}daptive Selectio{\bf N} Algorithm}\label{alg:ph:barcode}
    \SetKwRepeat{Do}{do}{while}
    
    {\bf Input}: {Estimated edge weights $\{\widehat{W}_e\}_{e \in \bar E}$} , FDR level $q$. 
       
    Initialize  $t=0$,  
     edge set and homology group at level $\mu_0$ as $ E^{(0)}=\hat E(\mu_0)$ and $\hat Z^{(0)} = \hat Z(\mu_0)$ applying DGS($\mu_0$) in Algorithm~\ref{GS:alg:2} \;
    
\While{    $\hat{Z}^{(t)}\neq \emptyset$}{


Determine the next change point using \eqref{GS:eqn:cp}\;   

Update the edge set $E^{(t+1)} = \hat E(\mu^{(t+1)})$ and the homology group $\hat{Z}^{(t+1)} = \hat Z(\mu^{(t+1)})$ by applying DGS($\mu^{(t+1)}$) \;

$t=t+1$\;

    }
 
 {\bf Output}: $\hat{Z}(\mu)=\hat{Z}^{(s)}$ for $\mu \in [\mu^{(s)},  \mu^{(s+1)})$  and $\hat{ Z}(\mu_0)=\hat{ Z}^{(0)}$, $\hat{ Z}(\mu)=\emptyset$ if $\mu>\mu^{(t)}$.
\end{algorithm}

\section{Theoretical Results} \label{GS:sec:thm}

In this section, we provide the theoretical results for the graph feature selection algorithm proposed  in Section \ref{sec:method}, and persistent homology selection algorithms proposed  in Section \ref{sec:ph}.
 
\subsection{Graph Feature Selection}\label{sec:the:subselect}

To ensure the FDR is adequately controlled in the  multiple testing problem \eqref{eqn:hypo} for general graphical models, it we need  to impose assumptions for the generic estimator $\hat{W}_e$ and $\hat{\sigma}_e$. In specific, we need the following  assumption on the asymptotic normality of $\hat W_e$ and the consistency of the variance estimator $\hat{\sigma}_e$.
\begin{restatable}{assumption}{assumThatcondition}\label{assum:That:condition}For some constant $C>0$, there exist i.i.d. random variables $\xi_{i}(e)$ for $i = 1, 2, \ldots, n$  with bounded $\psi_1$-orlicz norm and $W_e=\sum_{i=1}^n\xi_{i}(e)/n$ such that 
\begin{equation*}
   \sup_{e\in \bar E}\Bigl(\frac{\sqrt{n}\big|\hat{W}_e-W_e-W^*_e\big|}{\sigma_e\sqrt{\log d} }+\Big|\frac{\sigma_e}{\hat{\sigma}_e}-1\Big|\Bigr)\leq\frac{C}{\log d}+o_P(1),
 \end{equation*}
where  $\sigma_e^2 = \Var(\xi_1(e))$.
\end{restatable}

The assumption imposes that $\hat W_e$ can be decomposed into three parts: (1) the dominating error term $W_e$, which is asymptotically normal; (2) $W^*_e$, the true edge weight; and (3) the remainder term, which is at the order of $1/\sqrt{n\log d}$. 
Such decomposition has been achieved under different graphical models \citep{van2014asymptotically, neykov2019property}. 
Next we provide two examples of edge weight estimators $\hat W_e$ satisfying Assumption \ref{assum:That:condition} under the Gaussian graphical model  and the Ising model.

Recall that in Section \ref{GS:sec:GGM}, the Gaussian graphical model use the precision matrix as the weight $ W_{uv}^*=  {\Theta}^*_{uv} $ and estimated by $\hat{W}_{uv}= \widehat{\Theta}_{uv}^{\mathrm{d}}$ in \eqref{eq:thetad}.
 By \eqref{eqn:decomp}, the dominating error term  $W_{uv}$ becomes
\begin{equation}\label{eq:That:GGM2}  
 W_{uv}=\frac{1}{n }\sum_{i=1}^n\Theta^*_{u\cdot}(\bX_i\bX_i^\top\Theta^*_{\cdot v} - \mathbf{e}_{v}). 
  \end{equation}
We have the following proposition that validates Assumption~\ref{assum:That:condition} under the Gaussian graphical model. The proof of this proposition is deferred to   Section~\ref{sec:pflem:GGM} in the Supplementary Material.
\begin{restatable}[]{prop}{lemGGM}\label{lem:GGM} 
For GGM with parameter space $\Theta \in {\cal U}(s)$ defined in \eqref{eqn:us}, suppose $s^2\log^4(dn)/n=o(1)$, the estimator $\hat{W}_e$ and $\hat \sigma_{e}^2$ defined in \eqref{eq:thetad} with $W_{uv}$ in \eqref{eq:That:GGM2} satisfy Assumption~\ref{assum:That:condition}.
\end{restatable}

For the case of Ising model introduced in Section \ref{GS:sec:Ising}.
Recall that we  consider  the weight  $
     W^*_{uv}=\mathbb{E}[X_uX_v]-\tanh(\theta) 
$
and its estimator $\hat W_{uv}$ in \eqref{eq:That:Ising}. Then the dominating error term becomes \begin{equation}\label{eq:That:Ising2}
     W_{uv}=\frac{1}{n}\sum_{i=1}^{n} X_{iu} X_{iv}-\mathbb{E}[X_uX_v].
\end{equation}
We have the following theoretical guarantee for asymptotic normality of $\hat W_{uv}$.  The proof of this proposition is left to Section~\ref{sec:pflem:Ising} in the Supplementary Material.

\begin{restatable}[]{prop}{lemIsing}\label{lem:Ising}  For the Ising models in \eqref{eqn:Ising0} with the parameter space defined in \eqref{model:Ising:para},  the estimator $\hat{W}_e$ and $\hat \sigma_{e}^2$ defined in \eqref{eq:That:Ising} with $W_{uv}$ in \eqref{eq:That:Ising2} satisfy Assumption~\ref{assum:That:condition}.
 \end{restatable}

We next introduce the signal strength condition  required for our theoretical analysis. This condition helps to identify the subset of hypotheses with sufficiently strong signals within the context of  graph feature selection. Specifically, we define the set of hypotheses with strong signal strength under the two scenarios as 
\begin{align}
\begin{split}\label{eqn:h1tilde}
     \text{Scenario $(a)$: }&\widetilde{{\cal H}}_1=
\Big\{j\in {\cal H}_1:{ \min_{e\in E(F_j)}|W^*_e|} \geq C\sqrt{\log d/n} \Big\};\\
\text{Scenario $(b)$: }&\widetilde{{\cal H}}_1=
\Big\{j\in {\cal H}_1: \min_{e\in E(F_j)}{ W^*_e} \geq C\sqrt{\log d /n} \Big\}.
\end{split}
\end{align}

To characterize the dependence among the multiple tests in \eqref{eqn:hypo}, we introduce \begin{equation}\label{GS:eq:N}
   N_j=\{e\in E(F_j):W^*_e=0\}, \text{ for } j=1,2,\ldots,J,
\end{equation} which is the null edge set for the structure $F_j$ in the graph $G^*$. 
Recall that $\xi_{i}(e)$ is the summand of $W_e=\sum_{i=1}^n\xi_{i}(e)/n$ defined   in Assumption~\ref{assum:That:condition}.
{We define the dependence level as 
\begin{align}
\nonumber S = \big|\{({j_1},{j_2}): & ~j_{1}, j_{2} \in {\cal H}_{0}, j_1\neq j_2, \exists~ e_1 \in N_{j_1},e_2 \in N_{j_2}, \\
&\text{ s.t. } \big|\text{Cov}(\xi_1(e_1),\xi_1(e_2)) \big| \geq 
C(\log d)^{-2} (\log\log d)^{-1} \}\big|, \label{eq:S-dep}
\end{align}
where $C$ is a sufficiently large constant. The dependence level quantifies the dependency among multiple  hypotheses by the number of graph feature pairs $F_{j_1}$, $F_{j_2}$ whose null edges covariances are larger than certain level.
}


Next we present the assumption on the scaling conditions among the number of hypotheses and the the dependence level $S$. 
Recall that $\cH_0=\bigl\{1\leq j \leq J \biggiven  F_j \not\subseteq G^* \bigr\}$, ${\cal H}_1=\bigl\{1\leq j \leq J \biggiven  F_j \subseteq G^* \bigr\}$ and $\widetilde{{\cal H}}_1 \subseteq \cH_1$ is defined in \eqref{eqn:h1tilde}. 
We assume the following condition.
\begin{assumption}[Signal strength and dependence conditions]\label{ass:spadepend} 
We assume   
\[ 
\underbrace{ \frac{J}{|{\cal H}_0|\cdot |\widetilde{{\cal H}}_1|}}_{\text{Signal strength}} +\underbrace{\frac{JS}{|{\cal H}_0|^2|\widetilde{{\cal H}}_1|}}_{\text{\rm Dependence effect}}=o\Big(\frac{1}{\log d}\Big).
\]
\end{assumption}
Assumption \ref{ass:spadepend} comprises two types  conditions.  The first condition involves the numbers of null hypotheses and the alternatives with strong enough signal strength. In specific, 
if the number of the null hypotheses $|\cH_0|$ is of the same order of $J$, then it suffices to validate the first condition when the number of alternatives with strong signals has $|\tilde \cH_1| = \omega(\log d)$. 
The second condition balances the dependence level and the number of hypotheses. 
Again, if  $|\cH_0| \asymp J$, the second term is satisfied when $S = o\big({|{\cal H}_0| \cdot |\widetilde{{\cal H}}_1|}/{\log d}\big)$. {Similar assumption is also imposed by \cite{zhang2021startrek}, but our assumption is weaker as we only take into account the edges pairs whose covariance is significantly large in \eqref{eq:S-dep} but Eq.(5.3) in \cite{zhang2021startrek} counts any non-zero covariance  which makes their dependence level significantly larger.}
In the following, we will show that such condition can be validated    for a wide range of graphs under the Gaussian graphical model and the Ising model. We denote $M=\max_j|V(F_j)|$ as the maximum size of the graph features of interest and we can give an explicit upper bound  on $S$ to simplify the second condition into a direct sufficient condition on $|\widetilde{{\cal H}}_1|$.
\begin{proposition}\label{prop:GGM}   
For the Gaussian graphical model with the weights $\Theta \in {\cal U}(s)$ where $s$ is the maximum degree of the graph $G^*$, we have  $S=O(|{\cal H}_0|d^{M-2}s^2)$.
\end{proposition}
The proof is deferred to Section \ref{sec:pfprop:GGM} in the Supplementary Material.
As  $M=\max_j|V(F_j)|$, then the maximal possible number of tests $J = |\cH_0| + |\cH_1|\asymp d^M$.
By  Proposition
\ref{prop:GGM}, under GGM, if the number of null hypotheses is balanced with the alternatives, i.e., $|{\cal H}_0|\asymp |{\cal H}_1| \asymp d^M$, then the dependence effect in Assumption \ref{ass:spadepend}  can be implied when $|\widetilde{{\cal H}}_1| = \omega(s^2 d^{M-2} \log d)$.

The following proposition provides an upper bound of $S$ under the Ising model.

\begin{proposition}\label{prop:Ising}
For ferromagnetic Ising models, assume that the graph with $d$ nodes consists of $d/s$ trees and each tree has $s$ nodes. We have  $S=O(|{\cal H}_0|d^{M-1}s)$.
\end{proposition}
The proof is deferred to Section \ref{sec:pfprop:Ising} in the Supplementary Material.
Similarly, under the conditions of Proposition~\ref{prop:Ising}, if $|{\cal H}_0|\asymp  |{\cal H}_1|$, the dependence effect can be implied when $|\widetilde{{\cal H}}_1| = \omega(s d^{M-1} \log d)$.

\subsubsection{Theoretical Results for  Graph Feature Selection}

The following theorem shows that Algorithm~\ref{alg:1} can control the FDR of graph feature selection.
\begin{restatable}[General FDR control]{theorem}{thmgeneralfdp}\label{thm:generalfdp} 
 Suppose $\log d/\log n=O(1)$. Under Assumptions~\ref{assum:That:condition} and~\ref{ass:spadepend}, we have that the $\{\psi_j\}_{j \in [J]}$ determined by Algorithm \ref{alg:1} satisfies
\[
\mathrm{FDP}\leq q\frac{|{\cal H}_0|}{J}+o_{P}(1),\qquad \mathrm{FDR}\leq q\frac{|{\cal H}_0|}{J}+o(1).
\]
\end{restatable}
The proof of this theorem is deferred to Section \ref{sec:pfthm:generalfdp} in the Supplementary Material. On the other side, the next  theorem shows that the algorithm is also powerful.
\begin{restatable}[Power Analysis]{theorem}{thmpower}\label{thm:power}
For $\{\psi_j\}_{j \in [J]}$ determined by Algorithm \ref{alg:1}, under the same conditions of Theorem \ref{thm:generalfdp},  we have 
\begin{equation*}
\mathbb{P}\Bigl(\psi_j=1,\text{ for all }j\in \widetilde{{\cal H}}_1\Bigr) = 1 -  o(1),
\end{equation*}   
where $\widetilde{{\cal H}}_1$ is defined in \eqref{eqn:h1tilde}.
\end{restatable}
This theorem guarantees that the true alternative hypotheses satisfying certain signal strength condition will be selected with probability going to $1$. The proof of this theorem can be found in Section \ref{sec:pfthm:power} in the Supplementary Material.
 We have the following corollaries for the above general  results under the Gaussian and Ising graphical models. The  proofs are postponed to Sections~\ref{sec:cor:FDRPOW:GGM} and \ref{sec:pfcor:FDRPOW:Ising} in the Supplementary Material. 

\begin{restatable}[Graph feature selection under GGM]{corollary}{corFDRPOWGGM}\label{cor:FDRPOW:GGM} Under the same conditions of Proposition \ref{lem:GGM}, if  $\log d/\log n=O(1)$ and Assumption \ref{ass:spadepend} is satisfied, we have
\[
\mathrm{FDP}\leq q\frac{|{\cal H}_0|}{J}+o_{P}(1),\qquad \mathrm{FDR}\leq q\frac{|{\cal H}_0|}{J}+o(1),
\]
 and the uniform power control 
\begin{equation*}
\mathbb{P}\big(\psi_j=1,\text{ for all }j\in \widetilde{{\cal H}}_1\big) = 1 -  o(1).
\end{equation*}   

\end{restatable}

\begin{restatable}[Graph feature selection under Ising models]{corollary}{corFDRPOWIsing}\label{cor:FDRPOW:Ising} Under the same conditions of Proposition \ref{lem:Ising}, if  $\log d/\log n=O(1)$ and Assumption \ref{ass:spadepend} is satisfied, we have
\[
\mathrm{FDP}\leq q\frac{|{\cal H}_0|}{J}+o_{P}(1),\qquad \mathrm{FDR}\leq q\frac{|{\cal H}_0|}{J}+o(1),
\]
and the uniform power control
\begin{equation*}
\mathbb{P}\big(\psi_j=1,\text{ for all }j\in \widetilde{{\cal H}}_1\big) = 1- o(1).
\end{equation*}   
\end{restatable}

\subsection{Theoretical Results for Persistent Homology} \label{sec:thm:ph}

In this section, we provide the theoretical justification of the proposed KHAN algorithm. {Similar to the dependence level $S$ in Assumption~\ref{ass:spadepend}, we need to quantify the dependency among different homological features at different filtration levels. We first introduce the concept of critical edges which may change the rank of the homology group. 
Let the edges in $\bar E$ be ordered as $e_{(1)}, e_{(2)}, \ldots, e_{(|\bar E|)}$ such that $W^*_{e_{(1)}} \ge \ldots \ge W^*_{e_{(|\bar E|)}}$ and denote $E_{(i)} = \{e_{(1)},  \ldots, e_{(i)}\}$ for all $i   \in [|\bar{E}|]$ and $E_{(0)} = \emptyset$. Recall that $\bar J = \mathrm{rank}(Z(\bar{E}))$. To characterize the critical edges which may increase the rank in Algorithm~\ref{GS:alg:2}, we construct the critical edge list $\bar{e}_1,  \ldots, \bar{e}_{\bar J}$ by the following procedure.  
Let $\ell_i=\mathrm{rank}(Z(E_{(i)}))-\mathrm{rank}(Z(E_{(i-1)}))$ for all $i \in [|\bar{E}|]$. Then the $j$-th critical edge  $\bar{e}_j=e_{(i)}$ if $j$ has $\sum_{u=1}^{i-1}\ell_u<j\leq \sum_{u=1}^i\ell_u$. }
We then define the dependence level  at the filtration level $\mu$ by
\begin{equation}
\begin{aligned}
\bar{S}(\mu)&=\Big|\Bigl\{
 (j_1,j_2): j_1 \neq j_2, \bar{e}_{j_1},\bar{e}_{j_2}\notin E^*(\mu),  
 \big|\text{Cov}\big( \xi_1(\bar{e}_{j_1}),\xi_1(\bar{e}_{j_2}) \big) \big|  \geq 
\frac{C}{(\log d)^{2}\log (\log d)} \Bigr\}\Big|,
\end{aligned}    
\end{equation}
Similar to \eqref{eq:S-dep}, $\bar{S}(\mu)$ characterizes the dependence between the pair of critical edges $(\bar{e}_{j_1}, \bar{e}_{j_2})$.

{Similar to Assumption~\ref{ass:spadepend}, 
  we need to impose  the assumption on the scaling conditions of the number of hypotheses and dependence level. However, as the generators could be linearly dependant among each other in the homology group,  the homological features in  \eqref{eqn:3} are not as identifiable as the graph features in \eqref{eqn:hypo}.  Therefore, we bypass the identifiability issue by defining  the uFDP in \eqref{eq:ph:FDPFDR} through the identifiable rank difference instead of listing the null hypotheses explicitly. Following this idea, given any filtration level $\mu \in [\mu_0, \mu_1]$, we define the number of significant alternatives $|\widetilde{{\cal H}}_1(\mu)| =  \mathrm{rank}(Z(\mu + C\sqrt{\log d/n}))$ and the number of null hypotheses $   |\bar{{\cal H}}_0(\mu)|= \bar J -\mathrm{rank}(Z(\mu))$. Here we use the notation of cardinality without introducing the set simply in order to link the notations to the ones in Assumption~\ref{ass:spadepend}. We refer to Section~\ref{sec:pfthm:ph:generalfdppow} in the Supplementary Material for a detailed discussion on the definitions above.}
 \begin{assumption}\label{ass:phspadepend} For a sufficiently large constant $C$, we assume
\[  
 \underbrace{ \frac{\bar{J}}{|\bar{{\cal H}}_0(\mu_1)|\cdot|\widetilde{{\cal H}}_1(\mu_1)|}}_{\mathrm{Signal \;strength}} +\underbrace{\frac{\bar{J}\bar{S}(\mu_1)}{|\bar{{\cal H}}_0(\mu_1)|^2|\widetilde{{\cal H}}_1(\mu_1)|}}_{\mathrm{Dependence \;effect}}=o\Big(\frac{1}{\log d}\Big).
\]
\end{assumption}
Notice that it is sufficient to assume Assumption \ref{ass:phspadepend} only for $\mu=\mu_1$ 
without the uniform control on $\mu\in [\mu_0,\mu_1]$. This makes the assumption less stringent and easier to verify.
We have the following proposition to bound $\bar S(\mu)$ under the Gaussian graphical model and the Ising model. The proof is deferred to Section \ref{sec:disc:assm:phspadepend} in the Supplementary Material.
\begin{proposition}\label{prop:GGM-Ising}   
Define $R = \max_{i=1 }^{|\bar{E}|}\ell_i$ as the maximal increased rank.
For the Gaussian graphical model with $s(\mu)$ being the maximum degree of graph $G^*(\mu)$,  we  have
$$ \bar S(\mu)=O(|\bar{{\cal H}}_0(\mu)|\cdot s(\mu)^2 R).$$
For the ferromagnetic Ising model, assuming that the graph $G^*(\mu)$ with $d$ nodes consists of $d/s(\mu)$ trees, with each tree having $s(\mu)$ nodes, we   have
\[  \bar S (\mu) 
= O\big(|\bar{{\cal H}}_0(\mu)| \cdot s(\mu) dR \big).
\]

\end{proposition}

Based on Proposition \ref{prop:GGM-Ising}, under the Gaussian graphical model, it suffices to verify Assumption~\ref{ass:phspadepend}  if the number of strong signals  $|\widetilde{{\cal H}}_1(\mu_1)|=\omega( R s(\mu_1)^2\log d$). 
Under the Ising model, Assumption~\ref{ass:phspadepend} can be implied by $|\widetilde{{\cal H}}_1(\mu_1)|=\omega( R s(\mu_1)d\log d$).

The following theorem controls the uFDP in \eqref{eq:ph:FDPFDR}  for the KHAN algorithm. The proofs can be found in Section~\ref{sec:pfthm:ph:generalfdppow} in the Supplementary Material.

\begin{restatable}[uFDR and Power of KHAN]
{theorem}{thmphgeneralfdppow}\label{thm:ph:generalfdppow}
Suppose $\log d/\log n=O(1)$ and  Assumptions \ref{assum:That:condition} and \ref{ass:phspadepend}   are satisfied. Algorithm~\ref{GS:alg:3} has
\[
\mathrm{uFDP}\leq q| {\cal H}_0(\mu_1)|/\bar{J}+o_P(1),\quad \mathrm{uFDR}\leq  q| {\cal H}_0(\mu_1)|/\bar{J}+o(1). 
\]
We also have the uniform power of the selections  
\begin{equation*}
\mathbb{P}\Bigl(Z \big(\mu + C\sqrt{\log d/n} \big) \subseteq \hat{ Z}(\mu) \text{ for all }\mu\in[\mu_0,\mu_1] \Big) = 1- o(1).
\end{equation*}
\end{restatable}

We have the following corollaries under the Gaussian graphical model and  the Ising model.

\begin{restatable}
{corollary}{corbarcodeGGM}\label{cor:ph:barcodeGGM} Suppose Assumption \ref{ass:phspadepend} is satisfied   and $\log d/\log n=O(1)$. 
For Gaussian graphical model with $\Theta^*\in {\cal U}(s)$, we further assume $s^2\log^4(dn)/n=o(1)$, then   Algorithm \ref{alg:ph:barcode}  satisfies 
\[
\mathrm{uFDP}\leq q+o_P(1),\qquad \mathrm{uFDR}\leq q +o(1).
\]
and  the uniform power control
\begin{equation*}
\mathbb{P}\Bigl(Z \big(\mu + C\sqrt{\log d/n} \big) \subseteq \hat{ Z}(\mu) \text{ for all }\mu\in(\mu_0,\mu_1) \Bigr) = 1- o(1).
\end{equation*}   
\end{restatable}
\begin{restatable}
{corollary}{corbarcodeIsing}\label{cor:ph:barcodeIsing} Suppose Assumption \ref{ass:phspadepend} is satisfied   and $\log d/\log n=O(1)$. 
For Ising models with $w^*\in {\cal W}$, Algorithm \ref{alg:ph:barcode}  satisfies 
\[
\mathrm{uFDP}\leq q+o_P(1),\qquad \mathrm{uFDR}\leq q +o(1),
\]
and the uniform power control
\begin{equation*}
\mathbb{P}\Bigl(Z \big(\mu + C\sqrt{\log d/n} \big) \subseteq \hat{ Z}(\mu) \text{ for all }\mu\in(\mu_0,\mu_1) \Bigr) = 1- o(1).
\end{equation*}   
\end{restatable}
The proofs can be found in Section~\ref{sec:coro:ph} in the Supplementary Material.

\section{Application to SARS-CoV-2 Spike Protein Data}

The COVID-19 pandemic has had a devastating impact worldwide. Infection by the causative agent, severe acute respiratory syndrome coronavirus 2 (SARS-CoV-2) involves the attachment of the receptor-binding domain (RBD) of its spike proteins to the ACE2 receptors on the peripheral membrane of host cells \citep{ou2020characterization}. The binding process begins with a conformational change in the spike protein, transitioning from a downward orientation to an upward orientation, thereby exposing the RBD to the receptor. Recent studies that search for therapeutics have primarily focused on the RBD that is highly prone to mutations. However, it is worth considering that other residues within the spike protein, apart from the RBD, could also be potential targets for inhibition \citep{liu2020potent}.  To explore this possibility, we apply our method to identify these potential target residues based on persistent homology, taking into account the protein's complex network structure.

We use the trajectory data from \citep{ray2021distant}. The data encompasses three distinct states of the SARS-CoV-2 spike protein: closed, partially open, and fully open. The protein consists of three chains, with each chain comprising $1146$ residues and four domains: N-terminal domain (NTD), RBD, linker, and S2. At a fundamental level, conformational changes in protein structure originate from intricate combinations of transitions between different states of the protein backbone's torsional angles, specifically $\phi$ and $\psi$, which are used as features of the residues. Additionally, $173$ features of pairwise distances for residues from the RBD are provided, which are calculated between them and from other parts of the spike near the RBD. These distance features capture the down-to-up transition of the RBD, and residues are identified to be important for the conformational change if their backbone torsion angles are highly correlated to these features. In total, we obtain $d = 7048$ features consisting of backbone torsion angles and distances, with the sample sizes $n$ of $8497$, $14888$, $8131$ for the closed, partially open, and fully open states, respectively. We build a $d$-dimensional Gaussian graphical model for each state to analyze the relationship among these features simultaneously, to gain insights into the specific residues that have essential functions in the infection process. We convert the angle features by the transform $\Phi^{-1}(\hat{P_j}(\cdot))$, where $\hat{P_j}(\cdot)$ denotes the empirical distribution of the $j$-th angle feature.

 \begin{figure*}[htbp!]
    \centering
    \includegraphics[width=.8\linewidth]{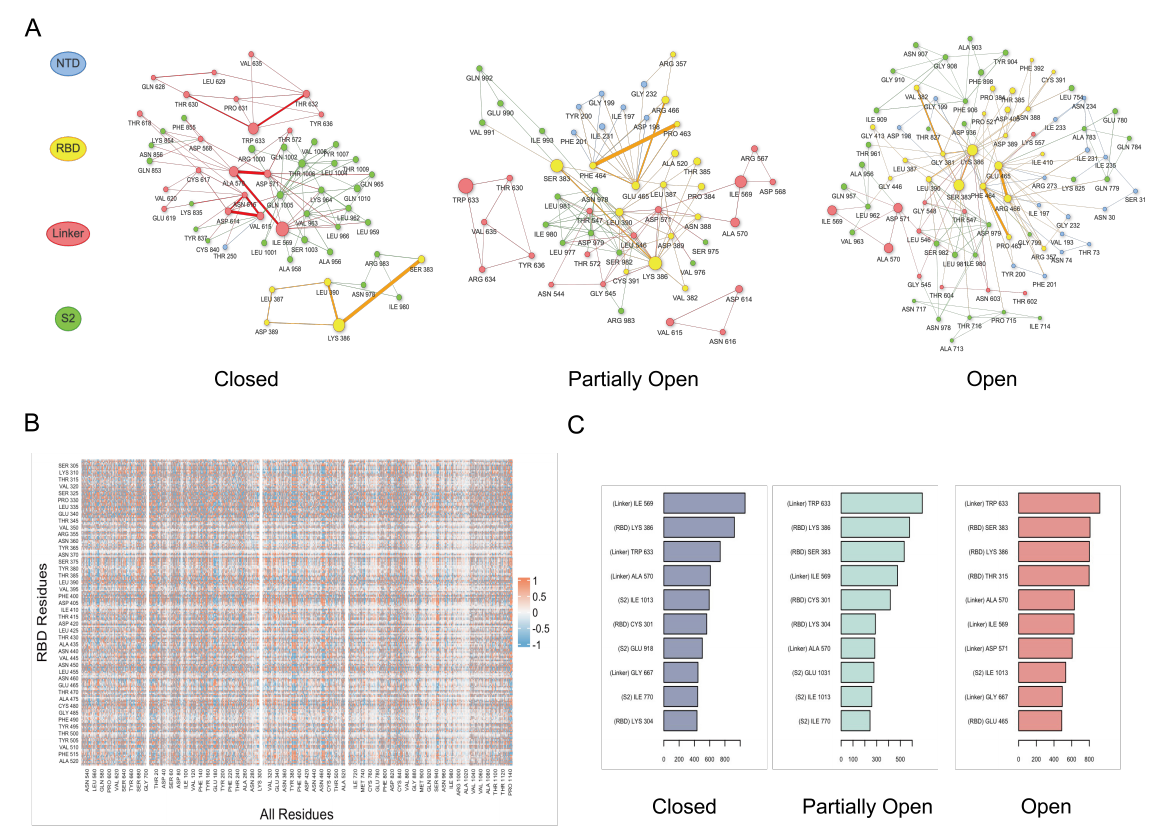}
    \caption{This figure presents a comprehensive visualization of the connectivity and correlation analysis among residues across different domains. Panel A highlights the strong connectivity of residues from the RBD (the yellow nodes) with those from other domains. Panel B showcases the correlation analysis, emphasizing the strong associations between residues within the RBD and those in the NTD, the linker domain, and a portion of the S2 domain. Panel C summarizes the top 10 residues identified at each stage based on their importance scores.}
    \label{fig:GraphHeatmapBar}
\end{figure*}

\begin{figure*}[htbp]	
\centering
\begin{tabular}{cccc}
\rotatebox{90}{{K386/S383}} & \includegraphics[width=0.30\textwidth]{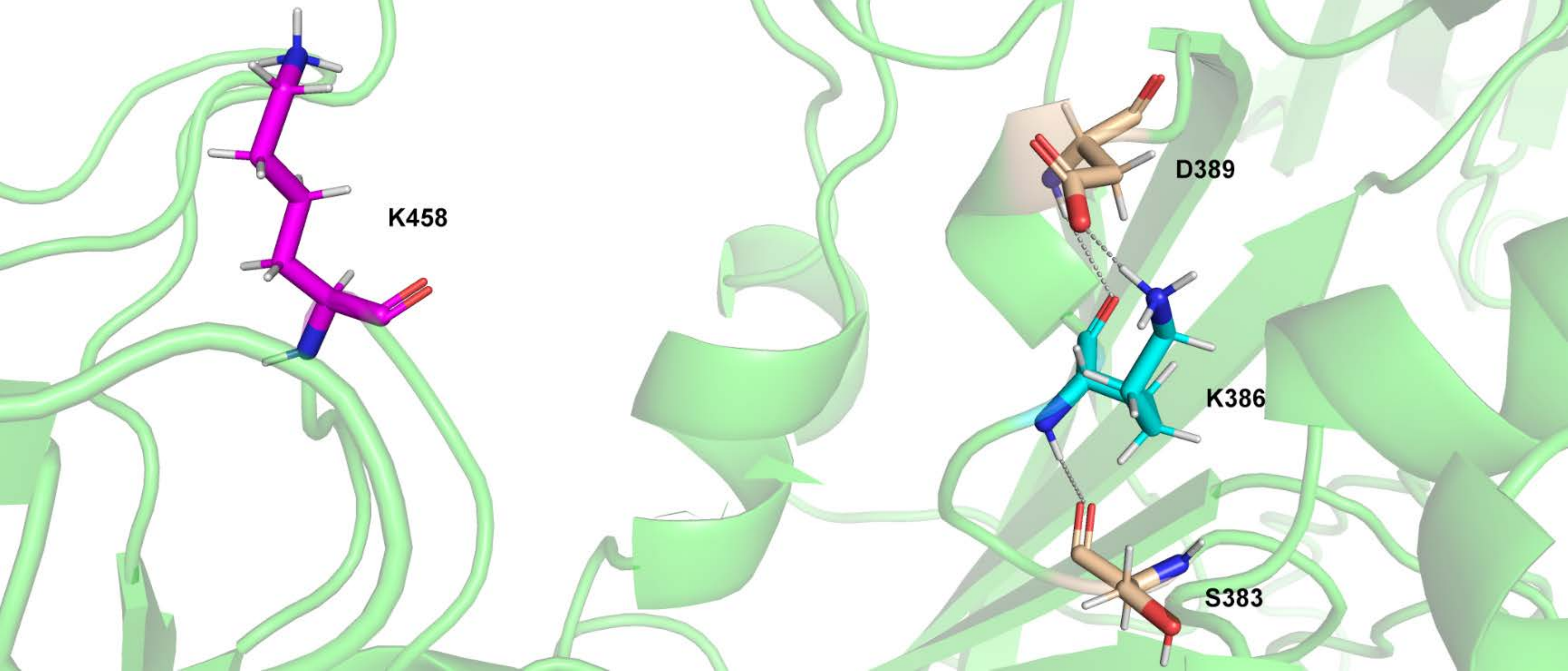}  &
 \includegraphics[width=0.30\textwidth]{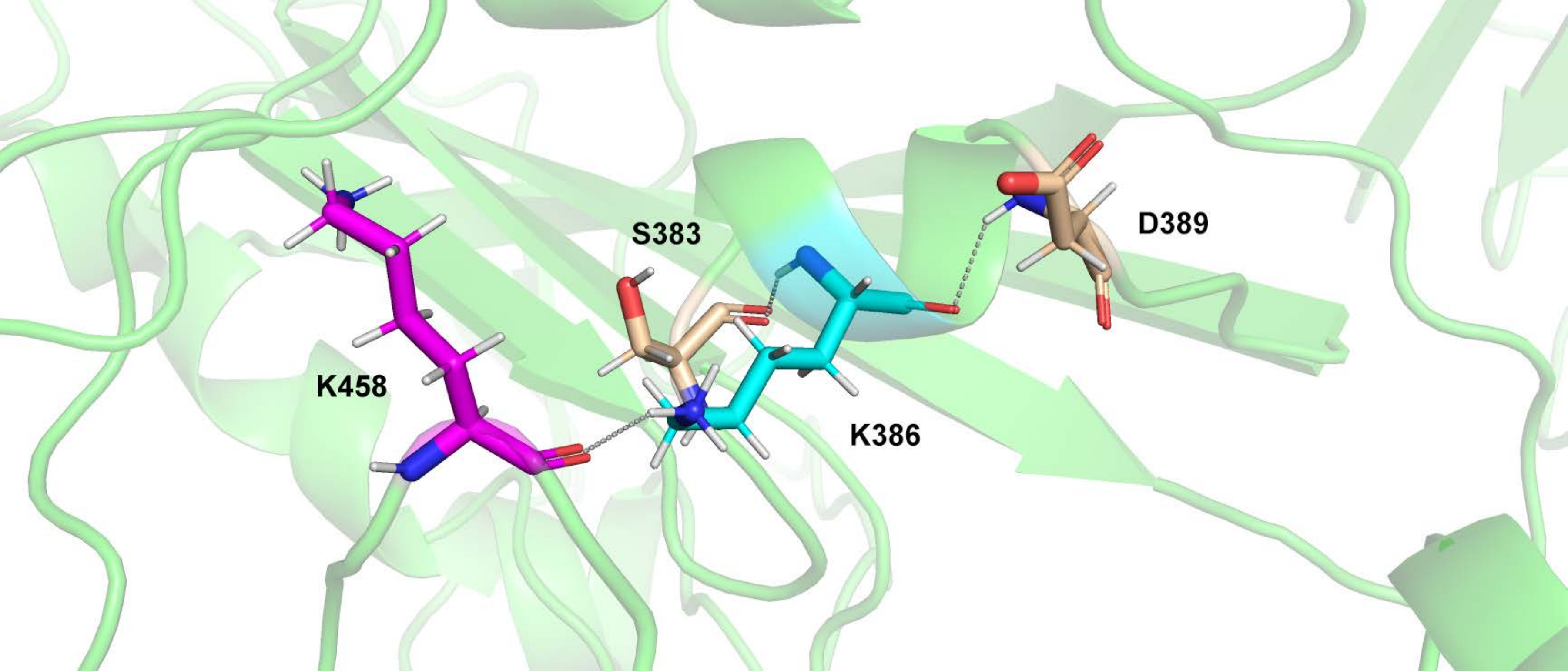} &  
 \includegraphics[width=0.30\textwidth]{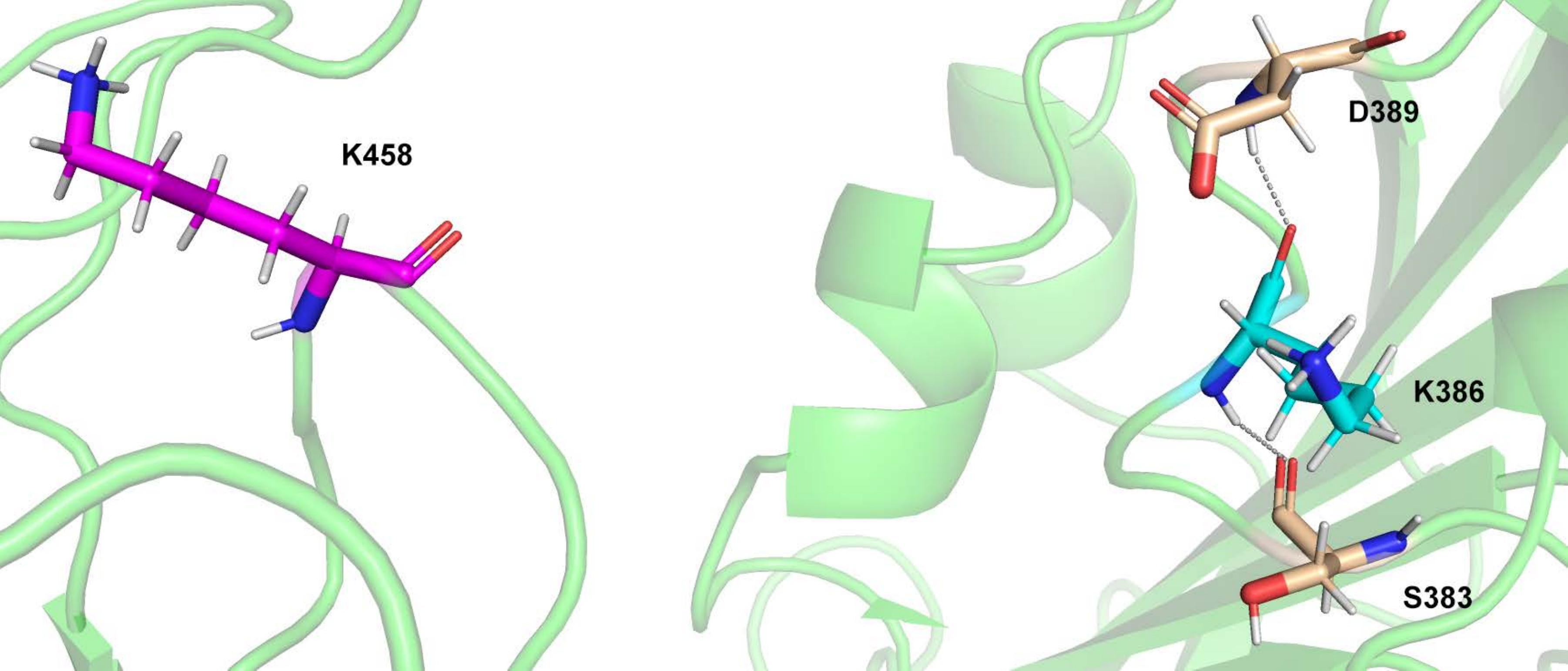}  \\ 
  \rotatebox{90}{ \centering{A570/I569} } & \includegraphics[width=0.30\textwidth]{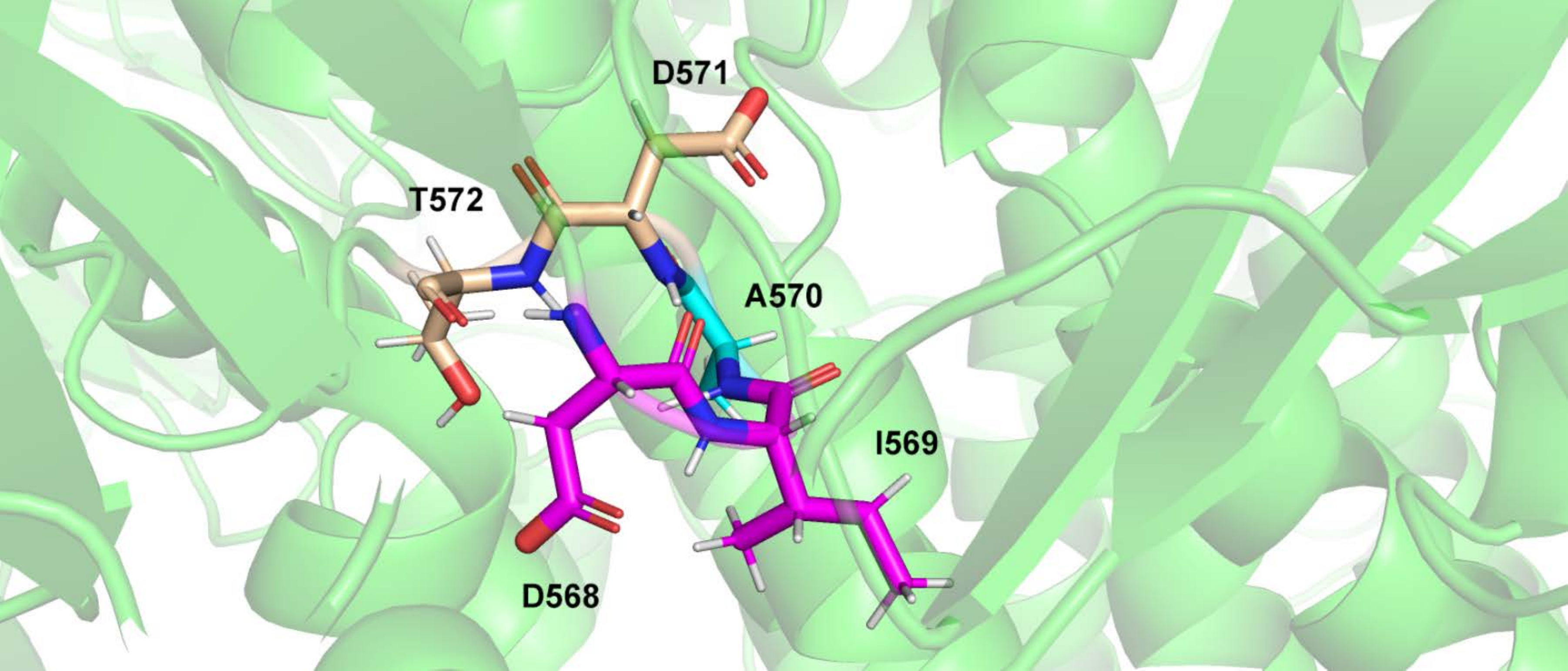}  &
 \includegraphics[width=0.30\textwidth]{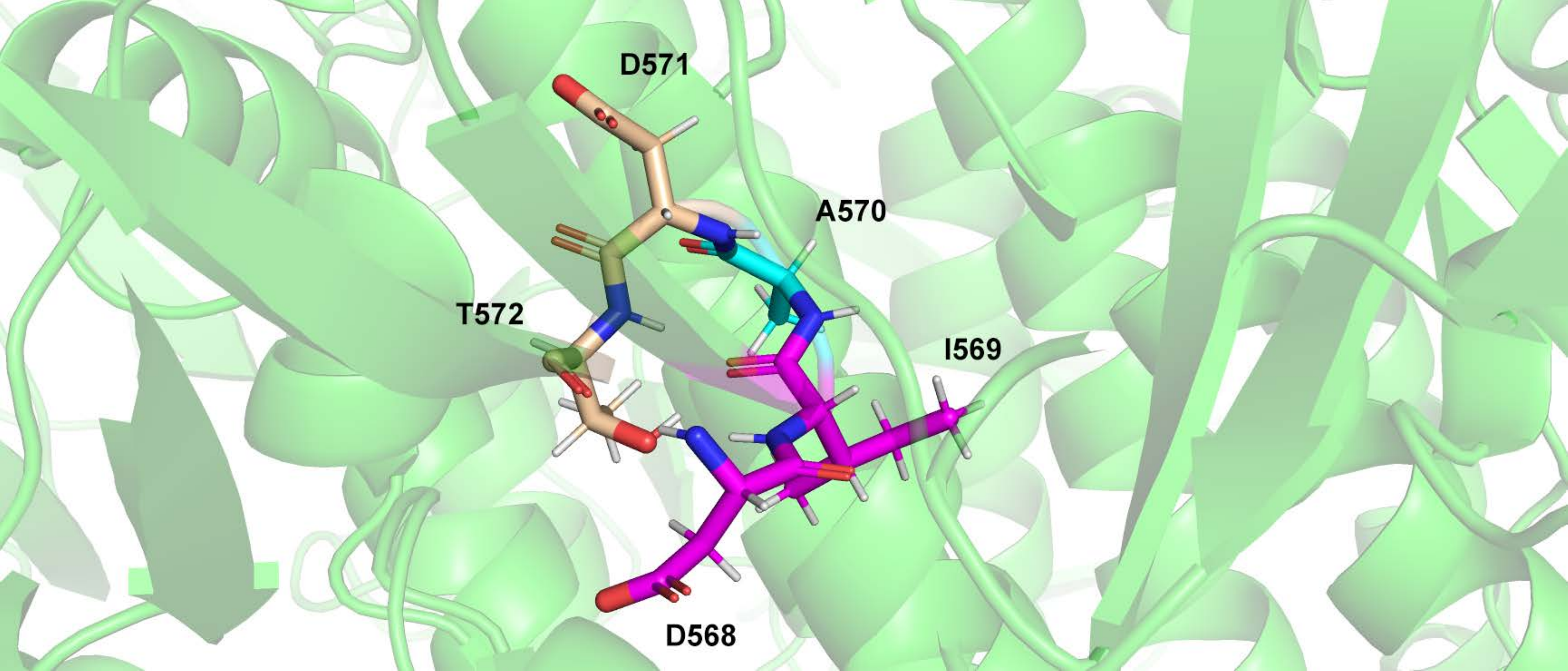} &  
 \includegraphics[width=0.30\textwidth]{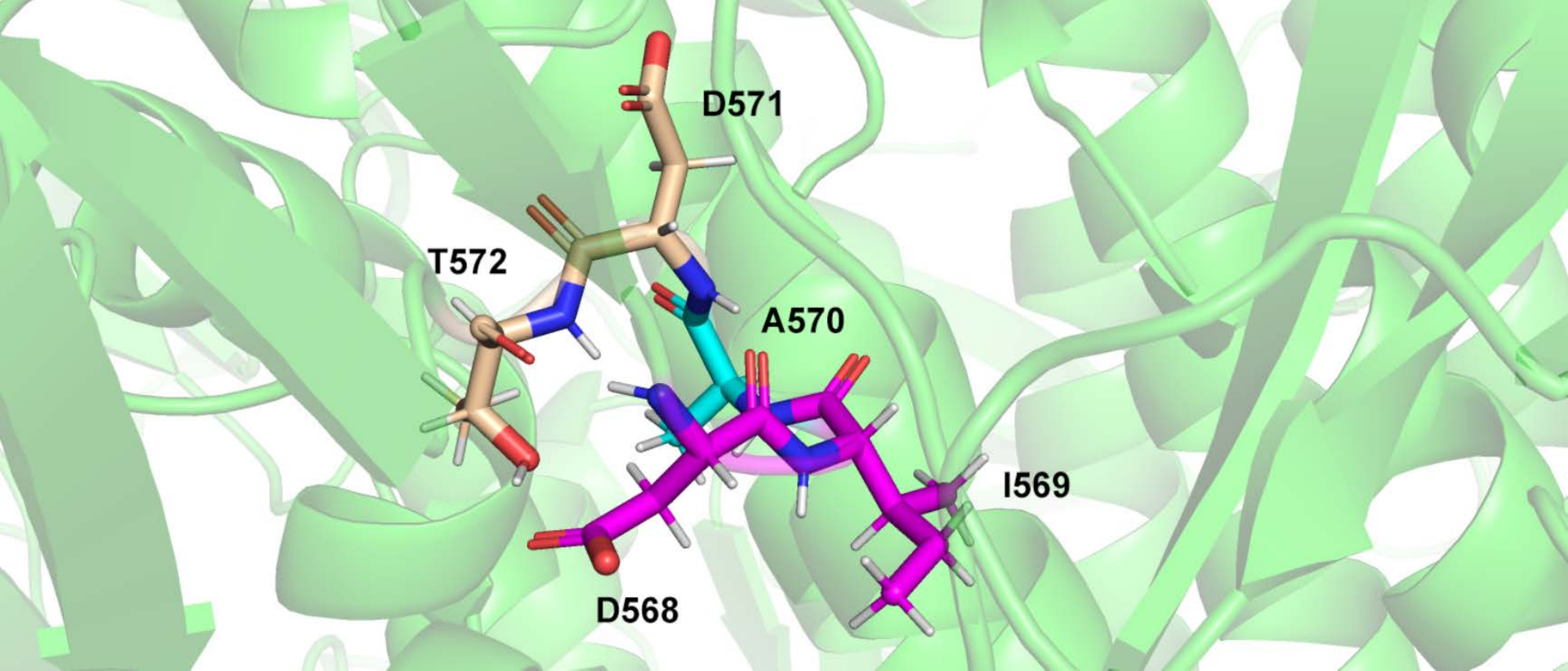}  \\ 
 &
 { (a) Closed}& 
  { (b) Partially Open} &
  { (c) Open}
\end{tabular}
\caption{Illustration of the conformational changes and disruption of hydrogen bonds in residues K386/S683 and A570/I569 across three stages.}
\label{fig:ResidueExample}
\end{figure*}

During each stage, residues exhibiting strong connections with numerous other residues can be considered as crucial. Specifically, we apply the proposed persistent homology selection algorithm to obtain the evolution of homology groups among the features. In order to evaluate the importance of the selected homological features, we define the importance score of each feature as the length of the filtration interval of each feature exists. Then, the importance score of each residue is calculated by summing over the importance scores of all the homological features containing the residue. The detailed algorithm for determining the importance scores is presented in Algorithm \ref{alg:hs}  in the Supplementary Material. 

A subset of residues outside the RBD  can exert significant influence on the conformational changes of the spike protein, thereby impacting viral infection. Specifically, residues located in the linker domain connecting the RBD and the S2 domain act as a pivotal hinge, facilitating the transition of the RBD from a downward to an upward position through torsional changes in the protein backbone. Residues within the S2 domain or the NTD play a crucial role in modulating the dynamics of the RBD or linker domain, owing to their spatial proximity to the RBD in the three-dimensional conformation. These phenomenon has been previously elucidated in \citep{ray2021distant}. 

To provide a more intuitive demonstration of the significant impact of residues from other domains on the RBD, we initially select $100$ homological features from each stage based on their edge weights and constructed graphs accordingly. By visualizing the graphs from the three stages as Figure~\ref{fig:GraphHeatmapBar}A, it is illustrated that the yellow nodes representing residues from the RBD exhibit strong connectivity with residues from other domains. This phenomenon is further evident in Figure~\ref{fig:GraphHeatmapBar}B, which is generated from the correlation analysis between residues within the RBD and residues located in other regions. Figure~\ref{fig:GraphHeatmapBar}B  demonstrates the strong associations and dependencies between these residues. In particular, there exists a significant correlation between residues originating from the RBD and residues located in the NTD, the linker domain, as well as a subset of residues within the S2 domain.

At each stage, we pick the top $10$ residues in terms of the important score, which are summarized in Figure~\ref{fig:GraphHeatmapBar}C. Our method  successfully identifies the residues previously reported in \citep{ray2021distant}, such as A570 and I569, confirming the effectiveness of our approach. Additionally, our method discovers novel residues that were not reported in \citep{ray2021distant}, primarily located in the linker and S2 domains.

\begin{figure*}[htbp]	
\centering
\begin{tabular}{ccc}
    \includegraphics[width=0.28\textwidth]{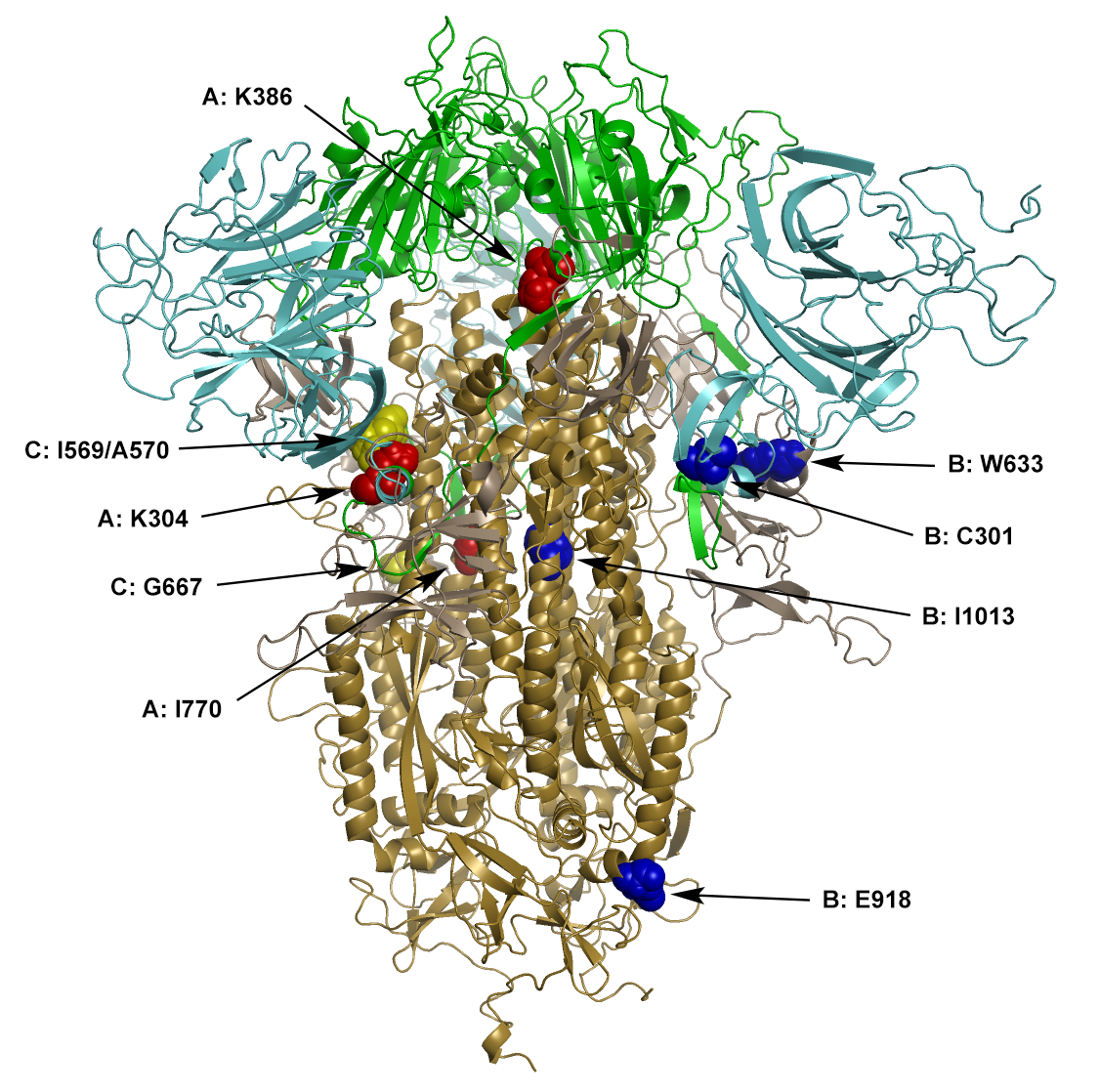}  &
 \includegraphics[width=0.28\textwidth]{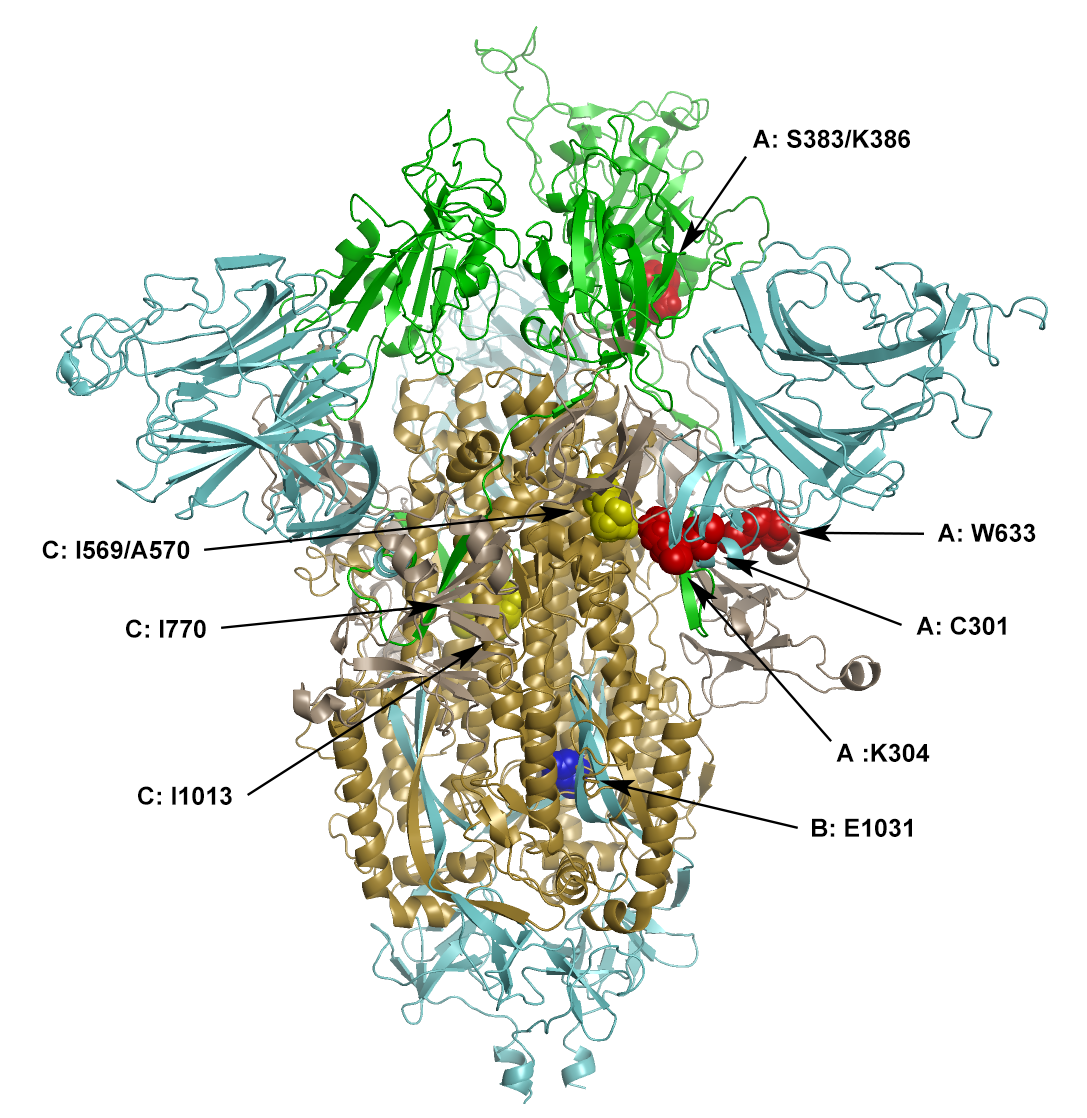} &   \includegraphics[width=0.28\textwidth]{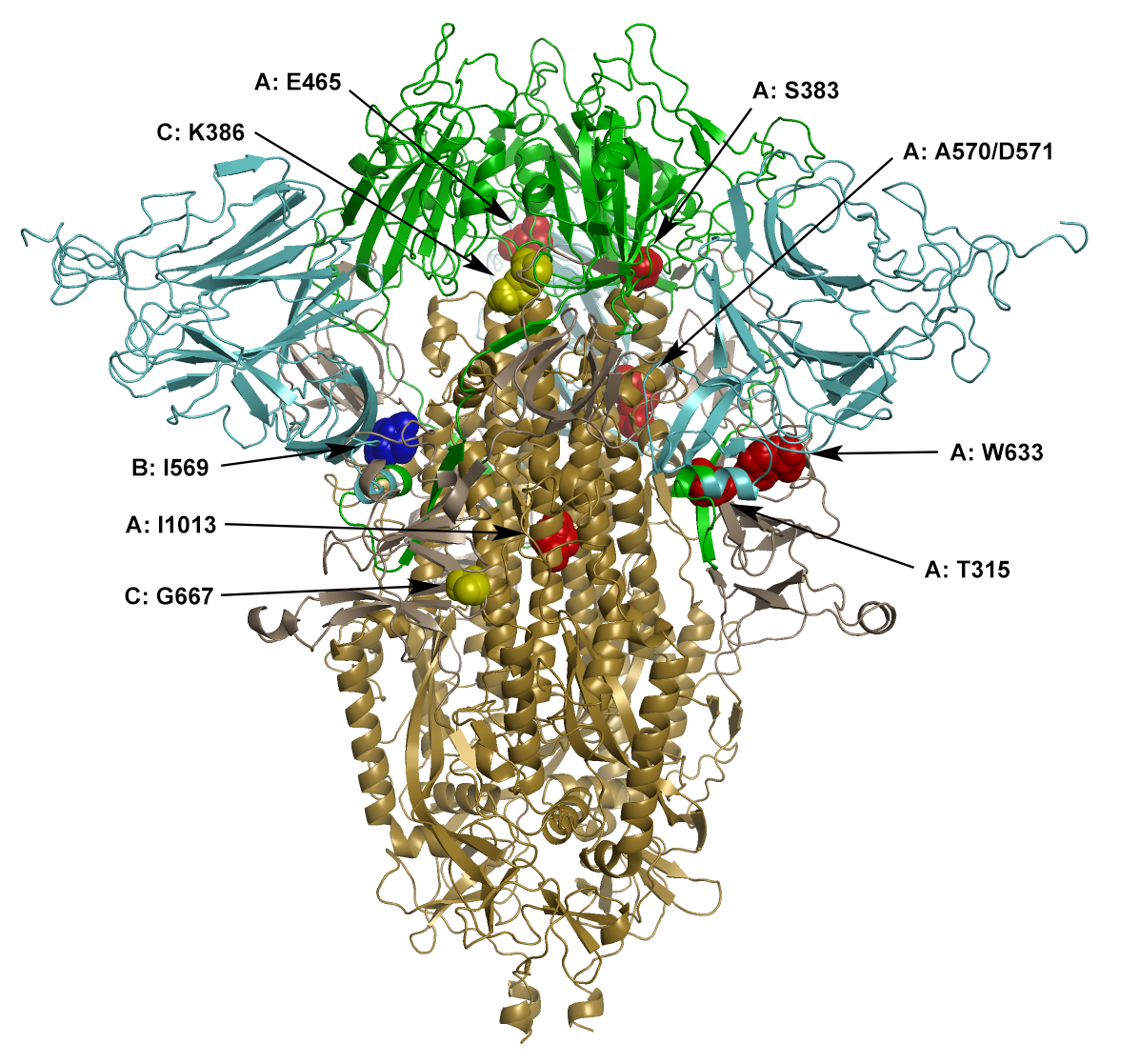} 
 \\
 { (a) Closed}& 
  { (b) Partially Open} &
  { (c) Open}
\end{tabular}
\caption{Visualization of the location of residues within the protein structure, emphasizing their proximity to the RBD.  Color coding: Chain A in red, Chain B in blue, Chain C in yellow. Domains are distinguished as follows: NTD in light blue, RBD in green, the linker domain in grey, and the S2 domain in khaki.}
\label{fig:protein}
\end{figure*}

For a more comprehensive analysis, we have illustrated the snapshots of residues K386/S683 and A570/I569 in Figure~\ref{fig:ResidueExample}. These snapshots reveal the disruption of hydrogen bonds and notable conformational changes between the three stages, emphasizing the significance of these residues in our investigation. Moreover, we have indicated the location of these residues in the protein structure, as depicted in Figure~\ref{fig:protein}, revealing their close proximity to the RBD. The spatial arrangement of these residues in the 3D structure demonstrates their potential influence on the dynamics of the RBD. This observation further strengthens the validity of our findings.

\begin{figure*}[htbp]	
\centering
\begin{tabular}{cccc} 
   \includegraphics[width=0.09\textwidth]{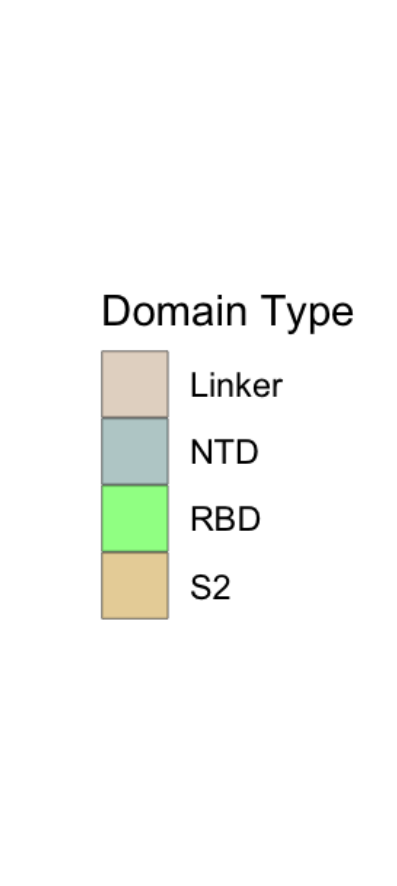} &
  \includegraphics[width=0.25\textwidth]{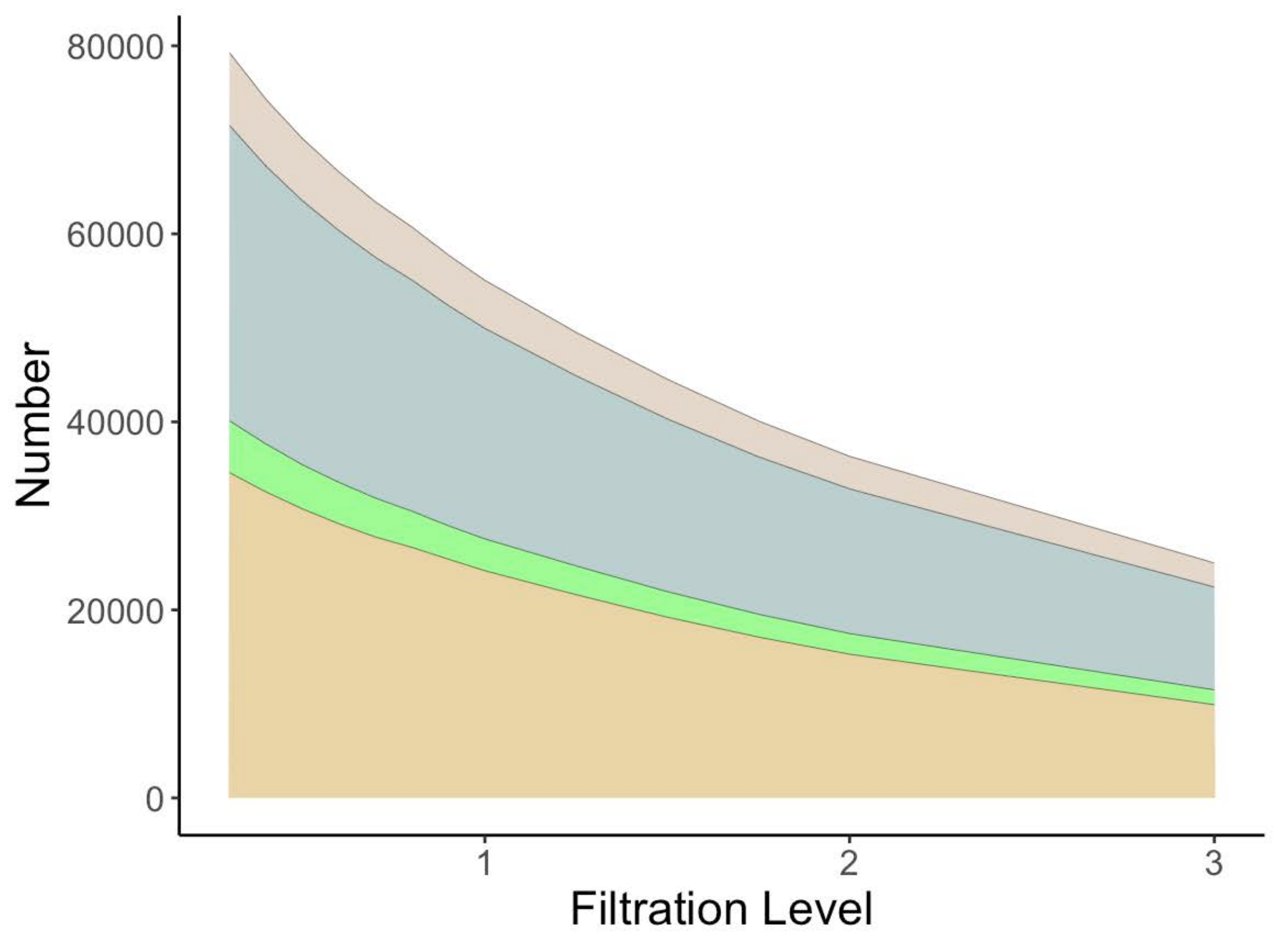}  &
 \includegraphics[width=0.25\textwidth]{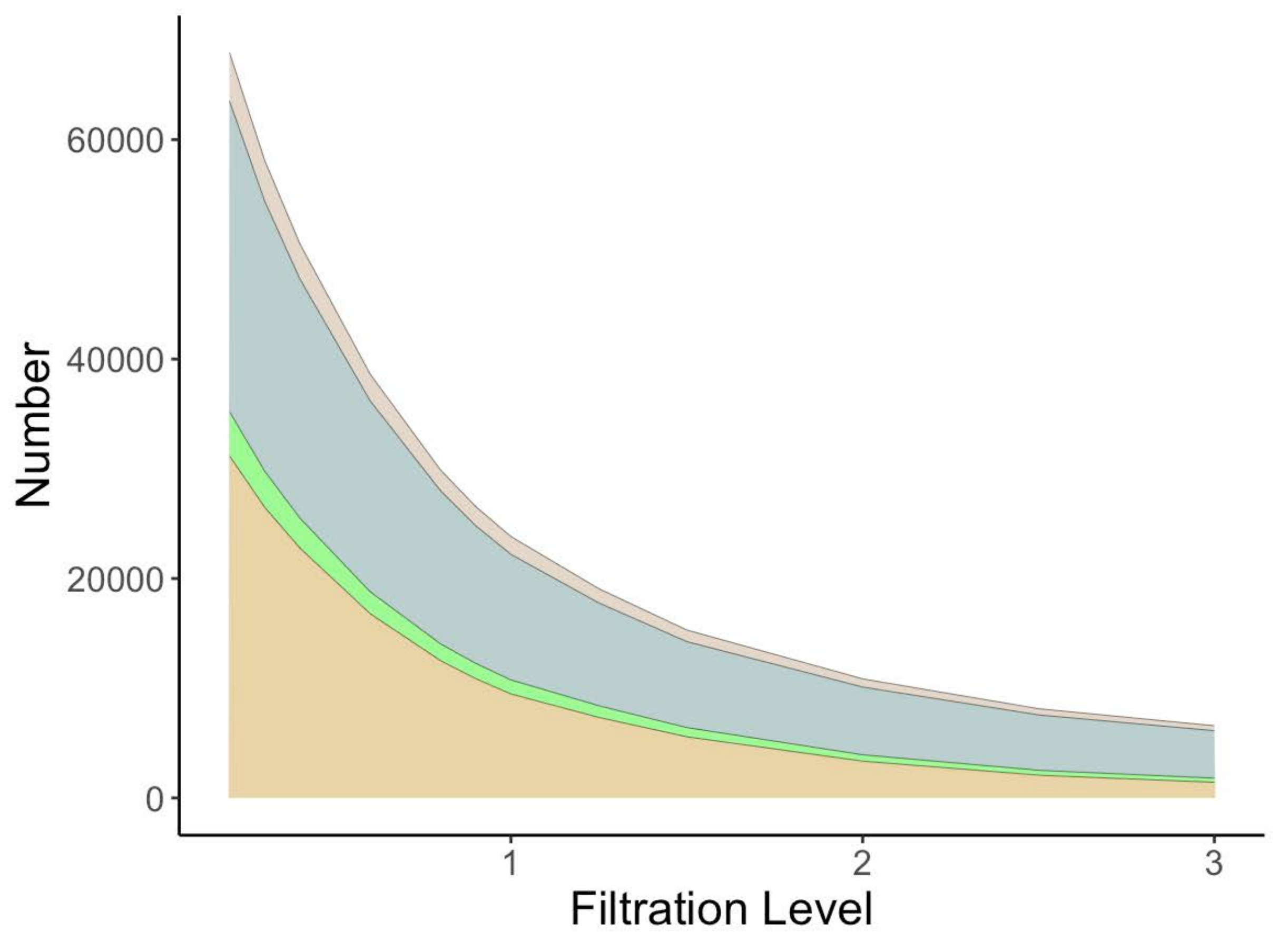} & 
  \includegraphics[width=0.25\textwidth]{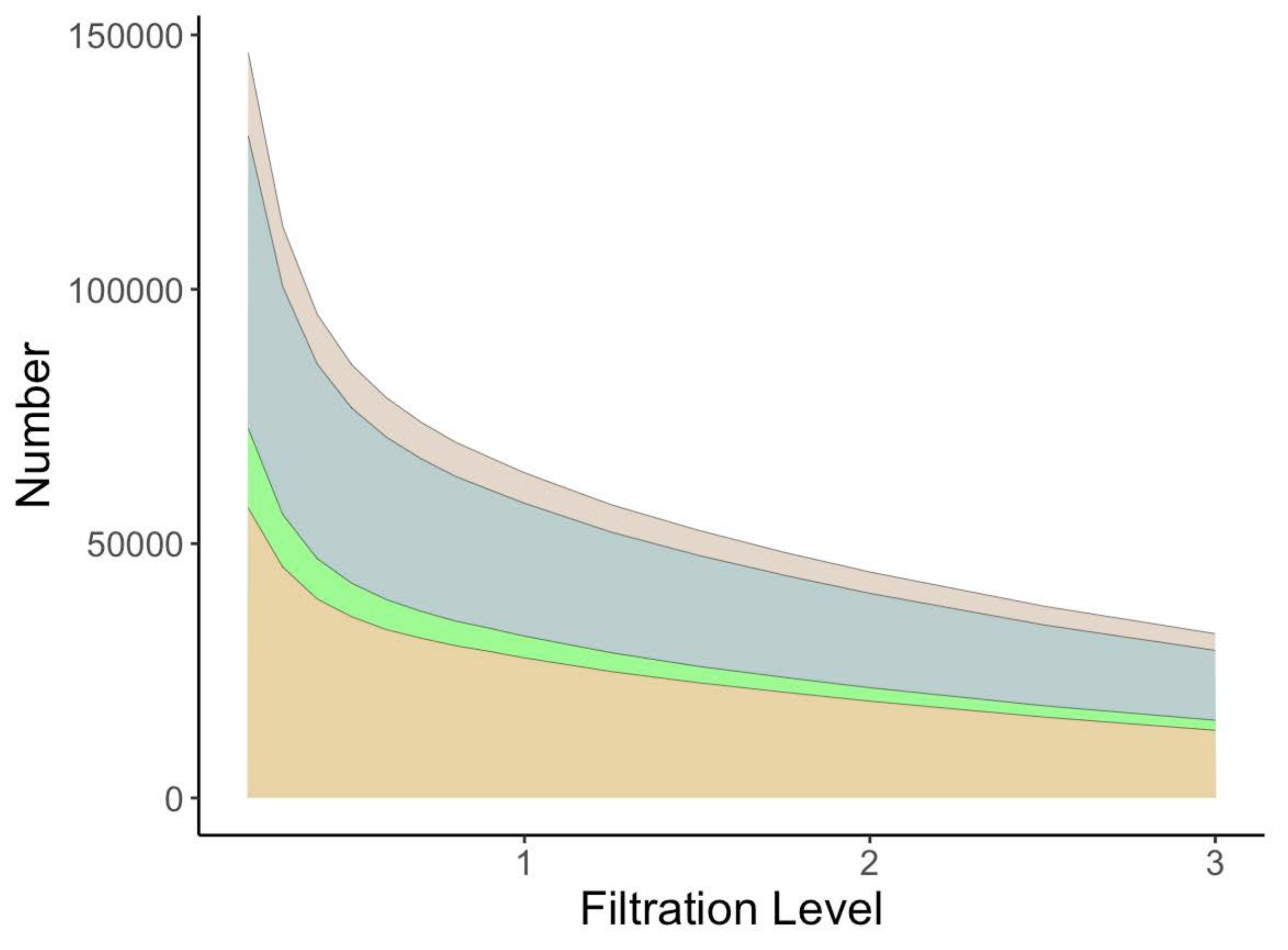} 
 \\
 &
 { (a) Closed}& 
  { (b) Partially Open} &
  { (c) Open}
\end{tabular}
\caption{The number of residues with non-zero important scores for each domain at each state.}
\label{fig:protein persis}
\end{figure*}

Finally, we present the number of residues with non-zero important scores for each domain at each state in Figure~\ref{fig:protein persis}. It shows that a substantial majority of residues are attributed to the S2 domain. However, as the filtration level gradually increases, there is a substantial reduction in residues originating from the S2 domain. In contrast, residues derived from the linker domain and RBD manifest a more gradual decrease, indicative of a higher level of persistence. These observations lend further support to our aforementioned findings.

\section{Conclusion}
This paper presents a flexible and general approach for the simultaneous detection of  graph features with false discovery rate controlled.  We also propose to select homological features from persistent homology with uniform false discovery rate controlled. Our method eliminates the need for maximal statistics and bootstrap for multiplicity control, making it simpler and faster than existing methods for simultaneous combinatorial inference.

\setlength{\bibsep}{0.85pt}{\small\bibliographystyle{ims}\bibliography{paper}}

\newpage

\begin{center}
\textit{\large Supplementary material to}
\end{center}
\begin{center}
\title{\LARGE The Wreaths of KHAN: Uniform Graph Feature Selection with False Discovery Rate Control}
\vskip10pt
\end{center}

\let\thefootnote\relax\footnotetext{
\!\!\!\!\!\! 
}

\setcounter{section}{0}
\renewcommand{\thesection}{\Alph{section}}
    
\maketitle
\begin{abstract}
  This document contains the supplementary material to the paper
  {``The Wreaths of KHAN: Uniform Graph Feature Selection with False Discovery Rate Control''.}   It is organized as follows.  
 In Section~\ref{sec:phapp}, we provide  detailed preliminaries on persistent homology. 
Section~\ref{pfsec:subgraph} proves the theoretical results of subgraph selection on FDR and power. Section~\ref{sec:pfthm:ph:generalfdppow} proves the theoretical results of on the uFDR and power of KHAN algorithm. Section~\ref{sec:ph:pfcor} proves auxiliary results in Section~\ref{GS:sec:thm}. Section~\ref{sec:B} proves  auxiliary results used in Sections \ref{pfsec:subgraph} - \ref{sec:ph:pfcor}. Section~\ref{sec:findconven} discusses the algorithm to determine homology  bases. Finally, Section~\ref{sec:simulation} presents the simulation results for synthetic data\end{abstract} 

\begin{appendix}

\section{Preliminaries of Persistent Homology on Graphical Models}\label{sec:phapp}

Persistent homology is a powerful tool to analyze the subgraph features on graphs. Comprehensive introductions to persistent homology and its foundational algebraic topology can be found in the literature \cite{horak2009persistent,aktas2019persistence}. The primary  goal of persistent homology is to measure the lifetime (the birth and death) of certain topological properties (e.g.,  
 loops, 
tetrahedron etc.,) along the evolution of the graph. As edges are added to the graph in order, subgraph features may appear or disappear. The sequence of growing (or diminishing) graphs constructed in the process is known as filtration , as illustrated  in Figure \ref{fig:filtration} (a) in the main paper. The lifespan of each homological feature with respect to the filtration levels can be represented as an interval with left and right end points being the birth and the death 
time, respectively. All intervals yielded by this filtration process can be visualized as the persistent barcodes (See  Figure \ref{fig:filtration} (b) in the main paper), which distinguish long-lasting subgraph features  from short-lived ones.



In order to understand homology group, each filtered graph or arbitrary graph  $G_0=(V,E_0)$ can be associated with its clique complex $K$, which is a collection of all the complete subgraphs of $G_0$.  By analyzing this filtered complex, we aim to understand its corresponding persistent homology groups across different filtration levels which capture the evolution of the underlying homological features. To that end, for each clique complex $K$ and 
order $k$, one needs to introduce  the associated chain groups $\mathrm{C}_k(E_0)$, the boundary operators $\partial_k:  \mathrm{C}_k(E_0) \rightarrow \mathrm{C}_{k-1}(E_0)$, the cycle groups $ Z_k(E_0)=\text{ker}~\partial_k$, 
where
$$\mbox{ker}~\partial_k=
\big\{z\in \mathrm{C}_k(E_0):\partial_k(z)=\emptyset\big\}.
$$
%
To understand these definitions, $\mathrm{C}_k(E_0)$ is a vector space  
defined on the set of $k$-dimensional oriented
simplices in $K$. The elements $c$ of the vector space $\mathrm{C}_k(E_0)$ are also denoted as $k$-chains, $c = \sum_i c_i\sigma_i$, with $\{\sigma_i\}_i$ the set of oriented $k$-simplices in $K$ and $\{c_i\in \mathbb{R}\}_i$ the particular coefficients that specify the vector $c$. The boundary operators $\partial_k$ are linear maps defined by their actions on the oriented simplices $\{\sigma_i\}_i$ in $ K$. More precisely, the image of $\partial_k$ acting on an oriented $k$-simplex $\sigma=[v_0,v_1,\ldots,v_k]$ consists of a linear combination of its oriented $k$ facets in $\mathrm{C}_{k-1}(E_0)$:
\[
\partial_{k} (\sigma)=\sum_{i}(-1)^{i}[v_{0}, v_{1}, \ldots, v_{i-1},v_i, \ldots, v_{k}].
\]
The definition of the vector subspaces cycle group $ Z_k(E_0)=\text{ker}~ \partial_k$ follows elementarily from the defination of the chain group. 

Essentially, for each  order $k$ and each graph $G_0$, $ Z_k(E_0)$ is the set of cycles which are $k$-chains with empty boundary. 
%
%
As mentioned in Section \ref{sec:setup}, Throughout this paper, the terms ``homology group'' and ``cycle group'' will be used interchangeably. We define persistent homology groups as the sequences of groups $Z_k(E(\mu))$ across various filtration levels.

\section{Proofs of Theorems on the Graph Feature Selection}\label{pfsec:subgraph}



\subsection{Proof of Theorem \ref{thm:generalfdp}}\label{sec:pfthm:generalfdp}

Recall that by definitions
\[
\text{FDP}=\frac{\sum_{j\in {\cal H}_0}\psi_j(\hat{\alpha})}{\sum_{j\in [J]}\psi_j(\hat{\alpha})}=\frac{\hat{\alpha}J}{\sum_{j\in [J]}\psi_j(\hat{\alpha})}\cdot \frac{\sum_{j\in {\cal H}_0}\psi_j(\hat{\alpha})}{J_0\hat{\alpha}}\cdot \frac{J_0}{J},
\]
where $\psi_{j}(\alpha)$ is the result of the $j$-th test at level $\alpha$ where $\psi_{j}(\alpha)=1$ if we reject $H_{0j}$ and $\psi_{j}(\alpha)=0$ otherwise. And $[J]=\{1,2,\ldots,J\}$ is the set of all hypotheses. 
The subsequent lemma articulates that $\hat \alpha$ defined in Algorithm \ref{alg:1} has an equivalent definition.
\begin{restatable}[Equivalence of $\alpha$]{lemma}{lemeqialha}\label{lem:eqialha} For the $\hat{\alpha}:=j_{\max}q/J$ yielded by the BH procedure, equivalently, we have
\[
\hat{\alpha}=\sup\Big\{\alpha>0:\frac{\alpha J}{\sum_{j\in [J]}\psi_j(\alpha)}\leq q \Big\}.
\]

\end{restatable}
Therefore, to prove Theorem~\ref{thm:generalfdp}, it suffices to show that
\begin{equation}\label{pfthm:generalfdp:1}
   \mathbb{P}\bigg( \frac{\sum_{j\in {\cal H}_0}\psi_j(\hat{\alpha})}{J_0\hat{\alpha}}\leq 1+\epsilon \bigg)=1-o(1).
\end{equation}
In order to deal with the randomness of $\hat{\alpha}$, 
we provide a probability bound on $\hat{\alpha}$ first, i.e.,
\begin{equation}\label{pfthm:generalfdp:2}
    \mathbb{P}\Big(\hat{\alpha}\in \big[q\widetilde{J}_1/J,q \big]\Big)=1-o(1),
\end{equation}
{where $\widetilde{J}_1=|\widetilde{{\cal H}}_1| $ and $\widetilde{{\cal H}}_1 $ is defined in \eqref{eqn:h1tilde}}.
To prove \eqref{pfthm:generalfdp:2}, note that by Lemma \ref{lem:eqialha} the definition of $\hat{\alpha}$  we have
\[ 
\mathbb{P}\Bigl(q\widetilde{J}_1/J \leq \hat{\alpha} \leq q\Bigr) \geq \mathbb{P}\Bigl(\frac{q\widetilde{J}_1 / J\cdot J}{\sum_{j } \psi_j(q\widetilde{J}_1 / J)} \leq q\Bigr)=\mathbb{P}\Bigl(\frac{\widetilde{J}_1}{\sum_{j } \psi_{j}( q\widetilde{J}_1 / J)} \leq 1\Bigr),
\]
where the RHS is $1-o(1)$ by \eqref{pfthm:power:1} in the proof Theorem \ref{thm:power} as  $\widetilde{{\cal H}}_1\subset [J]$.

Now, utilizing \eqref{pfthm:generalfdp:2}, we can obtain that
\begin{align*}
    \mathbb{P} \bigg(\frac{\sum_{j\in {\cal H}_0}\psi_j(\hat{\alpha})}{J_0\hat{\alpha}}\leq 1+\epsilon \bigg)\geq &\mathbb{P} \bigg(\hat{\alpha}\in [q\widetilde{J}_1/J,q],\frac{\sum_{j\in {\cal H}_0}\psi_j(\hat{\alpha})}{J_0\hat{\alpha}}\leq 1+\epsilon \bigg)\\
    \geq &\mathbb{P}\bigg(\hat{\alpha}\in [q\widetilde{J}_1/J,q],\sup_{q\widetilde{J}_1/J\leq \alpha\leq q}\frac{\sum_{j\in {\cal H}_0}\psi_j(\alpha)}{J_0\alpha}\leq 1+\epsilon \bigg)\\
    = &\mathbb{P} \Big(\hat{\alpha}\in [q\widetilde{J}_1/J,q] \Big)\cdot \mathbb{P}\bigg(\sup_{q\widetilde{J}_1/J\leq \alpha\leq q}\frac{\sum_{j\in {\cal H}_0}\psi_j(\alpha)}{J_0\alpha}\leq 1 + \epsilon \bigg).
\end{align*}
So to show \eqref{pfthm:generalfdp:1}, it suffices to show that
\begin{equation}\label{pfthm:generalfdp:3}
    \mathbb{P} \bigg(\sup_{q\widetilde{J}_1/J\leq \alpha\leq q}\frac{\sum_{j\in {\cal H}_0}\psi_j(\alpha)}{J_0\alpha}\leq 1+\epsilon \bigg)=1-o(1).
\end{equation}
To make the analysis of uniform control on an interval w.r.t., $\alpha$ tractable, we introduce \begin{equation}\label{GS:eqn:t}
    t_1=\bar{\Phi}^{-1}(q), t_m=\bar{\Phi}^{-1}(q\widetilde{J}_1/J)
\end{equation} and choose $t_i$, $i=2,3,\ldots,m-1$ such that $t_i-t_{i-1}\asymp 1/\sqrt{\log d}$. Here  $\bar{\Phi}$  is the cumulative
distribution function of standard Gaussian.   So we have
\begin{equation}\label{pfthm:generalfdp:3.1}
  \frac{\bar{\Phi}(t_i+O(1/\sqrt{\log d}))}{\bar{\Phi}(t_{i-1}+O(1/\sqrt{\log d}))}=1+o(1),\qquad \text{for all }i=1,2,\cdots,m-1.  
\end{equation}
Also, for any constant $C$ and any $\alpha \in \big[\bar{\Phi}(t_i),\bar{\Phi}(t_{i-1})\big]$, we have
\[
\frac{\sum_{j\in {\cal H}_0}\psi_j(\bar{\Phi}(t_i))}{J_0\bar{\Phi}(t_i)}\frac{\bar{\Phi}(t_i)}{\bar{\Phi}(t_{i-1})}\leq \frac{\sum_{j\in {\cal H}_0}\psi_j(\alpha)}{J_0\alpha}\leq \frac{\sum_{j\in {\cal H}_0}\psi_j(\bar{\Phi}(t_{i-1}))}{J_0\bar{\Phi}(t_{i-1})}\frac{\bar{\Phi}(t_{i-1})}{\bar{\Phi}(t_i)}.
\]
Hence, to prove \eqref{pfthm:generalfdp:3}, it suffices to show that for a constant $C$
\begin{equation}\label{pfthm:generalfdp:3.2}
   \max_{1\leq i\leq m} \frac{\sum_{j\in {\cal H}_0}\psi_j(\bar{\Phi}(t_i))}{J_0\bar{\Phi}(t_i)}\leq 1+o_P(1).
\end{equation}
We discuss two cases  $(a)$ and $(b)$ separately: note that for any $j\in[J]$,
\begin{align*}    
&(a):  \psi_j(\alpha)=1 \Longleftrightarrow \sqrt{n}|\hat{W}_e/\hat{\sigma}_e|>\bar{\Phi}^{-1}(\alpha/2), \text{ for any $e\in E(F_j)$};\\
     &(b):  \psi_j(\alpha)=1 \Longleftrightarrow  \sqrt{n}\hat{W}_e/\hat{\sigma}_e>\bar{\Phi}^{-1}(\alpha), \text{ for any $e\in E(F_j)$ }.
\end{align*}
Introduce $\{\widetilde{t}_i\}_{i=1}^m$ and $\{\bar{t}_i\}_{i=1}^m$ as follow,
\[
\widetilde{t}_i=\bigg\{\begin{array}{ll}
     \bar{\Phi}^{-1}(\bar{\Phi}(t_i)/2)& \text{if $(a)$};\\
     t_i&  \text{if $(b)$};
\end{array},\qquad \bar{t}_i:=\frac{\widetilde{t}_i}{1+C/\log d}-\frac{C}{\sqrt{\log d}}.
\]
 Now, by direct calculations, for  case $(a)$  we have with probability $1-o(1)$
\begin{align*}
  \max_i\frac{\sum_{j\in {\cal H}_0}\psi_j(\bar{\Phi}(t_i))}{J_0\bar{\Phi}(t_i)}\leq & \max_i\frac{\sum_{j\in {\cal H}_0}\mathbb{I}\big\{\sqrt{n}|\hat{W}_{e_j}/\hat{\sigma}_{e_j}|>\widetilde{t}_i\big\}}{J_02\bar{\Phi}(\widetilde{t}_i)}\\
  \leq & \max_i\frac{\sum_{j\in {\cal H}_0}\mathbb{I}\big\{\sqrt{n}|\hat{W}_{e_j}/\sigma_{e_j}|(1+C/\log d)>\widetilde{t}_i\big\}}{J_02\bar{\Phi}(\widetilde{t}_i)}\\
    \leq & \max_i\frac{\sum_{j\in {\cal H}_0}\mathbb{I}\big\{ \big(\sqrt{n}|W_{e_j}/\sigma_{e_j}|+\frac{C}{\sqrt{\log d}} \big)(1+C/\log d)>\widetilde{t}_i\big\}}{J_02\bar{\Phi}(\widetilde{t}_i)}\\
  =&\max_i\frac{\sum_{j\in {\cal H}_0}\mathbb{I}\big\{\sqrt{n}|W_{e_j}/\sigma_{e_j}|>\bar{t}_i\big\}}{J_02\bar{\Phi}(\bar{t}_i)}\frac{\bar{\Phi}(\bar{t}_i)}{\bar{\Phi}(\widetilde{t}_i)}\\
  \leq&\bigg(\max_i\frac{\sum_{j\in {\cal H}_0}\mathbb{I}\big\{\sqrt{n}W_{e_j}/\sigma_{e_j}>\bar{t}_i\big\}}{J_02\bar{\Phi}(\bar{t}_i)}+\frac{\sum_{j\in {\cal H}_0}\mathbb{I}\big\{\sqrt{n}W_{e_j}/\sigma_{e_j}<-\bar{t}_i\big\}}{J_02\bar{\Phi}(\bar{t}_i)}\bigg) \big(1+o(1)\big)
\end{align*}
 where the first inequality follows from $\psi_j(\alpha)=1 $ implies $|\sqrt{n} \hat{W}_e/\hat{\sigma}_e|>\bar{\Phi}^{-1}(\alpha/2)$ for case (a), and the second and the third inequalities follow from Assumption \ref{assum:That:condition}  with  $W^*_e= 0$ for any $e\in N_j$, $j\in {\cal H}_0$ where $N_j=\{e\in E(F_j):W^*_e=0\}$ under case $(a)$;  the last inequality is by noting that $|\bar{t}_i-\widetilde{t}_i|=O(1/\sqrt{ \log d})$ and \eqref{pfthm:generalfdp:3.1}.

For case $(b)$ where  $W^*_e = 0$ for any $e\in N_j$, $j\in {\cal H}_0$, similarly,   we have with probability $1-o(1)$
\[
\max_i\frac{\sum_{j\in {\cal H}_0}\psi_j(\bar{\Phi}\big(t_i\big))}{J_0\bar{\Phi}(t_i)}\leq  \max_i\frac{\sum_{j\in {\cal H}_0}\mathbb{I}\big\{\sqrt{n} \hat{W}_{e_j}/\hat{\sigma}_e>\widetilde{t}_i\big\}}{J_0\bar{\Phi}(\widetilde{t}_i)}\leq \max_i \frac{\sum_{j\in {\cal H}_0}\mathbb{I}\big\{ \sqrt{n}W_{e_j}/\sigma_{e_j}>\bar{t}_i\big\}}{J_0\bar{\Phi}(\bar{t}_i)}\big(1+o(1)\big).
\]

Thus,  for both $(a)$ and $(b)$, to show \eqref{pfthm:generalfdp:3.2}, it suffices to show
\[
\max_i \frac{\sum_{j\in {\cal H}_0}\mathbb{I}\{ \sqrt{n}W_{e_j}/\sigma_{e_j}>\bar{t}_i\}}{J_0\bar{\Phi}(\bar{t}_i)}\leq 1+o_P(1),\qquad \max_i \frac{\sum_{j\in {\cal H}_0}\mathbb{I}\{ \sqrt{n}W_{e_j}/\sigma_{e_j}<-\bar{t}_i\}}{J_0\bar{\Phi}(\bar{t}_i)}\leq 1+o_P(1).
\]
Further, note that $-W_e/\sigma_e=\frac{1}{n}\sum_{i=1}^n [-\xi_i(e)]$ where $-\xi_i(e)$ is mean-zero and has bounded $\psi_1$-orlicz norm, thus by symmetry and union bounds, it suffices to show that for some given choice of $\{e_j\}_j$.
\begin{equation}\label{pfthm:generalfdp:4}
    \mathbb{P}\bigg(\bigg|\frac{\sum_{j\in {\cal H}_0}\mathbb{I}\{ \sqrt{n}W_{e_j}/\sigma_{e_j}>\bar{t}_i\}}{J_0\bar{\Phi}(\bar{t}_i)}-1 \bigg|>\epsilon \bigg)=o(1/\log d), \qquad \text{for each }i=1,2,\cdots,m.
\end{equation}
Fix $i$. By Markov inequality,
\begin{equation}\label{pfthm:generalfdp:5}
    \mathbb{P}\bigg(\bigg|\frac{\sum_{j\in {\cal H}_0}\mathbb{I}\big\{ \sqrt{n}W_{e_j}/\sigma_{e_j}>\bar{t}_i\big\}}{J_0\bar{\Phi}(\bar{t}_i)}-1\bigg|>\epsilon \bigg)\leq \frac{\mathbb{E}\Big[\sum_{j\in {\cal H}_0}\big(\mathbb{I} \big\{ \sqrt{n}W_{e_j}/\sigma_{e_j}>\bar{t}_i\big\}-\bar{\Phi}(\bar{t}_i)\big)\Big]^2}{J_0^2\epsilon^2\bar{\Phi}^2(\bar{t}_i)}=\RN{1}+ \RN{2},
\end{equation}
where 
\begin{align*}
        \RN{1}=&\frac{\sum_{j_1,j_2\in {\cal H}_0}\mathbb{P}\big\{ \sqrt{n}W_{e_{j_1}}/\sigma_{e_{j_1}}>\bar{t}_i, \sqrt{n}W_{e_{j_2}}/\sigma_{e_{j_2}}>\bar{t}_i\big\}-\mathbb{P}\big\{ \sqrt{n}W_{e_{j_1}}/\sigma_{e_{j_1}}>\bar{t}_i\big\}\mathbb{P}\big\{ \sqrt{n}W_{e_{j_2}}/\sigma_{e_{j_2}}>\bar{t}_i\big\}}{J_0^2\epsilon^2\bar{\Phi}^2(\bar{t}_i)},\\
  \RN{2}=&\frac{\big(\sum_{j\in {\cal H}_0}\mathbb{P}\big\{ \sqrt{n}W_{e_j}/\sigma_{e_j}>\bar{t}_i\big\}-\bar{\Phi}(\bar{t}_i)\big)^2}{J_0^2\epsilon^2\bar{\Phi}^2(\bar{t}_i)}.
\end{align*}
For some sufficiently large  constant $C>0$, we introduce
\begin{align*}
   {\cal H}_{01}=&\Big\{(j_1,j_2)\in {\cal H}_0\times {\cal H}_0:j_1\neq j_2, \big|\text{Cov}(\xi_1(e_{j_1}),\xi_1(e_{j_2}))\big|\leq C(\log d)^{-2} (\log\log d)^{-1}\Big\},\\
   {\cal H}_{02}=&\Big\{(j_1,j_2)\in {\cal H}_0\times {\cal H}_0:j_1\neq j_2, \big|\text{Cov}(\xi_1(e_{j_1}),\xi_1(e_{j_2})) \big|> C(\log d)^{-2} (\log\log d)^{-1}\Big\}.
\end{align*}
We further decompose $\RN{1}$ as a sum of $\RN{1}_0$, $\RN{1}_1$ and  $\RN{1}_2$ where
\begin{align*}
    \RN{1}_0=&\sum_{j\in {\cal H}_0}\frac{\mathbb{P}\big\{ \sqrt{n}W_{e_j}/\sigma_{e_j}>\bar{t}_i \big\}-\mathbb{P}\big\{\sqrt{n}W_{e_j}/\sigma_{e_j}>\bar{t}_i \big\}^2}{J_0^2\epsilon^2\bar{\Phi}^2(\bar{t}_i)}\\
    \RN{1}_1=&\sum_{(j_1,j_2)\in {\cal H}_{01}}\frac{\mathbb{P}\big\{ \sqrt{n}W_{e_{j_1}}/\sigma_{e_{j_1}}>\bar{t}_i,\sqrt{n}W_{e_{j_2}}/\sigma_{e_{j_2}}>\bar{t}_i\big\}-\mathbb{P}\big\{ \sqrt{n}W_{e_{j_1}}/\sigma_{e_{j_1}}>\bar{t}_i\big\}\mathbb{P}\big\{ \sqrt{n}W_{e_{j_2}}/\sigma_{e_{j_2}}>\bar{t}_i\big\}}{J_0^2\epsilon^2\bar{\Phi}^2(\bar{t}_i)},\\
    \RN{1}_2=&\sum_{(j_1,j_2)\in {\cal H}_{02}}\frac{\mathbb{P} \big\{ \sqrt{n}W_{e_{j_1}}/\sigma_{e_{j_1}}>\bar{t}_i, \sqrt{n}W_{e_{j_2}}/\sigma_{e_{j_2}}>\bar{t}_i\big\}-\mathbb{P} \big\{ \sqrt{n}W_{e_{j_1}}/\sigma_{e_{j_1}}>\bar{t}_i\big\}\mathbb{P} \big\{ \sqrt{n}W_{e_{j_2}}/\sigma_{e_{j_2}}>\bar{t}_i\big\}}{J_0^2\epsilon^2\bar{\Phi}^2(\bar{t}_i)},\\
\end{align*}
To bound $\RN{1}_0$, $\RN{1}_1$, $\RN{1}_2$ and $\RN{2}$, by our assumptions that $\xi_1(e)$ has bounded $\psi_1$-orlicz norm and therefore has bounded moments for bounded orders,  and Lemma 6.1 of \cite{liu2013gaussian}, we have
\begin{align}
    \Big|\mathbb{P} \big( \sqrt{n}W_{e_j}/\sigma_{e_j}>\bar{t}_i \big)-\bar{\Phi}(\bar{t}_i)\Big|= & \, o\Big(\frac{\bar{\Phi}(\bar{t}_i)}{\log d}\Big), \qquad \text{ for any }j \in {\cal H}_0\label{pfthm:generalfdp:6};\\
     \Big|\mathbb{P} \big( \sqrt{n}W_{e_{j_1}}/\sigma_{e_{j_1}}>\bar{t}_i, \sqrt{n}W_{e_{j_2}}/\sigma_{e_{j_2}}>\bar{t}_i \big)-\bar{\Phi}(\bar{t}_i)^2\Big|= &\, o\Big(\frac{\bar{\Phi}^2(\bar{t}_i)}{\log d}\Big), \qquad \text{ for any }(j_1,j_2) \in {\cal H}_{01}\label{pfthm:generalfdp:7}.
\end{align}
Now by definitions and direct calculations with triangle's inequality 
\begin{align*}
    & |\RN{1}_0|\leq \frac{J_0\bar{\Phi}(\bar{t}_i)\big(1-o(\frac{1}{\log d})\big)\Big(1-\bar{\Phi}(\bar{t}_i)\big(1-o(\frac{1}{\log d})\big)\Big)}{J_0^2\epsilon^2 \bar{\Phi}^2(\bar{t}_i)}=O\Big(\frac{1}{J_0\bar{\Phi}(\bar{t}_i)}\Big),\\
    &  |\RN{1}_1|\leq \frac{J_0^2o(\frac{\bar{\Phi}^2(\bar{t}_i)}{\log d})}{J_0^2\epsilon^2 \bar{\Phi}^2(\bar{t}_i)}=o\Big(\frac{1}{\log d}\Big),\\
    & |\RN{1}_2|\leq\frac{\sum_{(j_1,j_2)\in {\cal H}_{02}} 2\mathbb{P}\big(\sqrt{n}W_{e_{j_1}}/\sigma_{e_{j_1}}>\bar{t}_i\big)}
    {J_0^2\epsilon^2 \bar{\Phi}^2(\bar{t}_i)}=O\Big(\frac{|{\cal H}_{02}|}{J_0^2\bar{\Phi}(\bar{t}_i)}\Big),\\
    & |\RN{2}|=\frac{\big(J_0\bar{\Phi}(\bar{t}_i)o(1/\log d)\big)^2}{J_0^2\epsilon^2 \bar{\Phi}^2(\bar{t}_i)}=o\Big(\frac{1}{\log^2 d}\Big).
\end{align*}
where the first and fourth inequalities hold by \eqref{pfthm:generalfdp:6};  the second inequality holds by \eqref{pfthm:generalfdp:7}; the third inequality holds by \eqref{pfthm:generalfdp:6} and  
\begin{equation}
\begin{aligned}
&  \mathbb{P} \big( \sqrt{n}W_{e_{j_1}}/\sigma_{e_{j_1}}>\bar{t}_i, \sqrt{n}W_{e_{j_2}}/\sigma_{e_{j_2}}>\bar{t}_i \big)\leq  \mathbb{P} \big( \sqrt{n}W_{e_{j_1}}/\sigma_{e_{j_1}}>\bar{t}_i \big),  \\
& \mathbb{P}\big( \sqrt{n}W_{e_{j_1}}/\sigma_{e_{j_1}}>\bar{t}_i \big)\mathbb{P} \big( \sqrt{n}W_{e_{j_2}}/\sigma_{e_{j_2}}>\bar{t}_i \big)   \leq \mathbb{P} \big( \sqrt{n}W_{e_{j_1}}/\sigma_{e_{j_1}}>\bar{t}_i \big).    
\end{aligned}    
\end{equation}

Note that there exist a selection of $e_j$ from $ N_j$, $1\leq j\leq d$ such that $|{\cal H}_{02}|=S$. Recall that  $\bar{\Phi}(\bar{t}_i)\geq q\widetilde{J}_1/J$  by \eqref{GS:eqn:t}, and our assumption on the dependence that $SJ/(J_0^2\widetilde{J}_1)=o(1/\log d)$ and sparsity that $J/(J_0\widetilde{J}_1)=o(1/\log d)$, we have
\[
\max \big\{|\RN{1}_0|, |\RN{1}_1|, |\RN{1}_2|, |\RN{2}|\big\} = o\Big(\frac{1}{\log d}\Big).
\]
Combining this with \eqref{pfthm:generalfdp:5} and \eqref{pfthm:generalfdp:4} , we finish the proof of the first claim of this theorem.

For the second claim of this theorem which is about the FDR control, note that for any $\epsilon>0$,
\begin{align*}
    \text{FDR}=\mathbb{E}[\text{FDP}]=&\int_0^1\mathbb{P}(\text{FDP}>t)d t\\
    =&\int_0^{q\frac{J_0}{J}+\epsilon}\mathbb{P}(\text{FDP}>t)d t+\int_{q\frac{J_0}{J}+\epsilon}^1\mathbb{P}(\text{FDP}>t)d t\\
    \leq & \Big(q\frac{J_0}{J}+\epsilon \Big)\cdot 1+\Big(1-q\frac{J_0}{J}-\epsilon\Big)\cdot\mathbb{P}\Big(\text{FDP}>q\frac{J_0}{J}+\epsilon\Big)\\
    =&q\frac{J_0}{J}+\epsilon+o\Big(\mathbb{P}\Big(\text{FDP}>q\frac{J_0}{J}+\epsilon\Big)\Big).
\end{align*}
By the above inequality and the first claim, we have
\[
\limsup_{n,d\rightarrow \infty}\text{FDR}\leq q\frac{J_0}{J}+\epsilon \qquad \text{for any }\epsilon >0,
\]
which finishes the proof.

\subsection{Proof of Theorem \ref{thm:power}}\label{sec:pfthm:power}

Note that $\{\psi_j\}$ hinges on random $\hat{\alpha}$, we first separate out the role of $\hat{\alpha}$ by the following
\begin{align*}
    \mathbb{P}\big(\psi_j=1,\text{ for all }j\in \widetilde{{\cal H}}_1\big)\geq&   \mathbb{P}\Bigl(\psi_j=1,\text{ for all }j\in \widetilde{{\cal H}}_1,\hat{\alpha}\in \big[q\widetilde{J}_1/J,q \big]\Bigr)\\
    \geq & \mathbb{P}\Bigl(\psi_j(\alpha)=1,\text{ for all }j\in \widetilde{{\cal H}}_1,\hat{\alpha},\alpha\in \big[q\widetilde{J}_1/J,q \big]\Bigr)\\
    =& \mathbb{P}\Bigl(\psi_j(\alpha)=1,\text{ for all }j\in \widetilde{{\cal H}}_1,\alpha\in \big[q\widetilde{J}_1/J,q\big]\Bigr)\mathbb{P}\big(\hat{\alpha}\in \big[q\widetilde{J}_1/J,q \big]\big).
\end{align*}
Therefore, to show the claim, it suffices to show for any $\alpha \in \big[q\widetilde{J}_1/J,q \big]$
\[
\mathbb{P}\big(\psi_j(\alpha)=1,\text{ for all }j\in \widetilde{{\cal H}}_1 \big)=1-o(1)\text{ and }\mathbb{P}\big(\hat{\alpha}\in \big[q\widetilde{J}_1/J,q \big]\big)=1-o(1).
\]
Note that the second equation is proved in \eqref{pfthm:generalfdp:2} in the proof of Theorem \ref{thm:generalfdp}.

Therefore, it suffices to show that  for any $\alpha \in \big[q\widetilde{J}_1/J,q \big]$,
\begin{equation}\label{pfthm:power:1}
    \mathbb{P}\big(\psi_j(\alpha)=1,\text{ for all }j\in \widetilde{{\cal H}}_1\big)=1-o(1).
\end{equation}
Consider the two scenarios 
\begin{equation}
\begin{aligned}
\text{$(a)$: }& p_e=2-2\Phi({\sqrt{n}  |\hat{W}_e}|/\hat{\sigma}_e) \\
\text{$(b)$: }& p_e=1-\Phi\big(\sqrt{n}  \hat{W}_e/\hat{\sigma}_e\big)     
\end{aligned}    
\end{equation}
separately.   
For Scenario $(a)$,  by definitions, for all $e\in E(F_j)$, $\sqrt{n}|\hat{W}_e/\hat{\sigma}_e|>\bar{\Phi}^{-1}(\alpha/2)$ implies $\psi_j(\alpha)=1$, so we have
\begin{align*}
\mathbb{P}\big(\psi_j(\alpha)=1,\text{ for all }j\in \widetilde{{\cal H}}_1\big)\geq & \mathbb{P}\Bigl(\sqrt{n}|\hat{W}_e/\hat{\sigma}_e|>\bar{\Phi}^{-1}(\alpha/2),\text{ for all }e\in E(F_j),j\in \widetilde{{\cal H}}_1\Bigr)\\
\geq & 1-\sum_{e\in \cup_{j\in \widetilde{{\cal H}}_1}E(F_j)}\mathbb{P}\Bigl(\sqrt{n}|\hat{W}_e/\hat{\sigma}_e|\leq \bar{\Phi}^{-1}(\alpha/2)\Bigr).    
\end{align*}
Therefore, to show \eqref{pfthm:power:1}, by  $|\cup_{j}E(F_j)|\leq |\bar E|\leq d^2$ and union bounds, it suffices to show that
\[
\mathbb{P}\Bigl(\sqrt{n}|\hat{W}_e/\hat{\sigma}_e|\leq \bar{\Phi}^{-1}(\alpha/2)\Bigr)=o(1/d^2),\qquad \text{ for all }e\in \cup_{j\in \widetilde{{\cal H}}_1} E(F_j).
\]
To see this, denoting $B:=\Big\{\big|\hat{W}_e-W_e-W^*_e\big|/\sigma_e\leq C/\log d,|\sigma_e/\hat{\sigma}_e-1|\leq C\sqrt{(\log d)/n} \Big\}$,   we have
\begin{align*}
    \mathbb{P}\Bigl(\big|\sqrt{n}\hat{W}_e/\hat{\sigma}_e \big|\leq \bar{\Phi}^{-1}(\alpha/2)\Bigr)
    \leq & \mathbb{P}\Bigl(\big|\sqrt{n}\hat{W}_e/\hat{\sigma}_e \big|\leq \bar{\Phi}^{-1}(\alpha/2),B\Bigr)+\mathbb{P}(B^c)\\
    \leq& \mathbb{P}\Bigl(\sqrt{n}|W_e+W^*_e|/\sigma_e\leq \bar{\Phi}^{-1}(\alpha/2)+C\frac{C}{\sqrt{\log d}}\Bigr)+o(1/d^2)\\
    \leq& \bar{\Phi}\Bigl(\sqrt{n}|W^*_e|/\sigma_e-\bar{\Phi}^{-1}(\alpha/2)- \frac{C}{\sqrt{\log d}}\Bigr)\big(1+o(1/\log d)\big)+o(1/d^2)
    \\
    =&  o(1/d^2),
\end{align*}
where the first inequality holds by $\mathbb{P}(A)=\mathbb{P}(AB)+\mathbb{P}(AB^c)\leq \mathbb{P}(AB)+\mathbb{P}(B^c)$ for any events $A,B$; the second inequality holds by event $B$ with $\log d /\log n=O(1)$; the third inequality holds by \eqref{pfthm:generalfdp:6} and the symmetry of $W$;  the last equality holds by $ \bar{\Phi}(t)= \exp(-t^2/2)/(\sqrt{2\pi}t)(1+o(1))$ for $t\gg 1$ and $|\sqrt{n}W^*_e|/\sigma_e>\sqrt{4\log d}$ for any $e\in E(F_j)$, $j\in \widetilde{{\cal H}}_1$.

For Scenario $(b)$,  by definitions, for all $e\in E(F_j)$, $\sqrt{n}\hat{W}_e/\hat{\sigma}_e>\bar{\Phi}^{-1}(\alpha)$ implies $\psi_j(\alpha)=1$, similarly, we have 
\begin{align*}
\mathbb{P}\big(\psi_j=1,\text{ for all }j\in \widetilde{{\cal H}}_1\big)\geq & \mathbb{P}\Bigl(\sqrt{n}\hat{W}_e/\hat{\sigma}_e>\bar{\Phi}^{-1}(\alpha),\text{ for all }e\in E(F_j),j\in \widetilde{{\cal H}}_1\Bigr)\\
\geq & 1-\sum_{e\in \cup_{j\in \widetilde{{\cal H}}_1}E(F_j)}\mathbb{P}\Bigl(\sqrt{n}\hat{W}_e/\hat{\sigma}_e\leq \bar{\Phi}^{-1}(\alpha)\Bigr),  
\end{align*}
and we have
\begin{align*} 
    \mathbb{P}\Bigl(\sqrt{n}\hat{W}_e\leq \bar{\Phi}^{-1}(\alpha)\Bigr)
    \leq& \bar{\Phi}\Bigl(\sqrt{n}W^*_e-\bar{\Phi}^{-1}(\alpha)- \frac{C}{\sqrt{\log d}}\Bigr)\Big(1+o\Big(\frac{1}{\log d}\Big)\Big)+o(1/d^2)
    =  o(1/d^2),
\end{align*}
where the equality is by $\bar{\Phi}(t)\sim \exp(-t^2/2)/(\sqrt{2\pi}t)$ for $t\gg 1$ and $\sqrt{n}W^*_e/\sigma_e>\sqrt{4\log d}$ for any $e\in E(F_j)$, $j\in \widetilde{{\cal H}}_1$. This finishes the whole proof.








\section{Proof of Theorem \ref{thm:ph:generalfdppow}}\label{sec:pfthm:ph:generalfdppow} 

For simplicity of notations, without loss of generality, we assume $\sigma_e = 1$, thereby simplifying $W_e / \sigma_e = W_e$.
The core argument of Theorem \ref{thm:ph:generalfdppow} hinges on the following lemma:
\begin{restatable}[Equivalence of Algorithms]{lemma}{lemeqialgph}\label{lem:eqi:algph} The persistent barcodes determined by Algorithm \ref{alg:ph:barcode} is equivalent to that determined by running Algorithm \ref{alg:ph:basis} through the whole filtration line $\mu\geq 0$ in the sense that at each filtration level $\mu$, $\hat{ Z}(\mu)$ are the same. 
\end{restatable}
Based on Lemma \ref{lem:eqi:algph}, it suffices to develop statistical controls on uFDP and uFDR for the homology group by Algorithm \ref{alg:ph:basis} on $\mu \in (\mu_0,\mu_1)$.

To prove Theorem \ref{thm:ph:generalfdppow}, we first introduce several notations. 
Based on the definition of critical edges $\{\bar{e}_j: j=1,2,\ldots \bar{J}\}$ in Section~\ref{sec:thm:ph} that $\bar{e}_j$ corresponds to an increase of rank in the homology group. So
we define $Z_j$ as the corresponding basis in the homology group. Note that essentially we are considering the following multiple testing problem 
\begin{equation}\label{eq:ph:mt:criedge}
     \bar H_{0j}(\mu):\bar{e}_j\notin E^*(\mu),\quad\text{v.s.}\quad \bar H_{1j}(\mu):\bar{e}_j\in E^*(\mu), \quad j\in [\bar{J}].
\end{equation}
Let $\bar{{\cal H}}_0(\mu)=\{j\in [\bar{J}]: \bar H_{0j}(\mu)\text{ is true}\}$ is the set of all true null hypotheses at filtration level $\mu$ and $\bar{J}_0(\mu)=|\bar{{\cal H}}_0(\mu)|$. 
We also formally define test results $ \psi_{j,\mu}(\alpha)$ at level $\alpha$, where
$\psi_{j,\mu}(\alpha) = 1$ if and only if $ Z_j \subseteq \hat{Z}(\mu)$, and $\psi_{j,\mu}(\alpha) = 0$ otherwise. For brevity, we adopt $\psi{j,\mu}$ in the subsequent proof, when no ambiguity arises.

Then \eqref{eq:ph:FDPFDR} is equivalent to define uFDP and FDP($\mu$) in the following traditional way: \[
\mathrm{uFDP}=\sup_{\mu \in [\mu_0,\mu_1]} \mathrm{FDP}(\mu), \text{ where } 
\mathrm{FDP}(\mu) = \frac{\sum_{j\in \bar{{\cal H}}_0(\mu)}\psi_{j,\mu}(\hat{\alpha}(\mu))}{\max \Bigl\{1, \sum_{j=1}^{\bar{J}} \psi_{j,\mu}(\hat{\alpha}(\mu))\Bigr\}}.
\]
 We also have $\widetilde{{\cal H}}_1(\mu)= \big\{j:  \bar{H}_{1j}(\mu)\text{ is true with }W^*_{\bar{e}_j}(\mu)>C\sqrt{\log d/n}\big\}$ and $\widetilde{J}_1(\mu)=|\widetilde{{\cal H}}_1(\mu)|$. 


However, due to the absence of prior knowledge regarding the underlying graph structure $G^*$ (which would allow us to directly identify the homology groups), it becomes challenging to ascertain which edges are critical and to distinguish between homological features in each $Z_{j-1}$. This complexity hampers our ability to devise a definitive decision rule for each hypothesis. Given a selected homological feature, the inherent complexity of the problem limits our capacity to associate it with its respective hypothesis pair, rendering the empirical values of $\{\psi_{j,\mu}\}_{j}$ unobservable. Despite these obstacles, theoretically, it is still possible to assess $\{\psi_{j,\mu}\}_{j}$  based on the selected homological features and ${Z_{j-1}}$, enabling us to examine the theoretical behavior of uFDP and uFDR.


We first have the following theorem on uFDP control.

\begin{restatable}[General guarantee on tail probability of FDP$(\mu)$  for persistent barcodes]{theorem}{thmphgeneralfdp}\label{thm:phgeneralfdp}
Under the same conditions of Theorem \ref{thm:ph:generalfdppow}, 
we have the persistent barcode  yielded by Algorithm \ref{alg:ph:basis} on a fixed filtration level $\mu \geq0$ satisfies
\[
\sup_{\mu \in [\mu_0,\mu_1]} \mathrm{FDP}(\mu)\leq \frac{q\bar{J}_0(\mu_1)}{\bar{J}}+o_P(1). 
\]


\end{restatable}

Then uFDR control follows by Theorem \ref{thm:phgeneralfdp}. Note that for any $\epsilon>0$, we have
\begin{align*}
    \text{uFDR} = &  \mathbb{E}\big[\sup_{\mu\in[\mu_0,\mu_1]}\text{FDP}(\mu)\big]=  \int_0^1\mathbb{P}\big(\sup_{\mu\in[\mu_0,\mu_1]}\text{FDP}(\mu)>t\big) \mathrm{d} t\\
    =&\int_0^{q\bar{J}_0(\mu_1)/\bar{J} +\epsilon}\mathbb{P}\big( \sup_{\mu\in[\mu_0,\mu_1]}\text{FDP}(\mu)>t\big)\mathrm{d} t
    + \int_{q\bar{J}_0(\mu_1)/\bar{J} +\epsilon}^1\mathbb{P}\big(\sup_{\mu\in[\mu_0,\mu_1]}\text{FDP}(\mu)>t\big)\mathrm{d} t\\
    \leq & (q\bar{J}_0(\mu_1)/\bar{J}+\epsilon)\cdot 1+(1-\epsilon)\cdot \mathbb{P}\big(\sup_{\mu\in[\mu_0,\mu_1]}\text{FDP}(\mu)>q\bar{J}_0(\mu_1)/\bar{J} +\epsilon \big) \\
    =&q\bar{J}_0(\mu_1)/\bar{J}+\epsilon+o(1),
\end{align*}
which gives
\[
\text{uFDR} \leq \frac{q\bar{J}_0(\mu_1)}{\bar{J}}+o(1)+\epsilon,\qquad \text{for any }\epsilon >0.
\]

Secondly for power analysis, we have the following theorem.
\begin{restatable}[General guarantee on power  for persistent barcodes]{theorem}{thmphgeneralpow}\label{thm:phgeneralpow}
Under the same conditions of Theorem \ref{thm:ph:generalfdppow},  we have the persistent barcode  yielded by Algorithm \ref{alg:ph:basis} on a fixed filtration level $\mu \geq0$ satisfies
\[
\mathbb{P}\Bigl(Z \big(\mu + C\sqrt{\log d/n} \big) \subseteq \hat{ Z}(\mu) \text{ for all }\mu\in[\mu_0,\mu_1] \Big) = 1- o(1).
\]
\end{restatable}


\subsection{Proof of Theorem \ref{thm:phgeneralfdp}}
 The structure of this proof is similar to the proof of Theorem \ref{thm:generalfdp}. 
For simplicity, we only highlight the key steps.

With $\hat{\alpha}$ defined in Algorithm \ref{alg:ph:basis}, we employ this notation $\hat{\alpha}(\mu)$ specifically when executing Algorithm \ref{alg:ph:basis} at the filtration level $\mu$.  Note that 
\[
\text{FDP}(\mu)-\frac{q\bar{J}_0(\mu_1)}{\bar{J}}=\frac{\sum_{j\in \bar{{\cal H}}_0(\mu)}\psi_{j,\mu}(\hat{\alpha}(\mu))}{\sum_{j\in [\bar{J}]}\psi_{j,\mu}(\hat{\alpha}(\mu))}-\frac{q\bar{J}_0(\mu_1)}{\bar{J}}=\Big(\frac{\hat{\alpha}(\mu)\bar{J}}{\sum_{j\in [\bar{J}]}\psi_{j,\mu}(\hat{\alpha}(\mu))}\cdot \frac{\sum_{j\in \bar{{\cal H}}_0(\mu)}\psi_{j,\mu}(\hat{\alpha}(\mu))}{\bar{J}_0(\mu_1)\hat{\alpha}(\mu)}-q\Big)\cdot \frac{\bar{J}_0(\mu_1)}{\bar{J}},
\]
where $\bar J_0(\mu)= |\bar{{\cal H}}_0(\mu)|$.
By Lemma \ref{lem:eqialha}, the BHq procedure is equivalent to have $$
\hat{\alpha}(\mu)=\sup\Big\{\alpha:\frac{\alpha \bar{J}}{\sum_{j\in [\bar{J}]} \psi_{j,\mu}(\alpha)}\leq q \Big\},$$
where $\psi_{j,\mu}(\alpha)$ is the test result at level $\alpha$.
Therefore, to prove FDP control, it suffices to show that
\begin{equation}\label{pfthm:phgeneralfdp:1}
   \mathbb{P}\bigg( \sup_{\mu\in[\mu_0,\mu_1]} \frac{\sum_{j\in \bar{{\cal H}}_0(\mu)}\psi_{j,\mu}(\hat{\alpha}(\mu))}{\bar{J}_0(\mu_1)\hat{\alpha}(\mu)}\leq 1+\epsilon \bigg)=1-o(1),\qquad \forall \epsilon >0.
\end{equation}
Similarly, to deal with the randomness of $\hat{\alpha}(\mu)$, 
 we follow the proofs of Theorem~\ref{thm:generalfdp} and Theorem \ref{thm:phgeneralpow}. We have
\begin{equation}\label{pfthm:phgeneralfdp:1.5}
    \mathbb{P}\bigg(\hat{\alpha}(\mu)\in \Big[\frac{q\widetilde{J}_1(\mu)}{\bar{J}},q \Big],\text{ for all }\mu\in[\mu_0,\mu_1]\bigg)=1-o(1),
\end{equation}
where recall that $\widetilde{J}_1(\mu)= |\widetilde {\cal H}_1(\mu)|$. 
Therefore, we have
\begin{align*}
      &\mathbb{P}\bigg( \sup_{\mu\in[\mu_0,\mu_1]} \frac{\sum_{j\in \bar{{\cal H}}_0(\mu)}\psi_{j,\mu}(\hat{\alpha}(\mu))}{\bar{J}_0(\mu_1)\hat{\alpha}(\mu)}\leq 1+\epsilon \bigg)\\
      \geq &   \mathbb{P}\bigg( \frac{\sum_{j\in \bar{{\cal H}}_0(\mu)}\psi_{j,\mu}(\hat{\alpha}(\mu))}{\bar{J}_0(\mu_1)\hat{\alpha}(\mu)}\leq 1+\epsilon ,\;\hat{\alpha}(\mu)\in \Big[\frac{q\widetilde{J}_1(\mu)}{\bar{J}},q \Big],\text{ for all }\mu \in [\mu_0,\mu_1]\bigg)\\
      \geq &   \mathbb{P}\bigg(\sup_{\substack{\alpha \in [q\widetilde{J}_1(\mu_1)/\bar{J},q]\\ \mu \in[\mu_0,\mu_1]}} \frac{\sum_{j\in \bar{{\cal H}}_0(\mu)}\psi_{j,\mu}( \alpha)}{\bar{J}_0(\mu_1)\alpha}\leq 1+\epsilon ,\;\sup_{\mu\in[\mu_0,\mu_1]} \hat{\alpha}(\mu)\in \Big[\frac{q\widetilde{J}_1(\mu)}{\bar{J}},q \Big]\bigg)\\
         \geq &   \mathbb{P}\bigg(\sup_{\substack{\alpha \in [q\widetilde{J}_1(\mu_1)/\bar{J},q]\\ \mu \in[\mu_0,\mu_1]}} \frac{\sum_{j\in \bar{{\cal H}}_0(\mu)}\psi_{j,\mu}( \alpha)}{\bar{J}_0(\mu_1)\alpha}\leq 1+\epsilon \bigg)-\mathbb{P}\bigg(\sup_{\mu\in[\mu_0,\mu_1]} \hat{\alpha}(\mu)\in \Big[\frac{q\widetilde{J}_1(\mu)}{\bar{J}},q \Big]\bigg)\\
          \geq &   \mathbb{P}\bigg(\sup_{\substack{\alpha \in [q\widetilde{J}_1(\mu_1)/\bar{J},q]\\ \mu \in[\mu_0,\mu_1]}} \frac{\sum_{j\in \bar{{\cal H}}_0(\mu)}\psi_{j,\mu}( \alpha)}{\bar{J}_0(\mu_1)\alpha}\leq 1+\epsilon \bigg)-o(1)
\end{align*}
where the first inequality is by narrowing down the event; the second inequality is because $\hat{\alpha}(\mu)\in [q\widetilde{J}_1(\mu)/\bar{J},q]$ implies $\hat{\alpha}(\mu)\subset [q\widetilde{J}_1(\mu_1)/\bar{J},q]$ for all $\mu \in[\mu_0,\mu_1]$; the last inequality is by  \eqref{pfthm:phgeneralfdp:1.5}.  

Now we see that in order to show  \eqref{pfthm:phgeneralfdp:1}, it suffices to show that for any $\mu\in[\mu_0,\mu_1]$ we have
\begin{equation}\label{pfthm:phgeneralfdp:1.6}
   \sup_{\mu\in[\mu_0,\mu_1]}\sup_{q\widetilde{J}_1(\mu_1)/\bar{J}\leq \alpha\leq q}\frac{\sum_{j\in \bar{{\cal H}}_0(\mu)}\psi_{j,\mu}(\alpha)}{\bar{J}_0(\mu_1)\alpha}\leq 1+o_P(1),
\end{equation}
where we note that the two supremums can be exchanged.

We discuss two cases  $(a)$ and $(b)$ separately: note that for any $j\in[\bar{J}]$,
\begin{equation}\label{GS:two-case}
\begin{aligned}    
    (a):&   \psi_{j,\mu}(\alpha)=1 \Longrightarrow \sqrt{n}(|\hat{W}_{\bar{e}_j}|-\mu)/\hat{\sigma}_{\bar{e}_j}>\bar{\Phi}^{-1}(\alpha/2), \text{ (e.g., Gaussian graphical models)}\\
    (b):&  \psi_{j,\mu}(\alpha)=1 \Longrightarrow \sqrt{n}(\hat{W}_{\bar{e}_j}-\mu)/\hat{\sigma}_{\bar{e}_j}>\bar{\Phi}^{-1}(\alpha),  \text{  (e.g., Ising models)}.
\end{aligned}
\end{equation}
Consider case $(b)$ first, where    $W^*_{\bar{e}_j}\leq \mu$ for $j\in \bar{{\cal H}}_0(\mu)$. Directly,  for each $j\in \bar{{\cal H}}_0(\mu)$
\[
 \psi_{j,\mu}(\alpha)\leq\mathbb{I}\{\sqrt{n}(\hat{W}_{\bar{e}_j}-\mu)/\hat{\sigma}_{\bar{e}_j}>\bar{\Phi}^{-1}(\alpha)\}\leq \mathbb{I}\{\sqrt{n}(\hat{W}_{\bar{e}_j}-W^*_{\bar{e}_j})/\hat{\sigma}_{\bar{e}_j}>\bar{\Phi}^{-1}(\alpha)\}.
\]
Combining this with $\bar{{\cal H}}_0(\mu)\subset \bar{{\cal H}}_0(\mu_1)$, it follows that
\[
\sup_{\mu\in[\mu_0,\mu_1]}\sum_{j\in \bar{{\cal H}}_0(\mu)} \psi_{j,\mu}(\alpha)\leq\sum_{j\in \bar{{\cal H}}_0(\mu_1)}  \mathbb{I}\{\sqrt{n}(\hat{W}_{\bar{e}_j}-W^*_{\bar{e}_j})/\hat{\sigma}_{\bar{e}_j}>\bar{\Phi}^{-1}(\alpha)\}.
\]
And  so to show  \eqref{pfthm:phgeneralfdp:1.6}, it suffices to show
\[
\sup_{q\widetilde{J}_1(\mu_1)/\bar{J}\leq \alpha\leq q}\frac{\sum_{j\in \bar{{\cal H}}_0(\mu_1)}  \mathbb{I}\{\sqrt{n}(\hat{W}_{\bar{e}_j}-W^*_{\bar{e}_j})/\hat{\sigma}_{\bar{e}_j}>\bar{\Phi}^{-1}(\alpha)\}}{\bar{J}_0(\mu_1)\alpha}\leq 1+o_P(1)
\]

Then, with $\{t_i\}_i$ introduced in  \eqref{pfthm:generalfdp:3.1}, by similar argument in the proof of Theorem \ref{thm:generalfdp}, for $(b)$, it suffices to show   
\begin{equation}\label{pfthm:phgeneralfdp:2}
  \text{ for $(b)$: } \max_{1\leq i\leq m}\frac{\sum_{j\in \bar{{\cal H}}_0(\mu_1)}  \mathbb{I}\{\sqrt{n}(\hat{W}_{\bar{e}_j}-W^*_{\bar{e}_j})/\hat{\sigma}_{\bar{e}_j}>t_i\}}{\bar{J}_0(\mu_1)\bar{\Phi}(t_i)} \leq 1+o_P(1),
\end{equation}
and similarly
\begin{equation}\label{pfthm:phgeneralfdp:2.5}
  \text{ for $(a)$: } \max_{1\leq i\leq m}\frac{\sum_{j\in \bar{{\cal H}}_0(\mu_1)}  \mathbb{I}\{\sqrt{n}(|\hat{W}_{\bar{e}_j}|-|W^*_{\bar{e}_j}|)/\hat{\sigma}_{\bar{e}_j}>\bar{\Phi}^{-1}(\bar{\Phi}(t_i)/2)\}}{\bar{J}_0(\mu_1)\bar{\Phi}(t_i)} \leq 1+o_P(1).
\end{equation}
To study the effect of the estimators,  we introduce $\{\widetilde{t}_i\}_{i=1}^m$ and $\{\bar{t}_i\}_{i=1}^m$ as follows:
\[
\widetilde{t}_i=\bigg\{\begin{array}{ll}
     \bar{\Phi}^{-1}(\bar{\Phi}(t_i)/2)& \text{if $(a)$};\\
     t_i&  \text{if $(b)$};
\end{array},\qquad \bar{t}_i=\frac{\widetilde{t}_i}{1+C\sqrt{(\log d)/n}}-\frac{C}{\sqrt{\log d}}.
\]
Similar to the analysis in Section \ref{sec:pfthm:generalfdp}. For case $(b)$ where  $W^*_{\bar{e}_j}\leq \mu$ for $j\in \bar{{\cal H}}_0(\mu)$,   we have with probability $1-o(1)$
 \begin{align*}
   \max_{1\leq i\leq m}\frac{\sum_{j\in \bar{{\cal H}}_0(\mu_1)}  \mathbb{I}\{\sqrt{n}(\hat{W}_{\bar{e}_j}-W^*_{\bar{e}_j})/\hat{\sigma}_{\bar{e}_j}>t_i\}}{\bar{J}_0(\mu_1)\bar{\Phi}(t_i)}=&  \max_i\frac{\sum_{j\in \bar{{\cal H}}_0(\mu)}\mathbb{I}\{\sqrt{n}(\hat{W}_{\bar{e}_j}-W^*_{\bar{e}_j})/\hat{\sigma}_{\bar{e}_j}>\widetilde{t}_i\}}{\bar{J}_0(\mu_1)\bar{\Phi}(\widetilde{t}_i)}\\
    \leq & \max_i \frac{\sum_{j\in \bar{\cal H}_0(\mu)}\mathbb{I}\{W_{\bar{e}_j}>\bar{t}_i\}}{\bar{J}_0(\mu_1)\bar{\Phi}(\widetilde{t}_i)}\\
       \leq &\max_i \frac{\sum_{j\in \bar{\cal H}_0(\mu)}\mathbb{I}\{W_{\bar{e}_j}>\bar{t}_i\}}{\bar{J}_0(\mu_1)\bar{\Phi}(\bar{t}_i)}\frac{\bar{\Phi}(\bar{t}_i)}{\bar{\Phi}(\widetilde{t}_i)},
 \end{align*}
 where second to last inequality holds by assumption on $\hat{W}$ and $\hat{\sigma}$ {    and the assumption that $\sigma_e=1$.
 }
 
On the other hand, for the case $(a)$ with that $|W^*_{\bar{e}_j}|\leq \mu$, $j\in \bar{\cal H}_0(\mu)$, similarly, we have with probability $1-o(1)$
\begin{align*}
    & \max_{1\leq i\leq m}\frac{\sum_{j\in \bar{{\cal H}}_0(\mu_1)}  \mathbb{I}\{\sqrt{n}(|\hat{W}_{\bar{e}_j}|-|W^*_{\bar{e}_j}|)/\hat{\sigma}_{\bar{e}_j}>t_i\}}{\bar{J}_0(\mu_1)\bar{\Phi}(t_i)}
  \\
  & \qquad \leq \max_i\bigg( \frac{\sum_{j\in \bar{\cal H}_0(\mu_1)}\mathbb{I}\{W_{\bar{e}_j}>\bar{t}_i\}}{\bar{J}_0(\mu_1)2\bar{\Phi}(\bar{t}_i)}+ \frac{\sum_{j\in \bar{\cal H}_0(\mu_1)}\mathbb{I}\{W_{\bar{e}_j}<-\bar{t}_i\}}{\bar{J}_0(\mu_1)2\bar{\Phi}(\bar{t}_i)} \bigg) \frac{\bar{\Phi}(\bar{t}_i)}{\bar{\Phi}(\widetilde{t}_i)}.        
\end{align*}

By $\frac{\bar{\Phi}(\bar{t}_i)}{\bar{\Phi}(\widetilde{t}_i)}=1+o(1)$ obtained in Section \ref{sec:pfthm:generalfdp}, it can be seen that  for either case $(a)$ or $(b)$, to show \eqref{pfthm:phgeneralfdp:2}, it suffices to show
\[
\max_i \frac{\sum_{j\in \bar{\cal H}_0(\mu_1)}\mathbb{I}\{W_{\bar{e}_j}>\bar{t}_i\}}{\bar{J}_0(\mu_1)\bar{\Phi}( \bar{t}_i)}\leq 1+o_P(1),
\quad 
\max_i \frac{\sum_{j\in \bar{\cal H}_0(\mu_1)}\mathbb{I}\{W_{\bar{e}_j}<- \bar{t}_i\}}{\bar{J}_0(\mu_1)\bar{\Phi}( \bar{t}_i)}\leq 1+o_P(1).
\]


The procedure to prove  the probability bounds on the RHS of the above inequalities is the same as that in the rest proof of Theorem \ref{thm:generalfdp} and so we omit the details.

\subsection{Proof of Theorem \ref{thm:phgeneralpow}}
The proof is  similar to the proof of Theorem \ref{thm:power} but they are conceptually different. Similarly, by the definition of $\psi_{j,\mu}(\alpha)$, it suffices to show that for any $\alpha \in [q\widetilde{J}_1(\mu)/\bar{J},q]$, 
\[
\mathbb{P}\big(\psi_{j,\mu}(\alpha)=1,\text{ for all } {    j \in \widetilde{\cH}_1(\mu) } ,\mu\in[\mu_0,\mu_1]\big) = 1- o(1).
\]
We consider Scenario $(b)$ in \eqref{GS:two-case} first. We have
\begin{align*}
    & \mathbb{P}\big(\psi_{j,\mu}(\alpha)=1,\text{ for all }{    j \in \widetilde{\cH}_1(\mu) },\mu\in[\mu_0,\mu_1]\big)\\
    \geq & \mathbb{P}\bigg(\sqrt{n}\hat{W}_{e}(\mu)/\hat{\sigma}_e>\bar{\Phi}^{-1}(\alpha),\text{ for all $e$ with $W^*_{e}\geq \mu+\sqrt{\frac{4 \log d}{n}}$},\mu\in[\mu_0,\mu_1]\bigg)\\
\geq &  \mathbb{P}\bigg(W_{e}+W^*_{e}-\frac{C}{\sqrt{n\log d}}-\mu>\frac{\bar{\Phi}^{-1}(\alpha)}{\sqrt{n}},\text{ for all $e$ with $W^*_{e}\geq \mu+\sqrt{\frac{4 \log d}{n}}$},\mu\in[\mu_0,\mu_1]\bigg)\\
\geq &\mathbb{P}\bigg(W_{e}+\sqrt{\frac{4 \log d}{n}}-\frac{C}{\sqrt{n\log d}}>\frac{\bar{\Phi}^{-1}(\alpha)}{\sqrt{n}},\text{ for all $e\in E$}\bigg)\\ 
\geq & 1-\sum_{e\in E}\mathbb{P}\bigg(W_{e}+\sqrt{\frac{4 \log d}{n}}-\frac{C}{\sqrt{n\log d}}<\frac{\bar{\Phi}^{-1}(\alpha)}{\sqrt{n}}\bigg),
\end{align*}
where the  second inequality is because $\hat{W}_e$ satisfies  $\mathbb{P}(\sup_{e\in E}\big|\hat{W}_e-W_e-W^*_e\big|>\frac{C}{\sqrt{n\log d}})<o(1)$ since we assume $\sigma_e=1$; the third inequality is by direct bounds; the last inequality is by union bound.

For Scenario $(a)$ in \eqref{GS:two-case},  similarly, we have 
\begin{align*}
    &\mathbb{P}\big(\psi_{j,\mu}(\alpha)=1,\text{ for all }j\in \widetilde{H}_1(\mu),\mu\in[\mu_0,\mu_1]\big)\\
    \geq & \mathbb{P}\bigg(\sqrt{n}|\hat{W}_{e}(\mu)/\hat{\sigma}_e>\bar{\Phi}^{-1}(\alpha/2),\text{ for all $e$ with $|W^*_{e}|\geq \mu+\sqrt{\frac{4 \log d}{n}}$},\mu\in[\mu_0,\mu_1]\bigg)\\
\geq & 
\mathbb{P}\bigg(|W^*_{e}|-|W_{e}|-\mu-\frac{C}{\sqrt{n\log d}}>\frac{\bar{\Phi}^{-1}(\alpha/2)}{\sqrt{n}}\text{ for all $e$ with $|W^*_{e}|\geq \mu+\sqrt{\frac{4 \log d}{n}}$},\mu\in[\mu_0,\mu_1]\bigg)\\
\geq &1-\sum_{e\in E}\mathbb{P}\bigg(|W_{e}|+\sqrt{\frac{4 \log d}{n}}-\frac{C}{\sqrt{\log d}}<\frac{\bar{\Phi}^{-1}(\alpha/2)}{\sqrt{n}},\text{ for all $e\in E$}\bigg),
\end{align*}
where the second inequality is by $\mathbb{P}\big(\sup_{e\in E}\big|\hat{W}_{e}-W_{e}-W^*_{e}\big|>\frac{C}{\sqrt{n\log d}} \big)<o(1)$; the third is by union bound.

At the same time, by  \eqref{pfthm:generalfdp:6} and elementary calculations with $\bar{\Phi}(t)\sim \exp(-t^2/2)/(\sqrt{2\pi }t)$ for $t\gg 1$
\[
\mathbb{P}\bigg(W_{e}+\sqrt{\frac{4 \log d}{n}}-\frac{C}{\sqrt{\log d}}<\bar{\Phi}^{-1}(\alpha)\bigg)=o(1/d^2),\qquad \text{for any }e\in E.
\]
It is seen that, for both Scenarios $(a)$ and $(b)$ with  any $\alpha \in [q\widetilde{J}_1(\mu)/\bar{J},q]$
\begin{equation}\label{pfthm:phgeneralpow:1}
   \mathbb{P}\big(\psi_{j,\mu}(\alpha)=1,\text{ for all }j\in \widetilde{H}_1(\mu),\mu\in[\mu_0,\mu_1]\big)\geq 1-\sum_{e\in E}o(1/d^2)=1-o(1). 
\end{equation}
This finishes the proof.



%







\section{Proofs of Propositions \ref{lem:GGM}-\ref{lem:Ising} and  Corollaries \ref{cor:ph:barcodeGGM}-\ref{cor:ph:barcodeIsing}}\label{sec:ph:pfcor}

\subsection{Proof of Proposition \ref{lem:GGM}}\label{sec:pflem:GGM}

By Lemma L.4 in \cite{neykov2019combinatorial} and that for GGM $\sigma_{uv}=\sqrt{\Theta^*_{uu}\Theta^*_{vv}}$ with $\|\Theta^*\|_1\leq M$, we directly obtain 
\[
\mathbb{P}\Bigl( \sup_{e\in \bar E}\big|\hat{W}_e-W_e-W^*_e\big|/\sigma_e\leq\frac{C}{\sqrt{n\log d}}\Bigr)=1-o(1).
\]
At the same time, adjusting the coefficient of $\log d$ in (D.36) in \cite{zhang2021startrek} can show that
\[
\sup_{e\in \bar E}|\hat{\Theta}^d_e-\Theta^*_e|=O_P\bigg(\sqrt{\frac{\log d}{n}}\bigg),
\]
which implies that
\[
\mathbb{P}\bigg( \sup_{e\in \bar E}\big|\sigma_e/\hat{\sigma}_e-1\big|\leq C\sqrt{\frac{\log d}{n}}\bigg)=1-o(1).
\]
Combining the above with union bound, we obtain the result of this proposition.



Last, by $\|X-\mathbb{E}[X]\|_{\psi_1}\leq C\|X\|_{\psi_1}$and $\|XY\|_{\psi_1}\leq \|X\|_{\psi_2} \|Y\|_{\psi_2}$, we have for any $(u,v)\in \bar E$,
\[
\|\xi_1((u,v))\|_{\psi_1}\leq C\bigg\|\frac{\Theta^*_{u\cdot}XX^\top\Theta^*_{\cdot v}}{\sqrt{\Theta^*_{uu}\Theta^*_{vv}}}\bigg\|_{\psi_1}\leq C \bigg\|\frac{\Theta^*_{u\cdot}X}{\sqrt{\Theta^*_{uu}}}\bigg\|_{\psi_2}\bigg\|\frac{\Theta^*_{\cdot v}X}{\sqrt{\Theta^*_{vv}}}\bigg\|_{\psi_2}=C.
\]
This finishes the whole proof.

\subsection{Proof of Proposition \ref{lem:Ising}}\label{sec:pflem:Ising} 


By definitions, for Ising model, we directly have
\[
\hat{W}_e-W_e-W^*_e=0.
\]
It is left to prove the claim on $\hat{\sigma}_e$.

 For simplicity of notations, we introduce, for any $(u,v)\in V\times V$, 
\[
\hat{m}_{uv}=\frac{1}{n}\sum_{i=1}^{n} X_{iu} X_{iv},\qquad m_{uv}=\mathbb{E}[X_uX_v].
\]
By definitions, for any $e\in \bar E$, we can write,
\[
\frac{\sigma_e}{\hat{\sigma}_e}-1=\frac{(\hat{m}_e-m_e)(\hat{m}_e+m_e)}{\sqrt{1-\hat{m}_e^2}(\sqrt{1-\hat{m}_e^2}+\sqrt{1-m_e^2})}.
\]
To bound the RHS, we first bound $|\hat{m}_e-m_e|$, which relies on the following lemma:
\begin{restatable}[Auxillary results for Ising models]{lemma}{lemIsingconres}\label{lem:Isingconres}For Ising model with $W^*\in \cW$ defined in \eqref{model:Ising:para}, and $ \log d/\log n=O(1)$, we have for some positive constants $C$,
\[
1-\sup_{e\in \bar E}m_e^2\geq C\log^3 d/n.
\]
\end{restatable}
This lemma allows us to utilize Bernstein inequality:
\[
|\hat{m}_e-m_e|=O_P\bigg(\sqrt{\frac{\log d}{n}(1-m_e^2)}\bigg).
\]
Further, combining this with $|m_e|,|\hat{m}_e|\leq 1$, we have
\[
\big|1-\hat{m}_e^2-(1-m_e^2)\big|=|\hat{m}_e-m_e|\cdot |\hat{m}_e+m_e|=O_P \bigg(\sqrt{\frac{\log d}{n}(1-m_e^2)} \bigg)=o_P(1-m_e^2),
\]
where the equality is because $1-m_e^2\geq C\log^3d /n$. 

It is seen that $1-\hat{m}_e^2=(1-m_e^2)(1+o(1))\geq C\log^3d /n$. And it follows that
\[
\Big|\frac{\sigma_e}{\hat{\sigma}_e}-1 \Big|=O_P\bigg(\frac{\log d}{n}\frac{1}{\sqrt{ \log^3 d/n}} \bigg)=O_P \bigg(\frac{1}{\log d}\bigg),
\]
which finishes the proof of the claim on $\hat{\sigma}$.

Last, for each $(u,v)\in \bar E$, to see $\xi_1((u,v))=\frac{X_uX_v-\mathbb{E}[X_uX_v]}{\sqrt{ \text{Var}(X_uX_v)}}$ has bounded $\psi_1$-orlicz norm. Note that
\[
\Var(X_uX_v)= 4(1-m_{uv}^2)\geq C \text{ and }\big|X_uX_v-\mathbb{E}[X_uX_v] \big|\leq 2,
\]
it can be seen that $\xi_1((u,v))$ is a bounded random variable and so has bounded $\psi_1$-orlicz norm. This finishes all the proofs of the lemma.


\subsection{Proof of Proposition \ref{prop:GGM}}\label{sec:pfprop:GGM}

For any valid dependent set $\mathcal{S}  \in \mathbb{S}$ and any $(j_1, j_2) \in \mathcal{S} $, there exist $(u,v) \in N_{j_1}$ and $(u',v') \in N_{j_2}$ such that $\big|\text{Cov}(\xi_1(u,v),\xi_1(u',v'))\big|\geq 
 C(\log d)^{-2} (\log\log d)^{-1}$,
and for the Gaussian graphical model,
we have
$$\big|\text{Cov}(\xi_1(u,v),\xi_1(u',v')) \big| = \frac{|\Theta^*_{uu'}\Theta^*_{vv'}+\Theta^*_{uv'}\Theta^*_{u'v}|}{\sqrt{\Theta^*_{uu}\Theta^*_{u'u'}\Theta^*_{vv}\Theta^*_{v'v'}}}. $$
Thus, we 
can bound $S$ by

\begin{equation}\label{pfprop:GGM:1} 
\begin{aligned}
S =&  \bigg|\Big\{({j_1},{j_2}):j_{1}, j_{2} \in {\cal H}_{0}, j_1\neq j_2,  \exists~ (u,v) \in N_{j_1},(u',v') \in N_{j_2} \text{ s.t. } \big|\text{Cov}(\xi_1(u,v),\xi_1(u',v')) \big| \geq 
\frac{C(\log d)^{-2} }{\log\log d} \Big\}\bigg|   \\
 = & \sum_{j_{1} \in {\cal H}_{0}}  
\sum_{j_{2} \in {\cal H}_{0}}  
\mathbb{I}\Big\{  \exists~ (u,v) \in N_{j_1}, (u',v') \in N_{j_2} \text{ s.t. }  \big|\text{Cov}(\xi_1(u,v),\xi_1(u',v')) \big| \geq 
 C(\log d)^{-2} (\log\log d)^{-1} \Big\}
\\
\leq & \big|\cH_0\big|   \max_{j_1 \in \cH_0, (u,v)\in N_{j_1}} 
\sum_{j_{2} \in {\cal H}_{0}} \mathbb{I}\Big\{  \exists~ (u',v') \in N_{j_2} \text{ s.t. }  \big|\text{Cov}(\xi_1(u,v),\xi_1(u',v')) \big| \geq 
 C(\log d)^{-2} (\log\log d)^{-1} \Big\}
\\
 \leq & \big|\cH_0\big|   \max_{j_1 \in \cH_0, (u,v)\in N_{j_1}}
 \sum_{j_{2} \in {\cal H}_{0}} \sum_{(u',v') \in N_{j_2}}  
 \mathbb{I}\Big\{ \big|\text{Cov}(\xi_1(u,v),\xi_1(u',v')) \big| \geq 
 C(\log d)^{-2} (\log\log d)^{-1} \Big\}\\
 = & \big|\cH_0\big|   \max_{j_1 \in \cH_0, (u,v)\in N_{j_1}}
 \sum_{(u',v') \in \bar E}   \sum_{j_{2} \in {\cal H}_{0}: (u',v')\in N_{j_2}}
 \mathbb{I}\Big\{ \big|\text{Cov}(\xi_1(u,v),\xi_1(u',v')) \big| \geq 
 C(\log d)^{-2} (\log\log d)^{-1} \Big\}\\
 = & \big|\cH_0\big|   \max_{j_1 \in \cH_0, (u,v)\in N_{j_1}}
 \sum_{(u',v') \in \bar E} \mathbb{I}\Big\{ \big|\text{Cov}(\xi_1(u,v),\xi_1(u',v')) \big| \geq 
 C(\log d)^{-2} (\log\log d)^{-1} \Big\}  \sum_{j_{2} \in {\cal H}_{0}: (u',v')\in N_{j_2}} 1
 .  
\end{aligned}
\end{equation}

So in the following proof, we assume $j_1$ and $(u,v)$ are fixed. 
Then, the edge $(u',v')$ satisfying $$\big|\text{Cov}(\xi_1(u,v),\xi_1(u',v')) \big| = \frac{|\Theta^*_{uu'}\Theta^*_{vv'}+\Theta^*_{uv'}\Theta^*_{u'v}|}{\sqrt{\Theta^*_{uu}\Theta^*_{u'u'}\Theta^*_{vv}\Theta^*_{v'v'}}} \geq 
 C(\log d)^{-2} (\log\log d)^{-1}$$ must satisfy $\Theta^*_{uu'}   \Theta^*_{vv'}+\Theta^*_{uv'} \Theta^*_{u'v} \neq 0$, which results in one of the two following cases for $u'$ and $v'$:
\[
\text{Case 1: }(u ,u') \in E^* \text{ and } (v, v') \in E^*; \qquad \text{Case 2: }(u,v') \in E^* \text{ and } (v, u') \in E^*.
\]
Since we fixed $(u,v)$ and the maximum degree of $G^*$ is $s$, we have $s$ choices of $u'$ and $s$ choices of $v'$ for any fixed $(u,v)$.

Then we fix $(u',v')$, and aim to bound $ \sum_{j_{2} \in {\cal H}_{0}: (u',v')\in N_{j_2}} 1$, which is  the number of $j_{2} \in {\cal H}_{0}$ satisfying $ (u',v')\in N_{j_2}$. We have
$$\sum_{j_{2} \in {\cal H}_{0}} \mathbb{I}\{(u',v') \in N_{j_2}\} 
\leq \sum_{j_{2} \in {\cal H}_{0}} \mathbb{I}\{(u',v') \in E(F_{j_2})\}
\leq \tbinom{d}{M-2},$$ where the first inequality holds by \eqref{GS:eq:N} the $N_j \subseteq E(F_j)$.
 Thus we finally have
$S=O(s^2d^{M-2}|{\cal H}_0|)$.





\subsection{Proof of Proposition \ref{prop:Ising}}\label{sec:pfprop:Ising}

For Ising models, 
we have
$$\big|\text{Cov}(\xi_1(u,v),\xi_1(u',v')) \big| = \Cov(X_uX_v,X_{u'}X_{v'}). $$
%
In Lemma 15 of \cite{nikolakakis2021predictive}, it is proved that for any tree structure $G = (V, E)$ and any even-sized set of nodes $V' \subseteq V$, $V'$ can be partitioned into $|V'|/2$ pairs of nodes, such that the paths connecting each pair are disjoint with each other. {So four nodes $u, u', v,v'$ can be divided into two pairs and gives two disjoint paths. We denote $\mathcal{CP}_{G}(u, u', v, v')$ as the collection of edges in the two edge-disjoint paths (so we have $\mathcal{CP}_{G}(u, u', v, v') \subseteq E$).}
Then by Theorem 10 in \cite{nikolakakis2021predictive}, if $u,v,u',v'$ are on the same  tree,
we immediately have 
\begin{equation}\nonumber
 \E(X_{u}X_{v}X_{u'}X_{v'})
 = \prod_{(i,j) \in \mathcal{CP}_{G}(u, u', v,v')} \E(X_i X_j) = \prod_{(i,j) \in \mathcal{CP}_{G}(u, u', v,v')}\tanh(w^*_{ij}),
\end{equation}
where the last equality follows from Lemma 14 in \cite{nikolakakis2021predictive}.

Further, for {forest} structure $G = \cup_\ell T_\ell$ where $T_\ell$ is the $\ell$-th tree of the forest $G$, we have 
\begin{equation}\nonumber
 \E(X_{u}X_{v}X_{u'}X_{v'})
 = \prod_{(i,j) \in \cup_\ell \mathcal{CP}_{T_\ell}(u, u', v,v')} \E(X_i X_j) = \prod_{(i,j) \in \cup_\ell \mathcal{CP}_{T_\ell}(u, u', v,v')}\tanh(w^*_{ij})
\end{equation}
from Section 4.4 in  \cite{nikolakakis2021predictive}.
In addition, in  {forest}-structu   graph models, for any node $i$ and $j$, Lemmas 13 and 14 in \cite{nikolakakis2021predictive} give
$$
 \E(X_i X_j) = \prod_{e \in \text{path}(i,j)} \tanh(w^*_e).
$$
Thus, we have
\begin{small}
\begin{equation}\label{eqn:E-E}
\begin{aligned}
\Cov(X_uX_v,X_{u'}X_{v'}) 
= & \, \E(X_{u}X_{v}X_{u'}X_{v'}) - \E(X_{u}X_{v} )\E( X_{u'}X_{v'}) \\
= &\, \prod_{(i,j) \in \cup_\ell\mathcal{CP}_{T_\ell}(u, u', v,v')}\tanh(w^*_{ij}) - \prod_{e \in \text{path}(u,v)} \tanh(w^*_e)
\prod_{e \in \text{path}(u',v')} \tanh(w^*_e).
\end{aligned}    
\end{equation}
\end{small}
Now, we are ready to give a bound for $S$. Similar to \eqref{pfprop:GGM:1}, $S$ can be bounded by
\begin{small}
\begin{equation}\nonumber
\begin{aligned}
S 
\leq &   \big|\cH_0\big|   \max_{j_1 \in \cH_0, (u,v)\in N_{j_1}}
 \sum_{(u',v') \in \bar E} \mathbb{I}\Big\{ \big|\text{Cov}(\xi_1(u,v),\xi_1(u',v')) \big|  \geq C(\log d)^{-2} (\log\log d)^{-1} \Big\}  \sum_{j_{2} \in {\cal H}_{0}: (u',v')\in N_{j_2}} 1 \\
\leq &   \big|\cH_0\big|   \max_{j_1 \in \cH_0, (u,v)\in N_{j_1}}
 \sum_{(u',v') \in \bar E} \mathbb{I}\Big\{ \big|\text{Cov}(\xi_1(u,v),\xi_1(u',v')) \big| \neq 0 \Big\}  \sum_{j_{2} \in {\cal H}_{0}: (u',v')\in N_{j_2}} 1 \\
 = &  \big|\cH_0\big|   \max_{j_1 \in \cH_0, (u,v)\in N_{j_1}}
 \sum_{(u',v') \in \bar E} 
 \mathbb{I}  \Big\{  \E(X_{u}X_{v}X_{u'}X_{v'}) - \E(X_{u}X_{v} )\E( X_{u'}X_{v'}) \neq 0 \Big\}   \sum_{j_{2} \in {\cal H}_{0}} 
 \mathbb{I}\{(u',v') \in N_{j_2}\}. 
\end{aligned}   
\end{equation}\end{small}
By \eqref{eqn:E-E}, for fixed $(u,v)$,
 the edge $(u',v')$ satisfying $\E(X_{u}X_{v}X_{u'}X_{v'}) - \E(X_{u}X_{v} )\E( X_{u'}X_{v'}) \neq 0$  must satisfy one of the two following cases:

(i) $u, u', v,v'$ are in two different trees with two nodes in each tree, and $u,v$ must be in the same tree, and $u', v'$ must be in the other tree.
So we have $O(\frac{d}{s})$ potential choices of the second tree for $u', v'$, and $O(s^2)$ choices of $u', v'$.

(ii) $u, u', v,v'$ are in the same tree, and $\text{path}(u,v) \cap \text{path}(u',v') \neq \emptyset$. We have $ \cup_\ell\mathcal{CP}_{T_\ell}(u, u', v,v') = \text{path}(u,u') \cup \text{path}(v,v')$ and thus $\E(X_{u}X_{v}X_{u'}X_{v'}) - \E(X_{u}X_{v} )\E( X_{u'}X_{v'}) \neq 0$. 
So we have  $O(s^2)$ choices of $u', v'$.
    
At the same time, when fixing $(u',v')$,  we have  $$\sum_{j_{2} \in {\cal H}_{0}} \mathbb{I}\{(u',v') \in N_{j_2}\} 
\leq \sum_{j_{2} \in {\cal H}_{0}} \mathbb{I}\{(u',v') \in E(F_{j_2})\}
\leq \tbinom{d}{M-2}.$$
 Thus we finally have
\[ S = O\Big( |\cH_0| \Big(\frac{d}{s} s^2 +s^2\Big) \tbinom{d}{M-2} \Big)
= O\big(|\cH_0|d^{M-1} s \big).
\]

\subsection{Proof of Proposition \ref{prop:GGM-Ising}   }\label{sec:disc:assm:phspadepend}
For the Gaussian graphical model, similar to the proof in Section~\ref{sec:pfprop:GGM}, we have under persistent homology
\begin{small}

\begin{align*}
   \bar S (\mu) &=\Big|\Bigl\{
 (j_1,j_2):  j_1 \neq j_2, \bar{e}_{j_1},\bar{e}_{j_2}\notin E^*(\mu),  
 \big|\text{Cov}\big( \xi_1(\bar{e}_{j_1}),\xi_1(\bar{e}_{j_2}) \big) \big|  \geq 
\frac{C}{(\log d)^{2}\log (\log d)} \Bigr\}\Big| \\
    \le &  \big|\big\{({j_1},{j_2}): j_1\neq j_2, (u,v)=\bar{e}_{j_1}\notin E^*,(u',v')=\bar{e}_{j_2}\notin E^*, \big|\text{Cov}(\xi_1(u,v),\xi_1(u',v')) \big|\neq 0\big\}\big|   \\
\leq & \sum_{j_{1}}\mathbb{I}\{(u,v)=\bar{e}_{j_1}\notin E^*\} 
 \max_{j_1 ,(u,v)}\sum_{(u',v') =\bar{e}_{j_2}\notin E^*}  \mathbb{I}\Big\{ \big|\text{Cov}(\xi_1(u,v),\xi_1(u',v')) \big|\big|\neq 0\Big\}  \sum_{j_{2} : (u',v') = \bar e_{j_2}\notin E^*} 1 \\
 \leq & \sum_{j_{1}}\mathbb{I}\{(u,v)=\bar{e}_{j_1}\notin E^*\} 
 \max_{j_1 ,(u,v)}\sum_{(u',v') =\bar{e}_{j_2}\notin E^*}  \mathbb{I}\Big\{ \Theta^*_{uu'}   \Theta^*_{vv'}+\Theta^*_{uv'} \Theta^*_{u'v} \neq 0 \Big\}  \sum_{j_{2} : (u',v') = \bar e_{j_2}\notin E^*} 1.  
\end{align*}
    
\end{small}
First we have $\sum_{j_{1}}\mathbb{I}\{(u,v)=\bar{e}_{j_1}\notin E^*\} = |\bar{{\cal H}}_0(\mu)|$ because only if $j_1 \in \bar{{\cal H}}_0(\mu)$ then $\bar{e}_{j_1}\notin E^*$.

Then we assume $j_1$ and $(u,v)$ are fixed. Similar to the proof in Section~\ref{sec:pfprop:GGM},  to satisfy $\Theta^*_{uu'}   \Theta^*_{vv'}+\Theta^*_{uv'} \Theta^*_{u'v} \neq 0$, we have $O( s(\mu)^2)$ choices of  critical edge $(u',v')$, where $s(\mu)$ is the maximum degree of $G^*(\mu)$.

Next we fix $(u',v')$, and aim to bound $\sum_{j_{2} : (u',v') = \bar e_{j_2}\notin E^*} 1$, which is  the number of $j_{2}$ having critical edge $ \bar e_{j_2} = (u',v')$. Recall that we defined $R = \max_{i=1 }^{|\bar{E}|}\ell_i$ Section \ref{sec:thm:ph}, so one critical edge can increase at most $R$ rank, i.e., correspond to at most $R$ hypotheses. So  $\sum_{j_{2} : (u',v') = \bar e_{j_2}\notin E^*} 1 = O(R)
$.

Thus, for the Gaussian graphical model,  we finally have
$$ \bar S(\mu)=O(|\bar{{\cal H}}_0(\mu)|\cdot s(\mu)^2 R).$$

For the Ising model, similar to the argument in Section~\ref{sec:pfprop:Ising},  we   have
\[  \bar S (\mu) 
= O\big(|\bar{{\cal H}}_0(\mu)| \cdot s(\mu) dR \big).
\]


\subsection{Proof of Corollary \ref{cor:FDRPOW:GGM}}\label{sec:cor:FDRPOW:GGM}

By Proposition  \ref{lem:GGM} for GGM, the three claims of this corollary  follows from Theorem \ref{thm:generalfdp}-\ref{thm:power} if  we set for any   $j\in {\cal H}_0$ and $e=(u,v)\in \bar E$
\[
 N_j=\big\{e\in E(F_j):\Theta^*_e=0\big\},\quad W_{uv}=\sum_{i=1}^n{\Theta^*_u}^\top(X_iX_i^\top\Theta^*_v-e_v)/n,\quad W^*_e= \Theta^*_e.
\]
It remains to check the sparsity and dependence conditions are matched. To see this, by Isserlis’ theorem (Theorem 1.1 in \cite{michalowicz2009isserlis}), and so
\begin{align*}
 \text{Cov}\big({\Theta^*_u}^\top (XX^\top \Theta^*_v-e_v)/\sqrt{\Theta^*_{uu}\Theta^*_{vv}},{\Theta^*_{u'}}^\top (XX^\top \Theta^*_{v'}-e_{v'})/\sqrt{\Theta^*_{u'u'}\Theta^*_{v'v'}}\big) =\frac{W^*_{uu'}W^*_{vv'}+W^*_{uv'}W^*_{u'v}}{\sqrt{W^*_{uu}W^*_{vv}}}.   
\end{align*}
This finishes the proof.

\subsection{Proof of Corollary \ref{cor:FDRPOW:Ising}} \label{sec:pfcor:FDRPOW:Ising}
Recall $m$ and $\hat{m}$ are introduced in Section~\ref{sec:pflem:Ising}.
By  Proposition~\ref{lem:Ising} for Ising model, the three claims of this corollary  follows from Theorems \ref{thm:generalfdp}-\ref{thm:power} if  we set  
\[
 N_j=
    \big\{e\in E(F_j):W^*_e\leq 0 \big\},\quad W_{uv}=\frac{1}{n}\sum_{i=1}^nX_{iu}X_{iv}-m_{uv}
\]
and
\[
W^*_{uv}=m_{uv}-\tanh(\theta)
\]
for any $j\in {\cal H}_0$, $(u,v)\in \bar E$.
The covariance calculations for Ising model are straightforward and so we omit it.  This finishes the proof.

\subsection{Proofs of Corollaries \ref{cor:ph:barcodeGGM}-\ref{cor:ph:barcodeIsing}}\label{sec:coro:ph}

The proofs of Corollaries \ref{cor:ph:barcodeGGM}-\ref{cor:ph:barcodeIsing} follow by Propositions \ref{lem:GGM}-\ref{lem:Ising} and Theorem \ref{thm:ph:generalfdppow}. The  procedures are the same as  that of Corollaries \ref{cor:FDRPOW:GGM}-\ref{cor:FDRPOW:Ising}. We omit the details.

\section{Proofs of Auxiliary Lemmas}\label{sec:B}





\subsection{Auxiliary Proposition for  Persistent Homology}


The following proposition gives a explicit formula for the dimension of the homology group of a complete graph.
\begin{restatable}{proposition}{propdim}\label{prop:dim} The dimension of the $k$-homology group of a $d$-clique is $(d-k)(d-k-1)/2$.
\end{restatable}
   By the above property, it can be seen that the dimension of each $k$-homology group is $O(d^2)$. In the ordinary subgraph selection, the number of testing triangles in a graph requires $\tbinom{d}{3}$ hypotheses. This suggests that persistent homology is a more reasonable way to select hole-structure graph features as it rules out redundancy.

\begin{proof}
For a $k$-th homology group of a $d$-clique, its dimension is equivalent to all ways to add the edges to the graph along the filtration. Denote $A_d^k=\text{dim}(Z_k(E_{d\times d}))$. We are going to develop a recurrence relation between $A_i^k$ and $A_{i-1}^k$. Note that for a $(i-1)$-clique, if we add a extra node and connect the node with the original $i-1$ nodes. The process is equivalent to the filtration process to complementing a $(i-1)$-clique with an isolated node to a $i$-clique. Noting that for $k$-th homology group, the first extra $k$-clique appears after the $k-1$ edges connecting the isolated node and the original $i-1$ nodes. After that each new edge induces a new $k$-th clique. There are $i-1$ less edges for a $(i-1)$-clique compared to a $i$-clique. Therefore, $Z_k(E_{i\times i})$ has extra $(i-1)-(k+1-1)$ dimension than $Z_k(E_{(i-1)\times (i-1)})$:
\[
A^k_i=A^k_{i-1}+i-k-1.
\]
Solving this recursive relation with $A^k_{k}=1$, we get
\[
A^k_d=(d-k)(d-k-1)/2,
\]
which finishes this proof.
\end{proof}


\subsection{Proof of Lemma \ref{lem:eqialha}}

The core of proving equivalence is to show 
\begin{equation}\label{GS:eqn:equ}
\sup\Big\{\alpha>0:\frac{\alpha J}{\sum_{j\in [J]}\psi_j(\alpha)}\leq q\Big\}\geq \frac{qj_{\max}}{J} 
\text{ and }  
\sup\Big\{\alpha>0:\frac{\alpha J}{\sum_{j\in [J]}\psi_j(\alpha)}\geq q\Big\}\geq \frac{qj_{\max}}{J}.
\end{equation}

First, noting that $j_{\max}=\sum_{j\in [J]}\psi_j(\frac{qj_{\max}}{J})$, we have
\[
\frac{\alpha J}{\sum_{j\in [J]}\psi_j(\alpha)}\bigg|_{\alpha=qj_{\max}/J}=q,
\]
which implies that $$\sup\Big\{\alpha>0:\frac{\alpha J}{\sum_{j\in [J]}\psi_j(\alpha)}\leq q \Big\}\geq \frac{qj_{\max}}{J}.$$

Next, we prove the second part of \eqref{GS:eqn:equ} by contradiction. 
Let us presume  \begin{equation}\label{GS:eqn:con}
    \widetilde{\alpha}=\sup\Big\{\alpha>0:\frac{\alpha J}{\sum_{j\in [J]}\psi_j(\alpha)}\leq q \Big\}>\frac{qj_{\max}}{J}.
\end{equation} By definitions,
\[
\widetilde{\alpha}< p_{(j_{\max}+1)}, \text{ the ($j_{\max}+1$)-th smallest $p$-value};
\]
otherwise $\frac{\alpha J}{\sum_{j\in [J]}\psi_j(\alpha)}\big|_{\alpha=p_{j_{(\max}+1)}}\leq q$ suggests that rejected hypotheses are at least $j_{\max}+1$ and lead to contradiction. By monotonicity of $\psi_j(\alpha)$ and definitions of $\widetilde{\alpha}$,
\[
\widetilde{\alpha} \leq \frac{q \sum_{j\in [J]}\psi_j(
\widetilde{\alpha})}{J}
<\frac{q \sum_{j\in [J]}\psi_j( p_{(j_{\max}+1}))}{J}=\frac{q (j_{\max}+1)}{J}.
\]
It remains to demonstrate the impossibility of $\frac{q (j_{\max})}{J}< \widetilde{\alpha}<\frac{q (j_{\max}+1)}{J}$. By \eqref{GS:eqn:con}, we have
\[
\frac{\alpha J}{\sum_{j\in [J]}\psi_j(\alpha)}\bigg|_{\alpha=\widetilde{\alpha}}>\frac{\frac{q (j_{\max})}{J} J}{\sum_{j\in [J]}\psi_j(\widetilde{\alpha})}=\frac{\frac{q (j_{\max})}{J} J}{\sum_{j\in [J]}\psi_j(\frac{q (j_{\max})}{J})}=q,
\]
which contradicts the definition of $\widetilde{\alpha}$ and finishes the proof.

\subsection{Proof of Lemma \ref{lem:eqi:algph}}

We aim to show that Algorithm \ref{alg:ph:barcode} and Algorithm \ref{alg:ph:basis} give the same homology group $\hat{Z}(\mu)$  for any filtration level $\mu$.
By definition, the edge sets constructing the homology group are identical for the two algorithms at $\mu=0$ and $\mu>\mu^{(t)}$. 
Recall that in Algorithm \ref{alg:ph:barcode}, we  obtain the homology group $\hat{Z}^{(t+1)} = \hat Z(\mu^{(t+1)})$.  Therefore, it remains to prove the following claims:

(i) In Algorithm \ref{alg:ph:basis}, the set of selected homology group  is right-continuous i.e., $\hat{Z}(\mu)$ stays the same for $\mu\in [\mu^{(s)},\mu^{(s+1)})$, for each $s=0,1,\ldots,t-1$;

    
To prove (i), recall that according to Algorithm \ref{alg:ph:barcode}, at iteration $s$, $\hat Z^{(s)}$ is the homology group  from the previous iteration. 
By the definitions of $\mu^{(s+1)}=\min_{e\in E^{(s)}}\text{lower bound of weight}_e(\mathrm{rank}(\hat Z^{(s)})q/\bar{J})$. Therefore, by the nature of BHq procedure, the selected homology group  between $[\mu^{(s)},\mu^{(s+1)})$ stays the same and we finishes the proof.




\subsection{Proof of Lemma \ref{lem:Isingconres}}

Consider the claim on $m_e$ first. By similar idea of the proof of Lemma C.3 in \cite{neykov2019property}, we can show that for any $e\in \bar E$
\begin{equation}\label{pflem:Isingconres:1}
    |m_e|\leq \tanh(w^*_e)+\sum_{m\geq 2}s(s-1)^{m-2}\tanh^m(\Theta)\leq \tanh(w^*_e)+\frac{s\tanh^2(\Theta)}{1-(s-1)\tanh(\Theta)}.
\end{equation}
By the conditions $\|w^*\|_{\infty}\leq \Theta$,  $s\tanh(\Theta)<\rho$ and $s\geq 2\rho/(1-\rho)$, we have 
\[
\tanh(w^*_e)+\frac{s\tanh^2(\Theta)}{1-(s-1)\tanh(\Theta)}\leq \tanh(\Theta)+\tanh(\Theta)\frac{\rho}{1-\rho} \leq \frac{\rho}{s(1-\rho)}\leq \frac{1}{2}.
\]
It follows that $1-m_e^2\geq 3/4\gg C\log^3d /n$. 
This finishes the proof.

\section{Algorithm for Determining the Importance Scores}\label{sec:findconven}

When searching for homology, edges are added sequentially and loops are searched for according to two principles. The first principle is to find the loop that contains the most nodes from the distance data. If multiple loops have the same number of nodes from the distance features, then search for the loop with the largest average edge weight. Record the loops found in this way.

\begin{algorithm}[H]
\caption{Homology Search Algorithm}
\label{alg:hs}

    \SetKwRepeat{Do}{do}{while}
    
    Input: Trajectory Data $X_\phi$ and $X_\psi$, Distance Data $X_d$, FDR level $q$, Filtration Level $\mu$. 

    Build a d-dimensional Gaussian graphical model based on $X_\phi$, $X_\psi$ and $X_d$.

    \For{$\mu_i \in \mu$}{
    
    Get Edge Set $E^{(i)}=\hat{E}(\mu_i)$ applying DGS\;

    Get estimated edge weights $\{\widehat{W}_e\}_{e \in E^{(i)}}$
    
    \For{$e \in E^{(i)}$}{

    Find all the loops linking the two nodes of $e$;

    Select the loop that has the most nodes from $X_d$;

    \If{Any two loops have the same number of nodes from the distance data}
    {
        Select the loop with the largest average edge weight according to $\{\widehat{W}_e\}_{e \in E^{(i)}}$
    }
    }
    }
    Output: Homologies Set $H(\mu)$
\end{algorithm}

\section{Numerical Experiments for Synthetic Data}\label{sec:simulation}

\begin{figure}[ht]
\centering 
\begin{tabular}{cc}
     \includegraphics[width=0.4\textwidth]{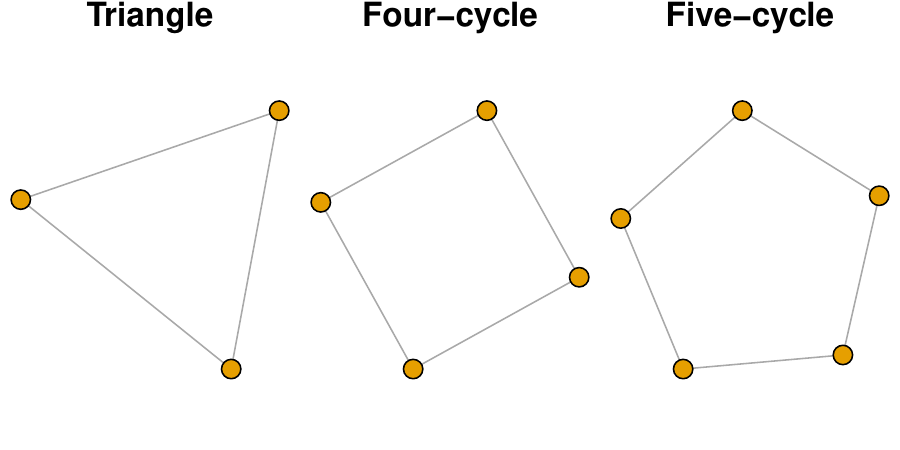} &
 \includegraphics[width=0.4\textwidth]{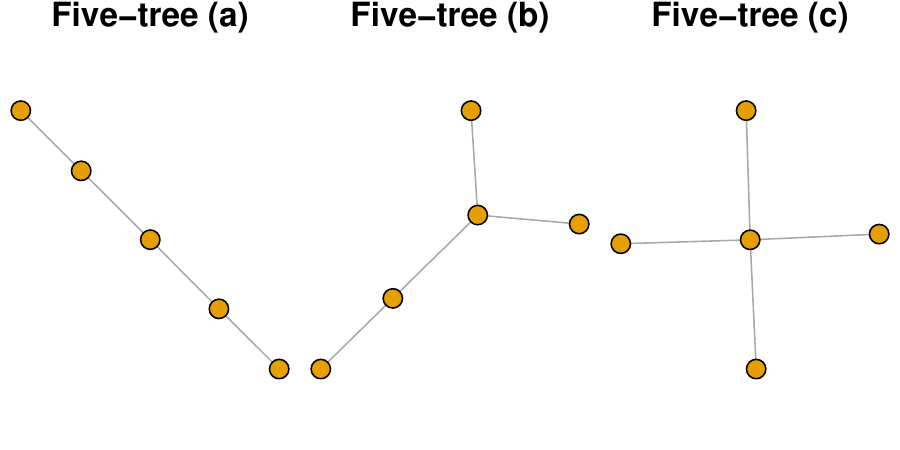}  \\
 (a) Sug-graphs for the Gaussian model & (b) Sug-graphs for the Ising model
\end{tabular}
\caption{Sub-graphs to test for the Gaussian graphical model and the Ising model. }
\label{fig:sim0}
\end{figure}

In this section, we examine the performance our method through numerical simulations. We first consider the general subgraph selection. 
For the Gaussian graphical model and the Ising model in  Sections \ref{GS:sec:GGM} and \ref{GS:sec:Ising}, in particular, we focus on the performance of the proposed method as we vary the dimension $d$, sample size $n \in \{300,350,400\}$, the nominal level $q \in \{0.05,0.1\}$ and the types of sub-graphs. For illustration, in the Gaussian graphical model, we focus on the detection of triangles, four-cycles and five-cycles, while in the Ising model, we focus on the detection of five-trees, which are trees consisting of five nodes with three possible shapes, as illustrated in Figure~\ref{fig:sim0}.  We estimate the empirical FDR and power based on $100$ repetitions.

\begin{figure}[ht]
\centering 
\begin{tabular}{cc}
     \includegraphics[width=0.2\textwidth]{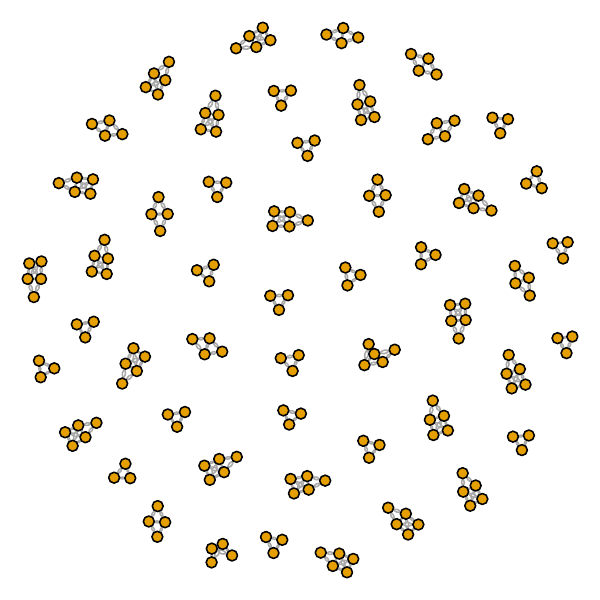} &
     \includegraphics[width=0.2\textwidth]{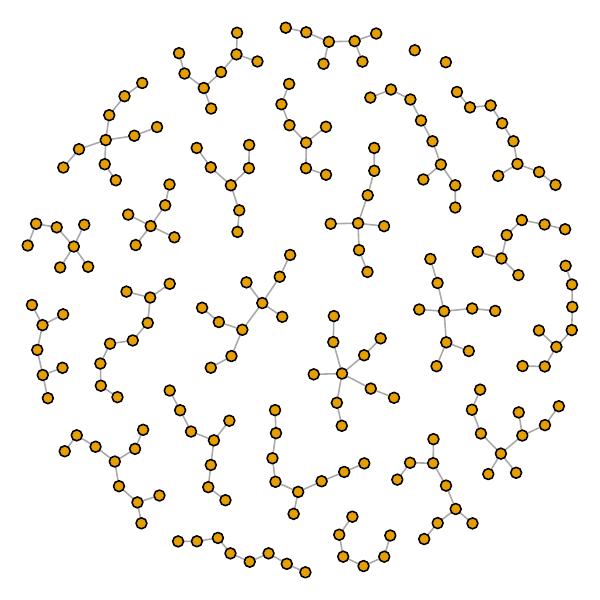}  \\
     (a) Graph of the Gaussian graphical model & (b) Graph of the Ising model
\end{tabular}
\caption{Examples of the adjacency graphs used for the Gaussian graphical model and the Ising model.}
\label{fig:sim}
\end{figure}

The Gaussian graphical model is generated as follows. First, we generate the adjacency matrix consisting of triangles, four-cycles and five-cycles: $3 m_1$ nodes form $m_1$  triangles, $4 m_2$ nodes form $m_2$ four-cycles, and  $5 m_3$ nodes form $m_3$ five-cycles, while these sub-graphs have no overlapping nodes with each other. Then dimension $d = 3m_1 + 4 m_2 +  5 m_3$. To make the problem more challenging, we randomly add one edge to the four-cycles and three edges to the five-cycles with an example of the adjacency matrix illustrated in Figure~\ref{fig:sim} (a). 
We then generate the entries of precision matrix $\Theta^*$ by sampling from ${\rm Uniform}(0.85,1)$ for edges in the adjacency graph and setting other entries as zero. Finally, we add a small value $v = 0.1$ together with the absolute value of the minimal eigenvalue of $\Theta^*$ to the diagonal elements of $\Theta^*$ to ensure its positive-definiteness. We take different combinations of $(m_1,m_2,m_3)$:  $(20,10,20)$, $(20,10,30)$, $(40,20,20)$ and $(60,30,10)$, which yield $d = 200, 250, 300$ and $350$, respectively. 


The Ising model is generated similarly where the adjacency graph  consisting of $m$-trees with $m$ sampled uniformly from $6$ to $10$ as illustrated by Figure ~\ref{fig:sim} (b). The dimension $d$ is also chosen from $\{200, 250, 300,350\}$ where the last few nodes (with size smaller than $m$) are set to be disconnected with all other nodes. We sample $w_e^{*}$ from ${\rm Uniform}(0.9,1)$ for all edges $e$'s in the adjacency graph and set $w_e^{*} = 0$ otherwise. The thresholding value $\theta$ in \eqref{eq:That:Ising} is set as $0.45$. 

The results for the Gaussian graphical model is presented in Table~\ref{tab1}. We can see that the FDRs of all configurations are well controlled by their nominal level $q$ while the power is pretty good for all configurations. Specifically, with the sample size $n$ increases, the FDRs tend to decrease. 
The FDRs for detection triangles are larger than those of four-cycles and five-cycles since there are more hypotheses to test for four cycles and five cycles compared to  triangles. 

The results for the Ising model is presented in Table~\ref{tab2}. The similar pattern is observed as the Gaussian graphical model. All FDRs are below the nominal level $q$ for different values of $n,d$ and $q$. The power increases quickly as the sample size $n$ increases. The results of the five-tree (a) and five-tree (b) are more similar compared to that of five-tree (c), as the numbers of hypotheses are same for testing five-tree (a) and five-tree (b).

\begin{table}[]
  \centering
  \caption{Empirical FDR and power averaged over $100$ repetitions in the detection of triangles, four-cycles and five-cycles under the Gaussian graphical model. }
      \vspace{10pt}
    {
    \resizebox{.8\columnwidth}{!}{
    \begin{tabular}{l|rrrrrrrrrrrr}
\hline 
  $q$        & \multicolumn{3}{c|}{ $0.05$ } & \multicolumn{3}{c|}{ $0.1$ } & \multicolumn{3}{c|}{ $0.05$ } & \multicolumn{3}{c}{ $0.1$ } \\ \hline 
   $n$     & \multicolumn{1}{r}{300} & \multicolumn{1}{r}{350} & \multicolumn{1}{r|}{400} & \multicolumn{1}{r}{300} & \multicolumn{1}{r}{350} & \multicolumn{1}{r|}{400} & \multicolumn{1}{r}{300} & \multicolumn{1}{r}{350} & \multicolumn{1}{r|}{400} & \multicolumn{1}{r}{300} & \multicolumn{1}{r}{350} & 400 \\ \hline 
          & \multicolumn{6}{c|}{FDR}                      & \multicolumn{6}{c}{Power} \\ \hline 
    $d$     & \multicolumn{12}{c}{(i) Triangle} \\ 
    200   & 0.031 & 0.025 & 0.020 & 0.040 & 0.033 & 0.025 & 0.999 & 1.000 & 1.000 & 0.999 & 1.000 & 1.000 \\
    250   & 0.036 & 0.034 & 0.033 & 0.046 & 0.045 & 0.042 & 0.999 & 1.000 & 1.000 & 1.000 & 1.000 & 1.000 \\
    300    & 0.029 & 0.030 & 0.030 & 0.037 & 0.037 & 0.039 & 0.999 & 1.000 & 1.000 & 1.000 & 1.000 & 1.000 \\
    400    & 0.029 & 0.025 & 0.026 & 0.036 & 0.031 & 0.032 & 1.000 & 1.000 & 1.000 & 1.000 & 1.000 & 1.000 \\ \hline 
    $d$     & \multicolumn{12}{c}{(ii) Four-cycle} \\ 
    200   & 0.011 & 0.008 & 0.005 & 0.014 & 0.011 & 0.007 & 0.996 & 1.000 & 1.000 & 0.997 & 1.000 & 1.000 \\
250    & 0.012 & 0.013 & 0.012 & 0.016 & 0.016 & 0.016 & 0.992 & 0.999 & 1.000 & 0.994 & 1.000 & 1.000 \\
    300    & 0.011 & 0.009 & 0.009 & 0.013 & 0.012 & 0.012 & 0.992 & 0.998 & 1.000 & 0.994 & 0.999 & 1.000 \\
    400   & 0.011 & 0.010 & 0.011 & 0.014 & 0.013 & 0.014 & 0.990 & 0.998 & 1.000 & 0.992 & 0.999 & 1.000 \\ \hline 
   $d$     & \multicolumn{12}{c}{(iii) Five-cycle} \\
    200   & 0.001 & 0.002 & 0.001 & 0.001 & 0.002 & 0.001 & 0.961 & 0.988 & 0.997 & 0.964 & 0.989 & 0.997 \\
    250   & 0.003 & 0.004 & 0.003 & 0.004 & 0.005 & 0.003 & 0.958 & 0.987 & 0.994 & 0.965 & 0.988 & 0.994 \\
    300  & 0.002 & 0.003 & 0.002 & 0.003 & 0.003 & 0.002 & 0.941 & 0.981 & 0.992 & 0.947 & 0.982 & 0.993 \\
    400  & 0.001 & 0.002 & 0.002 & 0.002 & 0.002 & 0.003 & 0.914 & 0.976 & 0.993 & 0.925 & 0.979 & 0.994 \\
    \bottomrule
    \end{tabular}%
    }
    }
  \label{tab1}%
\end{table}%

\begin{table}[htbp]
  \centering
  \caption{Empirical FDR and power averaged over $100$ repetitions in the detection of triangles, four-cycles and five-cycles under the Ising model.}
    \vspace{10pt}
    {
    \resizebox{.8\columnwidth}{!}{
    \begin{tabular}{l|rrrrrrrrrrrr}
   \hline 
  $q$        & \multicolumn{3}{c|}{ $0.05$ } & \multicolumn{3}{c|}{ $0.1$ } & \multicolumn{3}{c|}{ $0.05$ } & \multicolumn{3}{c}{ $0.1$ } \\ \hline 
   $n$     & \multicolumn{1}{r}{300} & \multicolumn{1}{r}{350} & \multicolumn{1}{r|}{400} & \multicolumn{1}{r}{300} & \multicolumn{1}{r}{350} & \multicolumn{1}{r|}{400} & \multicolumn{1}{r}{300} & \multicolumn{1}{r}{350} & \multicolumn{1}{r|}{400} & \multicolumn{1}{r}{300} & \multicolumn{1}{r}{350} & 400 \\ \hline  
          & \multicolumn{6}{c|}{FDP}                      & \multicolumn{6}{c}{Power} \\ 
        \hline 
    $d$     & \multicolumn{12}{c}{ Five-tree (a)} \\
    200    & 0.014 & 0.020 & 0.032 & 0.018 & 0.024 & 0.040 & 0.531 & 0.770 & 0.888 & 0.575 & 0.800 & 0.902 \\
    250    & 0.011 & 0.013 & 0.025 & 0.013 & 0.018 & 0.030 & 0.461 & 0.744 & 0.870 & 0.505 & 0.774 & 0.886 \\
    300    & 0.009 & 0.012 & 0.016 & 0.011 & 0.014 & 0.020 & 0.452 & 0.678 & 0.838 & 0.492 & 0.710 & 0.857 \\
    400  & 0.005 & 0.007 & 0.012 & 0.006 & 0.009 & 0.014 & 0.379 & 0.627 & 0.834 & 0.416 & 0.661 & 0.851 \\ \hline
    $d$     & \multicolumn{12}{c}{ Five-tree (b)} \\
    200    & 0.017 & 0.024 & 0.037 & 0.021 & 0.030 & 0.046 & 0.524 & 0.775 & 0.889 & 0.571 & 0.802 & 0.903 \\
    250    & 0.012 & 0.017 & 0.027 & 0.014 & 0.022 & 0.032 & 0.479 & 0.748 & 0.881 & 0.515 & 0.777 & 0.899 \\
    300    & 0.010 & 0.011 & 0.021 & 0.012 & 0.013 & 0.026 & 0.464 & 0.687 & 0.841 & 0.505 & 0.719 & 0.861 \\
    400    & 0.006 & 0.009 & 0.012 & 0.007 & 0.011 & 0.017 & 0.398 & 0.631 & 0.837 & 0.431 & 0.671 & 0.857 \\ \hline
    $d$     & \multicolumn{12}{c}{ Five-tree (c)} \\ 
    200    & 0.016 & 0.030 & 0.026 & 0.017 & 0.033 & 0.032 & 0.420 & 0.579 & 0.732 & 0.450 & 0.606 & 0.754 \\
    250   & 0.013 & 0.022 & 0.026 & 0.015 & 0.025 & 0.032 & 0.402 & 0.611 & 0.764 & 0.437 & 0.640 & 0.778 \\
    300   & 0.008 & 0.011 & 0.019 & 0.010 & 0.011 & 0.028 & 0.384 & 0.568 & 0.729 & 0.428 & 0.593 & 0.747 \\
    400   & 0.007 & 0.008 & 0.010 & 0.008 & 0.011 & 0.018 & 0.358 & 0.524 & 0.747 & 0.385 & 0.551 & 0.776 \\ \hline 
    \end{tabular}%
    }
    }
  \label{tab2}%
\end{table}%


We then consider the testing for persistent homology. Since the true correlation graph in Ising model does not have a explicit form, we only consider the Gaussian graphical model for the simulation in this part. 

\begin{figure}[htbp]
    \centering
\includegraphics[scale=0.15]{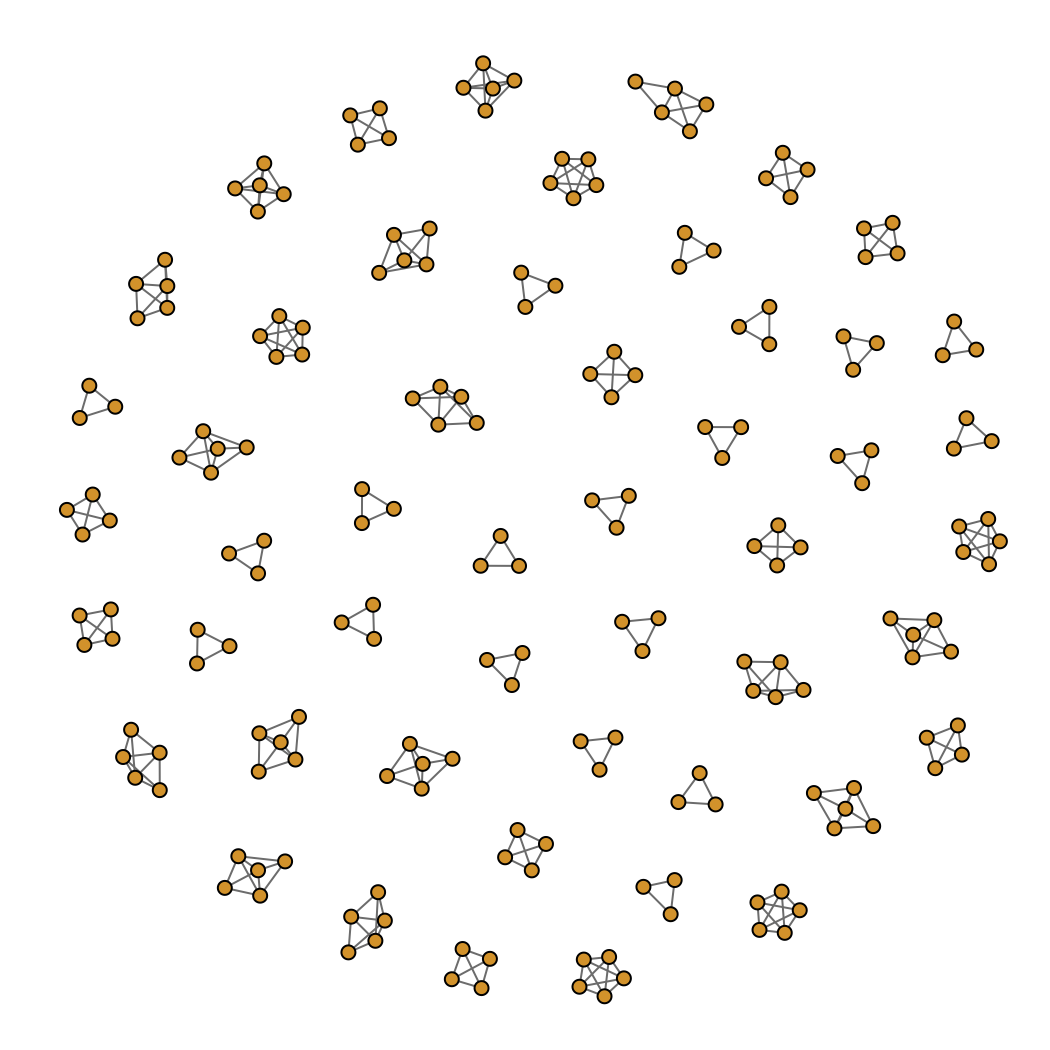}
    \caption{Graph of the Gaussian graphical model for persistent homology}
    \label{HomologyGraph}
\end{figure}

The model is generated in a similar way to the above setting. First, we generate the adjacency matrix consisting of triangles, four-cliques, and five-cliques: $3m_1$ nodes form $m_1$ triangles, $4m_2$ nodes form $m_2$ four-cliques, and $5m_3$ nodes form $m_3$ five-cliques, while these sub-graphs have no overlapping nodes with each other.  We then randomly reduce the edges of five cliques with probability $p = 0.1$. An example of the adjacency matrix is illustrated in Figure \ref{HomologyGraph}. The sample size $n$ is $400$ and the dimension $d$ is $3m_1 + 4m_2 + 5m_3$. We then generate the entries of the precision matrix $\Theta^*$ by sampling from Uniform$(0,10)$ for edges in the adjacency graph and setting other entries as zero. Finally, we add a small value $v = 0.25$ together with the absolute value of the minimal eigenvalue of $\Theta^*$ to the diagonal elements of $\Theta^*$ to ensure its positive definiteness. We set ($m_1$, $m_2$, $m_3$): $(10,30,10)$, which yields $d = 200$ respectively. The filtration level is set between $0$ and $1$.

\begin{figure}[htbp]
\centering
\subfigure
{
\includegraphics[scale=0.3]{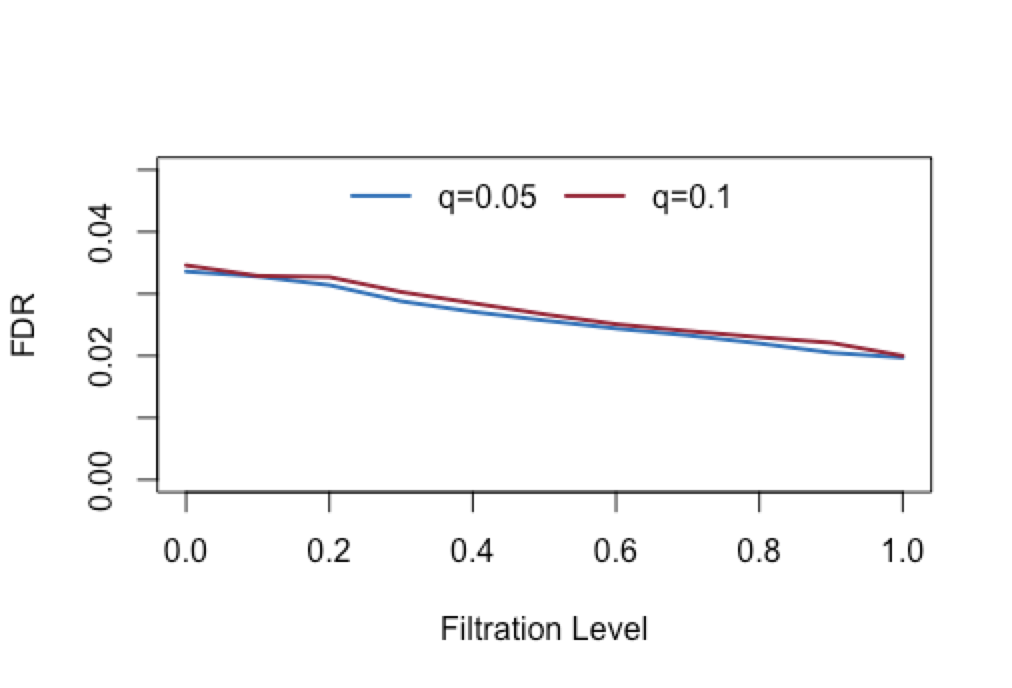} \label{1}
}
\quad
\subfigure
{
\includegraphics[scale=0.3]{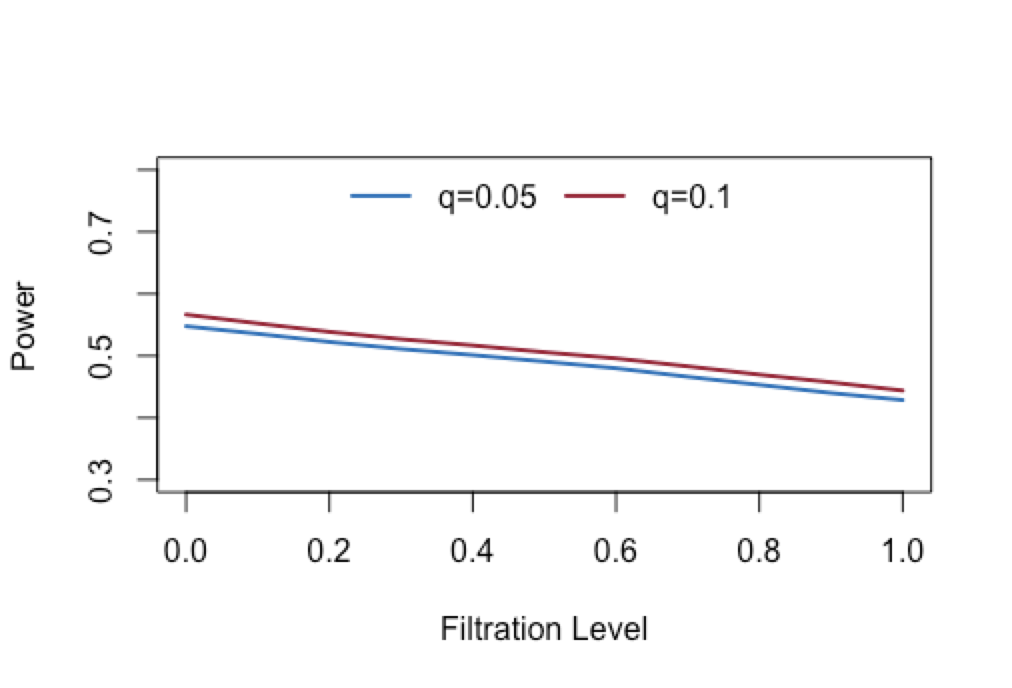} \label{2} 
}
\caption{Empirical FDR and power averaged over $100$ repetitions in the detection of homology under the Gaussian graphical model.}

\label{HomologyResult}
\end{figure}

The results for the Gaussian graphical model are presented in Figure \ref{HomologyResult}. We can see that the FDR is well controlled by their nominal level $q$ while the power is above $0.4$. Specifically, as the filter level increases, the FDR tends to decrease. 
 
\end{appendix}

\end{document}